\documentclass[12pt,a4paper]{article}
\usepackage[T2A]{fontenc}
\usepackage[cp1251]{inputenc}
\usepackage[russian]{babel}
\usepackage{amsmath,amssymb,amsthm,amsfonts,amscd}

\begin{document}
\sloppy
\date{}

\title{Геометрический подход к стабильным гомотопическим группам сфер II. Инвариант Кервера.}

\author{П.М.Ахметьев \thanks{Работа автора поддержена грантом the London Royal Society
(1998-2000), РФФИ 08-01-00663,  INTAS 05-1000008-7805.}}

\sloppy \theoremstyle{plain}
\newtheorem{theorem}{Теорема}
\newtheorem*{main*}{Основная Теорема}
\newtheorem*{theorem*}{Теорема}
\newtheorem{lemma}[theorem]{Лемма}
\newtheorem{proposition}[theorem]{Предложение}
\newtheorem{corollary}[theorem]{Следствие}
\newtheorem{conjecture}[theorem]{Гипотеза}
\newtheorem{problem}[theorem]{Проблема}

\theoremstyle{definition}
\newtheorem{definition}[theorem]{Определение}
\newtheorem{remark}[theorem]{Замечание}
\newtheorem*{remark*}{Замечание}
\newtheorem*{example*}{Пример}
\newtheorem{example}[theorem]{Пример}
\def\aa{\dot{a}}
\def\i{{\bf i}}
\def\j{{\bf j}}
\def\k{{\bf k}}
\def\hh{{\bf \dot{h}}}
\def\e{{\bf e}}
\def\f{{\bf f}}
\def\h{{\bf h}}
\def\dd{{\dot{d}}}
\def\Z{{\Bbb Z}}
\def\R{{\Bbb R}}
\def\RP{{\Bbb R}\!{\rm P}}
\def\N{{\Bbb N}}
\def\C{{\Bbb C}}
\def\A{{\bf A}}
\def\D{{\bf D}}
\def\Q{{\bf Q}}
\def\QQ{{\dot{\bf Q}}}
\def\E{{\bf E}}
\def\F{{\bf F}}
\def\J{{\bf J}}
\def\G{{\bf G}}
\def\I{{\bf I}}
\def\II{{\dot{\bf I}}}
\def\H{{\bf H}}
\def\fr{{\operatorname{fr}}}
\def\st{{\operatorname{st}}}
\def\mod{{\operatorname{mod}\,}}
\def\cyl{{\operatorname{cyl}}}
\def\dist{{\operatorname{dist}}}
\def\sf{{\operatorname{sf}}}
\def\dim{{\operatorname{dim}}}
\def\dist{\operatorname{dist}}

\maketitle

\begin{abstract}
Представлено решение Проблемы Кервера. Вводится понятие абелевой
структуры скошенно--оснащенного погружения, бициклической
структуры $\Z/2^{[3]}$--оснащенного погружения и бикватернионной
структуры $\Z/2^{[5]}$--оснащенного погружения.
 Используя введенные понятия  мы доказываем, что при достаточно большом   $n$,
$n=2^l-2$, произвольное скошенно--оснащенное погружение имеет
инвариант нулевой Кервера. В доказательстве использована теорема о
ретракции скошенно-оснащенного погружения в заданном классе
нормального кобордизма. Доказательство теоремы о ретракции  также
приводится в работе.
\end{abstract}

\section{Самопересечение погружений и Инвариант Кервера}

Проблема Инвариантов Кервера 1 является открытой проблемой в
Алгебраической топологии, см. по поводу алгебраического подхода к
решению [S], [B-J-M], [C-J-M]. Мы рассмотрим геометрический подход
к решению Проблемы, который основан на результатах П.Дж.Экклза
[E1]. По поводу другого геометрического подхода см. [C1],[C2].

Рассмотрим гладкое погружение  $f: M^{n-1} \looparrowright \R^n$,
$n= 2^l -2$, $l>1$ общего положения коразмерности 1. Обозначим
через  $g: N^{n-2} \looparrowright \R^n$ погружение многообразия
самопересечения.

\begin{definition}
 Инвариант Кервера погружения  $f$ определяется по формуле
\begin{eqnarray}\label{arf}
 \Theta_{sf}(f) = <\eta_N^{\frac{n-2}{2}}; [N^{n-2}] >,
\end{eqnarray}
 где
через $\eta_N = w_2(N^{n-2})$ обозначен двумерный нормальный класс
Штифеля-Уитни многообразия $N^{n-2}$.
\end{definition}

Инвариант Кервера является инвариантом класса регулярного
кобордизма погружения   $f$. Более того, инвариант Кервера
является гомоморфизмом
\begin{eqnarray}\label{1}
\Theta_{sf}: Imm^{sf}(n-1,1) \to  \Z/2.
\end{eqnarray}
 Нормальное расслоение  $\nu_g$ погружения  $g:
N^{n-2} \looparrowright \R^n$ является  2-мерным расслоением над
 $N^{n-2}$, которое снабжено $\D_4$--оснащением. Классифицирующее
 отображение этого расслоения (как и соответствующий
 характеристический класс) обозначается через
 $\eta_N: N^{n-2} \to K(\D_4,1)$.  Пара
$(g,\eta_N)$ представляет элемент в группе кобордизма
$Imm^{\D_4}(n-2,2)$. Гомоморфизм
\begin{eqnarray}\label{2}
\delta_{\D_4}: Imm^{sf}(n-1,1) \to Imm^{\D_4}(n-2,2)
\end{eqnarray}
 корректно определен.

Группа кобордизма $Imm^{sf}(n-k,k)$ обобщает группу кобордизма
$Imm^{sf}(n-1,1)$. Новая группа определена как группа кобордизма
троек   $(f,\Xi,\kappa_M)$, где $f: M^{n-k} \looparrowright \R^n$
погружение, причем задан изоморфизм  $\Xi: \nu(g) = k \kappa_M$,
который называется скошенным оснащением, через $\nu(f)$ обозначено
нормальное расслоение погружения $f$, $\kappa_M$ является заданным
линейным расслоением над  $M^{m-k}$, характеристический класс
которого обозначается также через $\kappa_M \in
H^1(M^{m-k};\Z/2)$. Отношение кобордизма на пространстве
рассматриваемых троек является стандартным.

Группа $Imm^{\D_4}(n-2,2)$ обобщается следующим образом. Определим
группы кобордизмов $Imm^{\D_4}(n-2k,2k)$. Каждый элемент группы
$Imm^{\D_4}(n-2k,2k)$ представлен тройкой $(g,\Psi,\eta_N)$, где
$g: N^{n-2k} \looparrowright \R^n$ -- погружение,
$\Psi$--диэдральное оснащение в коразмерности $2k$, т.е.
фиксированный изоморфизм $\Xi: \nu_g = k \eta_N$, где $\eta_N$
является 2-мерным расслоением над $N^{n-2k}$. Характеристическое
отображение этого расслоения, а также соответствующий
характеристический класс (соответствующий универсальный
характеристический класс) обозначается также через $\eta_N:
N^{n-2k} \to K(\D_4,1)$, $\eta_N \in H^2(N^{n-2k};\Z/2)$ ($\tau
\in H^2(K(\D_4,1);\Z/2)$). Отображение $\eta_N$ оказывается
характеристическим и для расслоения $\nu_g$, поскольку $\nu_g=k
\eta_N$.

По определению гомоморфизм Кервера ($\ref{1}$) задан композицией
гомоморфизма ($\ref{2}$) и гомоморфизма
\begin{eqnarray}\label{3}
\Theta_{\D_4} : Imm^{\D_4}(n-2,2) \to \Z/2.
\end{eqnarray}
Гомоморфизм ($\ref{3}$) называется инвариантом Кервера
$\D_4$--оснащенного погружения.

Гомоморфизм Кервера определен в более общей ситуации при помощи
прямого обобщения гомоморфизмов   ($\ref{1}$) и ($\ref{3}$):
\begin{eqnarray}\label{4a}
\Theta_{sf}^{k}: Imm^{sf}(n-k,k) \to \Z/2, \qquad \Theta_{sf}^k: =
\delta_{\D_4}^k \circ \Theta_{\D_4}^k.
\end{eqnarray}
\begin{eqnarray}\label{44}
\Theta_{\D_4}^{k} : Imm^{\D_4}(n-2k,2k) \to \Z/2, \qquad
\Theta_{\D_4}^{k}[(g,\Psi,\eta_N)] = \langle
\eta_N^{\frac{n-2k}{2}}; [N^{n-2k}] \rangle.
\end{eqnarray}
(Для $k=1$ новый гомоморфизм ($\ref{44}$) совпадает с
вышеопределенным гомоморфизмом ($\ref{3}$) при этом следующая
диаграмма коммутативна:
\begin{eqnarray}\label{5}
\begin{array}{ccccc}
Imm^{sf}(n-1,1) & \stackrel  {\delta_{\D_4}}{\longrightarrow} &
Imm^{\D_4}(n-2,2) & \stackrel{\Theta_{\D_4}}{\longrightarrow} & \Z/2  \\
\downarrow J_{sf}^k & &  \downarrow J_{\D_4}^{k}  &  &  \vert \vert \\
Imm^{sf}(n-k,k) & \stackrel{\delta_{\D_4}^{k}}{\longrightarrow} &
Imm^{\D_4}(n-2k,2k) &
\stackrel{\Theta_{\D_4}^k}{\longrightarrow} & \Z/2.  \\
\end{array}
\end{eqnarray}

Нам потребуется обобщить формулу ($\ref{44}$) для погружений с
оснащением более общего вида. Обозначим через $\Z/2^{[d]}$
сплетение $2^{d-1}$ экземпляров элементарной циклической группы.
Эта группа является подгруппой в группе  $O(2^{d-1})$, которая
определяется следующим условием:

-- Преобразования из $\Z/2^{[d]}$ оставляют инвариантными
следующие $d-1$ наборов $\Omega_{d}$, $\Omega_{d-1}$, $\dots$,
$\Omega_2$ координатных подпространств. Набор подпространств
$\Omega_i$, $2 \le i \le d$ состоит из $2^{i-1}$ координатных
подпространств, порожденных базисными векторами  $((\e_1, \dots
\e_{2^{d-i+1}}), \dots,  (\e_{2^d-2^{d-i+1}+1}, \dots,
\e_{2^d}))$.

В частности, $\D_4 = \Z/2^{[2]}$.
\[  \]

Рассмотрим погружение $g: N^{n-k2^{d-1}} \looparrowright \R^n$
общего положения коразмерности  $k2^{d-1}$. Скажем, что погружение
$g$ является $\Z/2^{[d]}$--оснащенным (c кратностью $k$), если
определен изоморфизм $\Psi: \nu_g = k \eta_N$ нормального
расслоение $\nu_g$ погружения $g$ c суммой Уитни $k$ экземпляров
$2^{d-1}$--мерного расслоения $\eta_N$ со структурной группой
$\Z/2^{[d]}$.

Расслоение $\eta_N$ классифицируется отображением $\eta_N:
N^{n-k2^{d-1}} \to K(\Z/2^{[d]},1)$. (Cоответствующий
характеристический класс) обозначим также через $\eta_N$).
Характеристический класс универсального $2^{d-1}$--мерного
$\Z/2^{[d]}$--расслоения над $K(\Z/2^{[d]},1)$ обозначим через
$\tau_{[d]}$. Таким образом, $\eta_N^{\ast}(\tau_{[d]})=\eta_N$.
Отображение $\eta_N$ оказывается характеристическим и для
расслоения $\nu_g$, поскольку $\nu_g=k \eta_N$.

Всевозможные тройки   $(g,\Psi,\eta_N)$, которые были описаны
выше, порождают группу кобордизма
$Imm^{\Z/2^{[d]}}(n-k2^{d-1},k2^{d-1})$.
 В некоторых рассуждениях
при обозначениях будет использован дополнительный индекс,
связанный  со структурной группой. Например, представитель группы
$Imm^{\D_4}(n-2k,2k)$ будем иногда обозначать через
   $(g_{[2]},\Psi_{[2]}, \eta_{N_{[2]}})$ и т.д.

Многообразие самопересечения произвольного
$\Z/2^{[d]}$--оснащенного  погружения является
$\Z/2^{[d+1]}$--оснащенным  погружением. Таким образом,
многообразие самопересечения представляет тройку $(h, \Lambda,
\zeta_L)$, где $h: L^{n-k2^{d}} \looparrowright \R^n$--погружение,
$\Lambda: k \zeta_L = \nu_h$, $\zeta_L: L^{n-k2^{d}} \to
K(\Z/2^{[d+1]},1)$--классифицирующее отображение $2^{d}$--мерного
расслоения $\zeta_L$.
 Корректно
определен гомоморфизм
\begin{eqnarray}\label{6}
\delta_{\Z/2^{[d+1]}}^{k} : Imm^{\Z/2^{[d]}}(n-k2^{d-1},k2^{d-1})
\to
 Imm^{\Z/2^{[d+1]}}(n-k2^{d},k2^{d}),
\end{eqnarray}
ставящий в соответствие классу нормального кобордизма
$[(g,\Psi,\eta_N)]$ класс нормального кобордизма $[(h, \Lambda,
\zeta_L)]$.

Определена подгруппа $i_{[d+1]}: \Z/2^{[d]} \subset \Z/2^{[d+1]}$,
как подгруппа преобразований подпространства $\R^{2^{d-1}} \subset
\R^{2^{d}}$, порожденного первыми $2^{d-1}$ базисными векторами.

Определена подгруппа
\begin{eqnarray}\label{7}
\bar i_{[d+1]}: \Z/2^{[d]} \oplus \Z/2^{[d]} \subset \Z/2^{[d+1]}
\end{eqnarray}
индекса 2, как подгруппа преобразований пространства  которая
состоит из преобразований, оставляющих инвариантным каждое
подпространство из набора $\Omega_2$.

Подгруппа  ($\ref{7}$) индуцирует двулистное накрытие $\pi_{[d+1]}
: K(\Z/2^{[d]} \oplus \Z/2^{[d]},1) \to K(\Z/2^{[d+1]},1)$.
Характеристическое отображение $\zeta_L: L^{n-k2^{d}} \to
K(\Z/2^{[d+1]},1)$ индуцирует двулистное накрытие $\pi_{[d+1],L} :
\bar L^{n-k2^{d}} \to L^{n-k2^{d}}$ из накрытия $\pi_{[d+1]}$ над
классифицирующим пространством. Двулистное накрытие
$\pi_{[d+1],L}$ можно определить геометрически, оно совпадает с
каноническим двулистном накрытием над многообразием $L^{n-k2^{d}}$
точек самопересечения $\Z/2^{[d]}$--оснащенного погружения
$(g,\Psi,\eta_N)$ (см. [A], раздел 1).

Определена проекция $p_{[d]}: \Z/2^{[d]} \oplus \Z/2^{[d]} \to
\Z/2^{[d]}$ на первое слагаемое, которая индуцирует
 отображение $p_{[d]}: K(\Z/2^{[d]} \oplus \Z/2^{[d]},1) \to
 K(\Z/2^{[d]},1)$.

Для многообразия самопересечения $(h,\Lambda,\zeta_L)$
произвольного $\Z/2^{[d]}$-оснащенного погружения $(g, \Psi,
\eta_N)$ рассмотрим двулистное накрытие $\bar \zeta_L: \bar
L^{n-k2^{d}}_{[d]} \to K(\Z/2^{[d]} \oplus \Z/2^{[d]},1)$, над
классифицирующим отображением $\zeta_L$, которое индуцировано
накрытием $\pi_{[d+1],L}$. Это накрытие совпадает с  каноническим
2-листным накрытием над классифицирующим отображением $\zeta_L:
L^{n-k2^d} \to K(\Z/2^{[d+1]},1)$, которое определено из
геометрических соображений. Характеристический класс $(p_{[d]}
\circ \bar \zeta)^{\ast}(\tau_{[d]}) \in H^{2^{d-1}}(\bar
L^{n-k2^{d}}_{[d]};\Z/2)$,
 $\tau_{[d]} \in H^{2^{d-1}}(K(\Z/2^{[d]},1)$
 обозначим через $\bar \zeta_{[d],L}$.

Определено отображение
 $i_{tot}=i_{[3]}\circ \dots \circ
i_{[d]}$:
\begin{eqnarray}\label{pp}
 K(\D_4,1) \stackrel{i_{[3]}}{\longrightarrow} K(\Z/2^{[3]},1)
\stackrel{i_{[4]}}{\longrightarrow}  \dots
\stackrel{i_{[d]}}{\longrightarrow} K(\Z/2^{[d]},1)
\stackrel{i_{[d+1]}}{\longrightarrow} K(\Z/2^{[d+1]},1).
\end{eqnarray}
 Определена башня 2-листных канонических накрытий
\begin{eqnarray}\label{pipi}
\bar L^{n-k2^{d}}_{[2]} \stackrel{\pi_{[3]}}{\longrightarrow} \bar
L^{n-k2^{d}}_{[3]} \stackrel{\pi_{[4]}}{\longrightarrow} \dots
\stackrel{\pi_{[d]}}{\longrightarrow} \bar L^{n-k2^{d}}_{[d]}
\stackrel{\pi_{[d+1]}}{\longrightarrow} L^{n-k2^d}.
\end{eqnarray}
Эта башня накрытий снабжена характеристическим отображением в
диаграмму ($\ref{pp}$) классифицирующих пространств. Обозначим
через
\begin{eqnarray}\label{tot}
\pi_{tot}=\pi_{[3]}\circ \dots \circ \pi_{[d]} : K(\D_4,1) \to
K(\Z/2^{[d]},1)
\end{eqnarray}
башню накрытий, индуцированную композицией гомоморфизмов в
диаграмме $(10)$.
 Определена
последовательность характеристических классов
$$
\bar \zeta_{[2],L} \in H^2(\bar L^{n-k2^{d}}_{[2]};\Z/2), \dots,
\bar \zeta_{[d],L} \in H^{2^{d-1}}(\bar
L^{n-k2^{d}}_{[d]};\Z/2),$$
\begin{eqnarray}\label{bzbz}
\zeta_{[d+1],L} \in H^{2^{d}}( L^{n-k2^{d}};\Z/2).
\end{eqnarray}
Каждый элемент в этой последовательности индуцирован из
характеристического класса соответствующего универсального
пространства в ($\ref{pp}$). Башня накрытий $(11)$ и
последовательность характеристических классов $(13)$ определена не
только для $\Z/2^{[d+1]}$--оснащенного многообразия, которое
является многообразием самопересечения некоторого
$\Z/2^{[d]}$--оснащенного погружения, но и для произвольного
$\Z/2^{[d+1]}$--оснащенного многообразия.

\begin{definition}
Инвариант Кервера  $\Theta_{\Z/2^{[d+1]}}^k$ произвольного
$\Z/2^{[d+1]}$--оснащенного  погружения $(h,\Lambda,\zeta_L)$
определим следующей формулой:
\begin{eqnarray}\label{thetaz2^3}
 \Theta_{\Z/2^{[d+1]}}^k  (h,\Lambda,\zeta_L) = \langle \bar \eta_{[2],L}^{\frac{n-k2^{d}}{2}};[\bar
N_{[2]}] \rangle,
\end{eqnarray}
где через $[\bar N_{[2]}]$ обозначен фундаментальный класс
накрывающего многообразия  в накрытии ($\ref{tot}$).
\end{definition}

Построенный инвариант определяет гомоморфизм
$\Theta_{\Z/2^{[d]}}^k: Imm^{\Z/2^{[d]}}(n,n-k2^{d-1}) \to \Z/2$,
который включен в следующую коммутативную диаграмму:

\begin{eqnarray}\label{77}
\begin{array}{ccc}
Imm^{\Z/2^{[d]}}(n-k2^{d-1},k2^{d-1}) & \stackrel{\Theta_{\Z/2^{[d]}}^k}{\longrightarrow} & \Z/2  \\
\downarrow             \delta_{\Z/2^{[d+1]}}^k  &  &  \| \\
Imm^{\Z/2^{[d+1]}}(n-k2^{d},k2^{d}) &
\stackrel{\Theta_{\Z/2^{[d+1]}}^{k}}{\longrightarrow} & \Z/2.  \\
\end{array}
\end{eqnarray}

 В разделе 1
определено понятие  $\J_b$--структуры (абелева структура)
скошенно-оснащенного погружения, представляющего элемент из группы
$Imm^{sf}(n-k,k)$. Доказана Теорема 6 о том, что при
соответствующих размерностных ограничениях и по модулю элементов
нечетного порядка произвольный класс кобордизма
скошенно-оснащенного погружения допускает $\J_b$--структуру.
Условие этой теоремы предполагает существование ретракции
характеристического класса скошенно-оснащенного многообразия, см.
Определение 5. Теорема о ретракции  доказана в разделе 8.

В разделе 3 сформулировано понятие $\I_a \oplus \II_a$--структуры
(бициклическая структура)   $\Z/2^{[3]}$--оснащенного погружения.
В Следствии 13 Теоремы 12 доказывается, что в условиях Теоремы 6
(где, в частности, определено натруральное число $n_s$)
произвольный элемент из группы

$$Im(\delta_{\Z/2^{[3]}}^{\frac{n-n_s}{32}}:
Imm^{sf}(n-\frac{n-n_s}{32},\frac{n-n_s}{32}) \to $$
\begin{eqnarray}\label{ImD4}
Imm^{\D_4}(n-\frac{n-n_s}{8},\frac{n-n_s}{8}))
\end{eqnarray}
представлен $\Z/2^{[3]}$--оснащенным погружением  с бициклической
структурой. Для такого погружения инвариант Кервера выражается
через $\I_a \oplus \II_a$--характеристический класс многообразия
самопересечения.

В разделе 4 сформулировано понятие $ \Q_a \oplus \QQ_a$--структуры
(бикватернионная структура)  для $\Z/2^{[5]}$--оснащенного
погружения. В Следствии 19 Теоремы 18 доказывается, что в условиях
Теоремы 7 произвольный элемент из группы
$$Im(\delta_{\Z/2^{[6]}}
^{\frac{n-n_s}{32}}: Imm^{sf}(n-\frac{n-n_s}{32},\frac{n-n_s}{32})
\to$$
\begin{eqnarray}\label{ImZZZ}
 Imm^{\Z/2^{[6]}}( n-\frac{n-n_s}{2},\frac{n-n_s}{2}))
\end{eqnarray}
представлен $\Z/2^{[6]}$--оснащенным погружением  с
бикватернионной структурой. Для такого погружения инвариант
Кервера выражается через $\Q_a \oplus \QQ_a$--характеристический
класс многообразия самопересечения.

Автор благодарит Проф. M.Mahowald'a (2005) и Проф. R.Cohen'a
(2007) за обсуждения, Проф. А.А.Воронова за приглашение с докладом
в Университет Миннесоты  (2005). Проф. В. Чернова за приглашение с
докладом в Дартмуский Колледж  (2009).

Работа была начата на семинаре М.М.Постникова в 1998 году. Работа
посвящается памяти Проф. Ю.П.Соловьева. Теорема о ретракции
доказана на семинаре А.С.Мищенко.

\section{Геометрический контроль многообразия самопересечения скошенно-оснащенных погружений}
В этом и в последующих разделах через $Imm^{sf}(n-k,k)$,
$Imm^{\D_4}(n-2k,2k)$, и.т.д. будут обозначаться не сами группы
кобордизма, а их $2$--компоненты. Если первый аргумент в скобках,
обозначающий размерность погруженного многообразия, строго
положителен, то указанная группа является конечной $2$--группой.

Диэдральная группа  $\D_4$ определяется своим копредставлением
$\{a,b\vert a^4 = b^2 = e, [a,b]=a^2\}$. Эта группа является
подгруппой в группе движений плоскости  $O(2)$, т.е. группой
преобразований стандартной плоскости, сохраняющей пару
направляющих прямых, порожденных базисными векторами
$\{\f_1,\f_2\}$. Элемент $a$ представлен вращением плоскости на
угол $\frac{\pi}{2}$. Элемент $b$ соответствует отражению
плоскости относительно прямой с образующим вектором  $\f_1 +
\f_2$.

Рассмотрим подгруппу $\J_b \subset \D_4$ диэдральной группы,
порожденную элементами  $\{a^2,b\}$. Заметим, что это элементарная
$2$--группа ранга 2. Это -- группа движений, сохраняющих по
отдельности каждую из прямых $l_1$, $l_2$ с направляющими
векторами $\f_1+\f_2$, $\f_1-\f_2$ соответственно. Группа
когомологий $H^1(K(\J_b,1);\Z/2)$ также является элементарной
$2$--группой с двумя образующими. Определим образующую $\kappa_d$
($\kappa_c$) группы $H^1(K(\J_b,1);\Z/2)$, отвечающую за
преобразования прямой $l_1$ (прямой $l_2$) соответственно при
симметриях плоскости. Когомологический класс $\kappa_c \kappa_d$
является элементом из группы $H^2(K(\J_b,1)\Z/2))$, который
обозначается через $\tau_b$.

Рассмотрим отображение $i_{b}: K(\J_b,1) \to K(\D_4,1)$ и
рассмотрим обратный образ $i_b^{\ast}(\tau_{[2]})$ эйлерового
класса $\tau_{[2]} \in H^2(K(\D_4,1);\Z/2)$ универсального
2-расслоения. Нетрудно проверить равенство
\begin{eqnarray}\label{ibeta}
 i_b^{\ast}(\tau_{[2]}) = \tau_b.
\end{eqnarray}
 Определена расщепляющаяся проекция $p_{b,d}: \J_b \to \J_d$ на элементарную 2-группу
 $\J_d$.
Образующая образа представлена классом смежности элемента $a^2$
(центральная симметрия). Расщепляющее вложение обозначим через
$i_{d,b}: \J_d \subset \J_b$. Это вложение определено так, что
образующая группы $\J_d$ при этом вложении переходит в подгруппу
симметрий плоскости, порожденную симметрией прямой $l_1$.

\begin{definition}
Пусть скошенно-оснащенное погружение   $(f,\Xi,\kappa_M)$, $f:
M^{n-k} \looparrowright \R^n$ представляет элемент $x \in
Imm^{sf}(n-k,k)$. Пусть $\D_4$--оcнащенное погружение
$(g,\Psi,\eta_N)$, $g: N^{n-2k} \looparrowright \R^n$ представляет
элемент $y=\delta_{\D_4}^{k}(x) \in Imm^{\D_4}(n-2k,2k)$. Скажем,
что скошенно-оснащенное погружение $(f,\Xi,\kappa_M)$ является
$\J_b$--погружением (абелевым погружением), если структурное
отображение $\eta_N: N^{n-2k} \to K(\D_4,1)$ представлено в виде
композиции $\eta_{b,N}: N^{n-2k} \to K(\J_b,1)$ и отображения
$i_{b}: K(\J_b,1) \to K(\D_4,1)$.
\end{definition}

\begin{definition}

Пусть скошенно-оснащенное погружение   $(f,\Xi,\kappa)$, $f:
M^{n-k} \looparrowright \R^n$ представляет элемент $x \in
Imm^{sf}(n-k,k)$, причем $n > 32k$. Пусть $\D_4$--оcнащенное
погружение $(g,\Psi,\eta_N)$, $g: N^{n-2k} \looparrowright \R^n$
является погружением многообразия самопересечения погружения $f$ и
представляет элемент $y=\delta_{\D_4}^{k}(x) \in
Imm^{\D_4}(n-2k,2k)$. Скажем, что скошенно-оснащенное погружение
$(f,\Xi,\kappa_M)$ допускает абелеву структуру
($\J_b$--структуру), если существует отображение $\eta_b: N^{n-2k}
\to K(\J_b,1)$, удовлетворяющее следующему уравнению:
\begin{eqnarray}\label{etaa}
\Theta_{\D_4}^{k}(y) = \langle \eta_N^{15k}
\eta_{b,N}^{\frac{n-32k}{2}};[N] \rangle.
\end{eqnarray}
В этом уравнении характеристический класс
$\eta_{b,N}^{\ast}(\tau_b) \in H^2(N^{n-2k};\Z/2)$ снова обозначен
через $\eta_{b,N}$, $[N]$-- фундаментальный класс многообразия
$N^{n-2k}$, $\eta_N \in H^2(N^{n-2k};\Z/2)$ -- характеристический
класс $\D_4$--оснащения $\Psi$ многообразия $N^{n-2k}$,
характеристическое число $\Theta_{sf}^{k}$ определено по формуле
($\ref{4a}$).
\end{definition}

\subsubsection*{Пример}

Пусть скошенно-оснащенное погружение   $(f,\Xi,\kappa_M)$, $f:
M^{n-k} \looparrowright \R^n$ представляет элемент $x \in
Imm^{sf}(n-k,k)$, $n > 32k$ и является $\J_b$--погружением.
 Тогда скошенно-оснащенное погружение
$(f,\Xi,\kappa_M)$ допускает абелеву структуру.

\begin{definition}
Пусть $(f,\Xi,\kappa_M) \in Imm^{sf}(n-k,k)$,
 $f: M^{n-k} \looparrowright \R^n$, $\kappa_M \in H^1(M^{n-k};\Z/2)$
скошенного оснащения $\Xi$. Скажем, что пара $(M^{n-k},\kappa_M)$
допускает ретракцию порядка  $q$, если отображение $\kappa_M :
M^{n-k} \to \RP^{\infty}$ представлено следующей композицией
$\kappa = I \circ \kappa_M' : M^{n-k} \to \RP^{n-k-q-1} \subset
\RP^{\infty}$. Скажем, что элемент $[(f, \Xi, \kappa_M)]$
допускает ретракцию порядка $q$, если в этом классе кобордизма
 существует тройка  $(M'^{n-k},
\Xi', \kappa_M')$, допускающая ретракцию порядка  $q$.
\end{definition}

\begin{theorem}
Пусть $n_s$, $n>n_s$, натуральное число вида  $2^s-2$, $s \ge 6$.
Предположим, что элемент  $\alpha \in
Imm^{sf}(n-\frac{n-n_s}{32},\frac{n-n_s}{32})$ допускает ретракцию
порядка  $q=\frac{n_s}{2}+1$. Тогда элемент $\alpha$ допускает
$\J_b$--структуру.
\end{theorem}

Докажем следующую лемму.
\begin{lemma}
Для произвольного натурального $k'$, $k' \equiv 1 \pmod{2}$, $k'
\ge 7$, существует $PL$-отображение $d: \RP^{n-k'} \to \R^n$
общего положения, для которого  многообразие $N(d)$ с
особенностями с краем точек самопересечения отображения $d$
допускает отображение $\kappa_{N(d)}: N(d) \to K(\J_d,1)$,
ограничение которого на край $\partial N(d)$ (край этого
многообразия с особенностями состоит из критических точек
отображения $d$) совпадает с композицией $\partial N(d) \to
\RP^{n-k'} \subset \RP^{\infty} = K(\J_d,1)$.
\end{lemma}

\subsubsection*{Конструкция  отображения    $d: \RP^{n-k'} \to
\R^n$}

 Обозначим через $J_0$ стандартную $(n-k')$--мерную
сферу коразмерности $k'$, которая представлена как джойн
$\frac{n-k'+1}{2}=r$ копий окружности $S^1$. Обозначим стандартное
вложение $J_0$ в $\R^n$ через $i_{J_0}: J_0 \subset \R^n$.

Определено отображение $p': S^{n-k'} \to J_0$, которое получается
в результате взятия джойна $r$ копий стандартных 2-листных
накрытий $S^1 \to \RP^1$. Стандартное антиподальное действие $\Z/2
\times S^{n-k'} \to S^{n-k'}$,  коммутирует с отображением $p'$.
Тем самым, определено отображение  $p: \RP^{n-k'} \to J_0$.

Рассмотрим композицию $i_{J_0} \circ p: \RP^{n-k'} \to J_0 \subset
\R^n$. Отображение $d$ определено в результате малой
$\delta$-деформации общего положения этого отображения. Деформация
$i_{J_0} \circ p \mapsto d$ и ее калибр выбираются в процессе
доказательства.

\subsubsection*{Доказательство Леммы 7}
Доказательство аналогично доказательству Леммы 24 из [A].
\[  \]

\subsubsection*{Доказательство Теоремы 6}
Определим $k=\frac{n-n_s}{32}$. Пусть элемент $\alpha$ представлен
скошено-оснащенным погружением $(f,\Xi,\kappa_M)$, $f: M^{n-k}
\looparrowright \R^n$. По предположению существует отображение
$\kappa_M': M^{n-k} \to  \RP^{n-k-q-1}$, такое, что композиция
$M^{n-k} \to \RP^{n-k-q-1} \subset K(\J_d,1)$ совпадает с
отображением $\kappa: M^{n-k} \to K(\J_d,1)$, $q=\frac{n_s}{2}+1$.
Заметим, что $k+q+1$ нечетно, обозначим это число через $k'$ и
рассмотрим отображение $d: \RP^{n-k'} \to \R^n$, построенное в
Лемме 7.

Выберем положительное число  $\varepsilon$ меньшим радиуса
регулярной (погруженной) окрестности  $U_{reg}$ регулярных точек
отображения $d$ и  меньшим радиуса регулярной окрестности
$U_{\partial N(d)}$ (погруженной) критических точек этого
отображения. Рассмотрим погружение $f_1: M^{n-k} \looparrowright
\R^n$ в классе регулярной гомотопии погружения $f$ такое, что
выполнено неравенство
$$\dist(d \circ \kappa';f_1)_{C^0} < \varepsilon. $$
Погружение $f_1$ является скошенно-оснащенным посредством $\Xi_1$
с тем же характеристическим классом скошенного оснащения
$\kappa_M$, при этом $[(f_1,\Xi_1,\kappa_M)]=\alpha$, см. [A],
Следствие 23.

Докажем, что скошенно-оснащенное погружение $f_1$ допускает
$\J_b$--структуру. Обозначим через $N_1^{n-2k}$ многообразие
самопересечения погружения $f_1$. Определим отображение $\eta_b:
N_1^{n-2k} \to K(\J_b,1)$. Рассмотрим представление многообразия
$N_1^{n-2k}$ в виде объединения двух многообразий по общей части
границы: $N^{n-2k}_1 = N^{n-2k}_{N(d)} \cup_{\partial}
N^{n-2k}_{reg}$. В этой формуле $N^{n-2k}_{N(d)}$ -- многообразие
с краем, погруженно в регулярную (погруженную) окрестность
$U_{N(d)}$ многообразия с особенностями с краем $N(d)$ точек
самопересечения отображения $d$. Многообразие $N^{n-2k}_{reg}$ с
краем, погруженно в погруженную окрестность $U_{reg}$. Общий край
многообразий  $N^{n-2k}_{N(d)}$  $N^{n-2k}_{reg}$ является
замкнутым многообразие размерности $n-2k-1$, это многообразие
погружено в границу погруженной окрестности $U_{\partial N(d)}$.

Определим на многообразии $N^{n-2k}_1$ когомологический класс
$\kappa_{d,N_1} \in H^1(N^{n-2k}_1;\Z/2)$ отображением $N^{n-2k}_1
\stackrel{\kappa_{d,N_1}}{\longrightarrow} K(\J_d,1)$. На
подмногообразии $N^{n-2k}_{reg}$ определим отображение
$\kappa_{d,N_{reg}}: N^{n-2k}_{reg} \to K(\J_d,1)$ как композицию
проекции $N^{n-2k}_{reg} \to \RP^{n-k'}$ и включения $\RP^{n-k'}
\subset \RP^{\infty} = K(\J_d,1)$. На подмногообразии
$N^{n-2k}_{N(d)}$ определим отображение $\kappa_{N(d)}:
N^{n-2k}_{N(d)} \to K(\J_d,1)$ как композицию проекции
$N^{n-2k}_{N(d)} \to N(d)$ и отображения $\kappa_{d}: N(d) \to
K(\J_d,1)$, построенного в Лемме 7. Ограничения отображений
$\kappa_{N(d)}$, $\kappa_{d,N_{reg}}$ на края $\partial
N^{n-2k}_{N(d)}$ и $\partial N^{n-2k}_{reg}$ гомотопны, поскольку
отображение $\kappa_{N(d)}: N(d) \to K(\J_d,1)$ удовлетворяет на
$\partial N(d)$ граничному условию. Следовательно, определено
отображение $\kappa_{d,N_1}: N^{n-2k}_1 \to K(\J_d,1)$ как
результат склейки отображений $\kappa_{N(d)}$ и
$\kappa_{d,N_{reg}}$.

Когомологический класс  $\tau_{c,N_1} \in H^1(N^{n-2k}_1;\Z/2)$
определим как характеристический класс канонического двулистного
накрытия $\bar N^{n-2k}_1 \to N^{n-2k}_1$. Пара когомологических
классов $\kappa_{c,N_1}$, $\tau_{d,N_1}$ определяют искомое
отображение $\eta_{b,N_1}: N^{n-2k}_1 \to K(\J_b,1)$. Это
отображение однозначно характеризуется тем, что
$\eta_{b,N_1}^{\ast}(\kappa_d)=\kappa_{d,N_1}$,
$\eta_{b,N_1}^{\ast}(\kappa_c)=\kappa_{c,N_1}$.

Проверим уравнение ($\ref{etaa}$). Рассмотрим отображение
$\kappa_M': M^{n-k} \to \RP^{n-k-q-1}$ и рассмотрим отображение
$\kappa_{M_2}: M^{n-k}_2 \to  \RP^{n-k-q-1}$, которое определено
как ограничение $\kappa_M'$ на полный прообраз $M^{n-16k}_2=
\kappa'^{-1}(\RP^{n-16k-q-1})$ проективного подпространства
$\RP^{n-16k-q-1} \subset \RP^{n-k-q-1}$ коразмерности $15k$, в
предположении, что отображение $\kappa_M'$ трансверсально вдоль
этого подмногообразия.

Рассмотрим погружение $f_2: M^{n-k} \looparrowright \R^n$, которое
определено как ограничение погружения $f_1$ на подмногообразие
$i_{M} : M^{n-16k}_2 \subset M^{n-k}$. Погружение является
скошенно-оснащенным в коразмерности $16k$. Обозначим скошенное
оснащение этого погружения через $\Xi_2$, а характеристический
класс этого скошенного оснащения через $\kappa_{M_2} \in
H^1(M^{n-16k}_2;\Z/2)$. По построению $\kappa_{M_2} =
i_{M}^{\ast}\kappa_M$. Тройка $(f_2,\Xi_2,\kappa_{M_2})$
определяет элемент $J_{sf}^{16k}(\alpha) \in Imm^{sf}(n-16k,16k)$.

Рассмотрим многообразие точек самопересечения погружения $f_2$ и
обозначим это многообразие через $N^{n-32k}_2$. Определено
естественное вложение
\begin{eqnarray}\label{N2}
i_{N_1}: N^{n-32k}_2 \subset N^{n-2k}_1.
\end{eqnarray}
 Фундаментальный
класс рассматриваемого подмногообразия представляет цикл $i_{N_1,
\ast}([N_2]) \in H_{n-32k}  (N^{n-2k}_1;\Z/2)$. Этот цикл
двойственен в смысле Пуанкаре циклу $\eta_{N_1}^{15k} \in
H^{30k}(N^{n-2k}_1;\Z/2)$.

Докажем, что подмногообразие ($\ref{N2}$) целиком содержится в
подмногообразии $N^{n-2k}_{reg} \subset N^{n-2k}_1$, т.е
справедливо включение
\begin{eqnarray}\label{N2reg}
Im(i_{N_1}( N^{n-32k}_2)) \subset  N^{n-2k}_{reg} \subset
N^{n-2k}_1.
\end{eqnarray}
 Рассмотрим структурное отображение
$\eta_{N(d)}: (N(d),\partial N(d)) \to (K(\D_4,1),K(\I_b,1))$ и
обозначим снова через $\eta_{N(d)} \in H^2(N(d);\Z/2)$
характеристический класс, который индуцирован из универсального
класса $\eta_{[2]} \in H^2(K(\D_4,1);\Z/2)$ (из эйлерового класса
универсального $\D_4$-расслоения) отображением $\eta_{N(d)}$.
Класс гомологий $(\eta_{N_1}^{15k})^{(op)} \in
H_{n-32k}(N^{n-2k}_1;\Z/2)$ представлен циклом,  ограничение
которого на подмногообразие с краем $N^{n-2k}_{N(d)} \subset
N^{n-2k}_1$ совпадает с полным прообразом относительного цикла
$((\eta_{N(d)})^{15k})^{(op)}$ при проекции $N^{n-2k}_{N(d)} \to
N(d)$ на центральное подмногообразие с особенностями с краем в
регулярной погруженной окрестности.

По соображениям размерности, поскольку размерность
$dim(N(d))=n-2k-2q-2=n - \frac{n-n_s}{16}-n_s-2$ меньше
коразмерности подмногообразия $(\ref{N2})$, которая равна $30k=
\frac{15(n-n_s)}{16}$, относительный гомологический класс,
двойственный коциклу $(\eta_{N(d)})^{15k}$, представлен в $N(d)$
пустым многообразием. Это доказывает формулу $(\ref{N2reg})$.

Теперь для доказательства ($\ref{etaa}$) достаточно заметить, что
в силу формулы $(\ref{N2reg})$, коцикл $\eta_{b,N_1} \in
H^2(N^{n-2k}_1;\Z/2)$, ограниченный на подмногообразие
($\ref{N2}$), совпадает с ограничением коцикла $\eta_{N_1} \in
H^2(N^{n-2k}_1;\Z/2)$ на это же многообразие. Тем самым, формула
($\ref{etaa}$) и Теорема 6 доказаны.

\section{$\I_a \oplus \II_a$--структура (бициклическая структура)
 на $\Z/2^{[3]}$--оснащенном погружении}

Вспомним определение циклической подгруппы
 $\I_a \subset \D_4$, см. [A], раздел 2.
 Определим подгруппу
\begin{eqnarray} \label{iaa}
i_{a \oplus \aa}: \I_a \oplus \II_a \subset
 \Z/2^{[4]}.
\end{eqnarray}

В евклидовом пространстве $\R^8$ определен ортонормальный базис
$(\e_1, \dots, \e_8)$, при помощи которого определялась группа
$\Z/2^{[4]}$.

Обозначим образующие слагаемых группы  $\I_a \oplus \II_a$ через
$a$, $\aa$ соответственно. Опишем преобразования  из $\Z_2^{[4]}$,
которые соответствуют каждой образующей. Рассмотрим новый базис
$\{\f_1, \dots, \f_8\}$, определенный формулами
$\f_{2i-1}=\e_{2i-1}+ \e_{2i}$, $\f_{2i}= \e_{2i-1}- \e_{2i}$,
$i=1, \dots, 4$. Образующая $a$ порядка $4$ представлена поворотом
в каждой плоскости $(\f_1,\f_3)$, $(\f_5,\f_7)$ на угол
$\frac{\pi}{2}$ и одновременной центральной симметрией в плоскости
$(\f_2-\f_4, \f_6-\f_8)$. Образующая $\aa$ представлена поворотом
в плоскостях $(\f_2-\f_4, \f_6-\f_8)$, $(\f_2+\f_4,\f_6+\f_8)$ на
угол $\frac{\pi}{2}$ и одновременной центральной симметрией в
плоскости $(\f_1-\f_5,\f_3-\f_7)$.

Покажем, что группа преобразований  $\I_a \oplus \II_a$ имеет
инвариантные
 $(2,2,2,2)$-мерные подпространства, которые мы обозначим через  $\R^2_{a,+}$, $\R^2_{a,-}$,
 $\R^2_{\aa,+}$, $\R^2_{\aa,-}$.

Подпространство $\R^2_{a,+}$ порождено парой векторов
$(\f_1+\f_5,\f_3+\f_7)$. Подпространство $\R^2_{a,-}$ порождено
парой векторов $(\f_1-\f_5,\f_3-\f_7)$. Подпространство
$\R^2_{\aa,+}$ порождено парой векторов $(\f_2+\f_4,\f_6+\f_8)$.
Подпространство $\R^2_{\aa,-}$ порождено парой векторов
$(\f_2-\f_4,\f_6-\f_8)$

Образующая $a$ действует поворотом на угол $\frac{\pi}{2}$ в
каждой плоскости $\R^2_{a,+}$, $\R^2_{a,-}$ и симметрией в
плоскости $\R^2_{\aa,-}$, которая, очевидно, коммутирует с
действием образующей $\aa$ в этой плоскости. Образующая $\aa$
действует поворотом на угол $\frac{\pi}{2}$ в каждой плоскости
$\R^2_{\aa,+}$ $\R^2_{\aa,-}$ и одновременно центральной
симметрией в плоскости $\R^2_{a,-}$, которая, очевидно,
коммутирует с действием образующей $a$ в этом плоскости. Подгруппа
($\ref{iaa}$) определена.

Удобно перейти к новому базису $$(\h_{1,+}, \h_{2,+}, \h_{1,-},
\h_{2,-}, \hh_{1,+}, \hh_{2,+}, \hh_{1,-}, \hh_{2,-}).$$ Пары
векторов $(\h_{1,+},\h_{2,+})$, $(\h_{1,-}, \h_{2,-})$ задают
базисы в подпространствах $\R^2_{a,+}$,  $\R^2_{a,-}$
соответственно. Далее пары векторов $(\hh_{1,+},\hh_{2,+})$,
$(\hh_{1,-}, \hh_{2,-})$ задают базисы в подпространствах
$\R^2_{\aa,+}$,  $\R^2_{\aa,-}$ соответственно.

 Рассмотрим подгруппу
 $i_{a \oplus \dd,a \oplus \aa}: \I_a \oplus \II_d \subset \I_a \oplus \II_a$, которая определена
 прямой суммой группы $\I_a$ c элементарной подгруппой $\II_d$ второго слагаемого.
 Определено вложение $i_{a \oplus \dd}: \I_a \oplus \II_d \subset
 \Z/2^{[3]}$, согласованное с вложением ($\ref{iaa}$).  При этом определена коммутативная диаграмма:

\begin{eqnarray}
\label{a,aa}
\begin{array}{ccc}
\qquad \J_b & \stackrel {i_{b}} {\longrightarrow}& \qquad \D_4 \\
i_{b,a \oplus \dd} \downarrow &  & i_{[3]} \downarrow \\
 \qquad \I_a \oplus \II_d &  \stackrel {i_{a \oplus \dd}}
{\longrightarrow}& \qquad \Z/2^{[3]} \\
i_{a \oplus \dd,a \oplus \aa} \downarrow &  &  i_{[4]} \downarrow \\
\qquad \I_a \oplus \II_a & \stackrel{i_{a \oplus \aa}}{\longrightarrow} &  \qquad \Z/2^{[4]}.\\
\end{array}
\end{eqnarray}

 Определим также включение
$i_{b,a \oplus \dd}= i_{b,a} \oplus p_{b,\dd}: \J_b \subset \I_a
\oplus \II_{d}$. Гомоморфизм $i_{a,b}: \J_b \to \II_a$ определен
как композиция гомоморфизма проекции $\J_b \to \II_d$ и
гомоморфизма включения $\II_d \subset \II_a$. При этом определена
коммутативная диаграмма:

\begin{eqnarray} \label{a,dd}
\begin{array}{ccccc}
\J_b & \stackrel  {i_{b,a \oplus \dd}}{\longrightarrow} &
\D_4 &  \stackrel {i_{diag,\D_4}}{\longrightarrow}& \D_4 \oplus \D_4  \\
i_{b,a \oplus \dd} \bigcap \qquad & & &  & \bar i_{[3]} \bigcap \\
\I_a \oplus \II_d & & \stackrel{i_{a \oplus \dd}}{\longrightarrow} & & \Z/2^{[3]}.\\
\end{array}
\end{eqnarray}

\begin{definition}
 Пусть $\Z/2^{[3]}$--оcнащенное ($\Z/2^{[4]}$--оcнащенное)
погружение $(h,\Lambda,\zeta_L)$, $h: L^{n-4k} \looparrowright
\R^n$
 ($h: L^{n-8k} \looparrowright \R^n$)  представляет элемент
$z \in Imm^{\Z/2^{[3]}}(n-4k,4k)$ ($z \in
Imm^{\Z/2^{[4]}}(n-8k,8k)$). Скажем, что это
$\Z/2^{[3]}$--оcнащенное  ($\Z/2^{[4]}$--оcнащенное) погружение
является $\I_a \oplus \II_d$--оснащенным ($\I_a \oplus
\II_a$--оснащенным) погружением, если структурное отображение
$\zeta_L: L^{n-4k} \to K(\Z/2^{[3]},1)$ ($\zeta_L: L^{n-8k} \to
K(\Z/2^{[4]},1)$) представлено в виде композиции отображения
$\zeta_{a \oplus \dd,L}: L^{n-4k} \to K(\I_a \oplus \II_d,1)$
$\quad$ ($\zeta_{a \oplus \aa,L}: L^{n-8k} \to K(\I_a \oplus
\II_a,1)$) и отображения $i_{a\oplus \dd}: K(\I_a \oplus \II_d,1)
\to K(\Z/2^{[3]},1)$ $\quad$ ($i_{a\oplus \aa}: K(\I_a \oplus
\II_a,1) \to K(\Z/2^{[4]},1)$).
\end{definition}

Рассмотрим аналоги соотношения ($\ref{ibeta}$) для групп $\I_a
\oplus \II_d$ и $\I_a \oplus \II_a$--оснащенных погружений
соответственно.

Группа когомологий $H^4(K(\I_a \oplus \J_d,1);\Z/2)$ ($H^8(K(\I_a
\oplus \II_a,1);\Z/2)$ содержит элемент $\tau_{a \oplus \dd}$,
($\tau_{a \oplus \aa}$), который определяется нижеследующем
уравнением ($\ref{i3}$) (($\ref{i4}$)).

Рассмотрим отображение $i_{a \oplus \dd}: K(\I_a \oplus \II_d,1)
\to K(\Z/2^{[3]},1)$ ( $i_{a \oplus \aa}: K(\I_a \oplus \II_a,1)
\to K(\Z/2^{[4]},1)$) и рассмотрим обратный образ $i_{a \oplus
\dd}^{\ast}(\tau_{[3]})$  ($i_{a \oplus \aa}^{\ast}(\tau_{[4]})$)
характеристического эйлерового класса  $\tau_{[3]} \in
H^4(K(\Z/2^{[3]}),1);\Z/2)$ ($\tau_{[4]} \in
H^8(K(\Z/2^{[4]}),1);\Z/2)$) универсального расслоения. Определим
\begin{eqnarray}\label{i3}
 i_{a \oplus \dd}^{\ast} (\tau_{[3]}) = \tau_{a \oplus \dd},
\end{eqnarray}
\begin{eqnarray}\label{i4}
 i_{a \oplus \aa}^{\ast} (\tau_{[4]}) = \tau_{a \oplus \aa}.
\end{eqnarray}

 В разделе 1, для $\Z/2^{[d+1]}$--оснащенного погружения $(h,\Lambda,\zeta_L)$,
 наряду с $2^{d}$-мерным характеристическим классом $\zeta_L \in H^{2^{d}}(L^{n-k2^{d}};\Z/2)$,
 рассматривался также 2-мерный
характеристический класс $\bar \zeta_{[2],L} \in H^2(\bar
L_{[2]}^{n-k2^{d}k};\Z/2)$.
 Для отображения $\zeta_{a \oplus \dd}:
L^{n-4k} \to K(\I_a \oplus \II_d,1)$ ($\zeta_{a \oplus \aa}:
L^{n-8k} \to K(\I_a \oplus \II_a,1)$) аналогом характеристического
класса $\bar \zeta_{[2],L}$ служит характеристический класс $\bar
\zeta_{b,L} \in H^2(\bar L_{b}^{n-4k};\Z/2)$, при $d=3$ ($\bar
\zeta_{b,L} \in H^2(\bar L_{b}^{n-8k};\Z/2)$, при $d=4$).
Определим этот 2-мерный характеристический класс.

 Характеристический класс $\bar \zeta_{b,L}$ индуцирован из универсального класса $\tau_b \in
H^2(K(\J_b,1);\Z/2)$ при отображении $\bar \zeta_{b,L}: \bar
L^{n-4k}_{b} \to K(\J_b,1)$ ($\bar \zeta_{b,L}: \bar L^{n-8k}_{b}
\to K(\J_b,1)$). Отображение $\bar \zeta_{b,L}$ определено как
2-листное накрытие отображения  $\zeta_{N}$ относительно подгруппы
$i_{b,a \oplus \dd}: \J_b \subset \I_a \oplus \II_d$ (как
4-листное накрытие относительно подгруппы $i_{b,a \oplus \aa}:
\J_b \subset \I_a \oplus \II_a$). Это 2-листное (4-листное)
накрытие над многообразием $L^{n-4k}$ ($L^{n-8k}$) относительно
подгруппы $i_{b,a \oplus \dd}$ ($i_{b,a \oplus \aa}$) обозначается
через $\pi_{b,a \oplus \dd}: \bar L^{n-4k}_b \to L^{n-4k}$
($\pi_{b,a \oplus \aa}: \bar L^{n-8k}_b \to L^{n-8k}$).

\begin{definition}
Пусть $\D_4$--оснащенное погружение   $(g,\Psi,\eta_N)$, $g:
N^{n-2k} \looparrowright \R^n$ представляет элемент $y \in
Imm^{\D_4}(n-2k,2k)$, причем $n > 32k$. Пусть
$\Z/2^{[3]}$--оcнащенное погружение $(h,\Lambda,\zeta_L)$, $h:
L^{n-4k} \looparrowright \R^n$ является погружением точек
самопересечения погружения $g$ и представляет элемент
$z=\delta^{\Z/2^{[3]},k}(y) \in Imm^{ \Z/2^{[3]}}(n-4k,4k)$.
 Скажем, что $\D_4$--оснащенное  погружение
$(g,\Psi,\eta_N)$ допускает $\I_a \oplus \II_d$--структуру, если
существует отображение $\zeta_{a \oplus \dd,L}: L^{n-4k} \to
K(\I_a \oplus \II_d,1)$, удовлетворяющее уравнению:

\begin{eqnarray}\label{etaad}
\Theta_{\D_4}^k(y) = \langle \pi^{\ast}_{b,a \oplus
\dd,L}(\zeta_L^{7k}) \bar
  \zeta_{b,L}^{\frac{n-32k} {2}};[\bar L_{b}] \rangle,
\end{eqnarray}
где $[\bar L_{b}]$-- фундаментальный класс многообразия $\bar
L^{n-4k}_{b}$, характеристическое число $\Theta_{\D_4}^k$
определено по формуле ($\ref{44}$).
\end{definition}

\begin{definition}
Пусть $\Z/2^{[3]}$--оснащенное погружение   $(g,\Psi,\eta_N)$, $g:
N^{n-4k} \looparrowright \R^n$ представляет элемент $y \in
Imm^{\Z/2^{[3]}}(n-4k,4k)$, причем $n > 32k$. Пусть
$\Z/2^{[4]}$--оcнащенное погружение $(h,\Lambda,\zeta_L)$, $h:
L^{n-8k} \looparrowright \R^n$ является погружением точек
самопересечения погружения $g$ и представляет элемент
$z=\delta^{\Z/2^{[4]},k}(y) \in Imm^{\Z/2^{[4]}}(n-8k,8k)$.
 Скажем, что $\Z/2^{[3]}$--оснащенное  погружение
$(g,\Psi,\eta_N)$ допускает бициклическую структуру ($\I_a \oplus
\II_a$--структуру), если существует отображение $\zeta_{a \oplus
\aa,L}: L^{n-8k} \to K(\I_a \oplus \II_a,1)$, удовлетворяющее
уравнению:

\begin{eqnarray}\label{zetaad}
\Theta_{\Z/2^{[3]}}^k(y) = \langle \pi^{\ast}_{b,a \oplus
\aa,L}(\zeta_L)^{3k} \bar
  \zeta_{b,L}^{\frac{n-32k} {2}};[\bar L_{b}] \rangle,
\end{eqnarray}
где $[\bar L_{b}]$-- фундаментальный класс многообразия $\bar
L^{n-8k}_{b}$, характеристическое число $\Theta_{\Z/2^{[3]}}^k$
определено по формуле ($\ref{thetaz2^3}$).
\end{definition}

\subsubsection*{Пример}

Пусть $\D_4$-оснащенное ( $\Z/2^{[3]}$-оснащенное) погружение
$(g,\Psi,\eta_N)$, $g: N^{n-2k} \looparrowright \R^n$ ($g:
N^{n-4k} \looparrowright \R^n$) представляет элемент $y \in
Imm^{\D_4}(n-2k,2k)$ ($y \in Imm^{\Z/2^{[3]}}(n-4k,4k)$) и
является $\J_b$--оснащенным ($\I_a \oplus \II_d$--оснащенным)
погружением, причем $n
> 32k$.
Пусть $\Z/2^{[3]}$-оснащенное ($\Z/2^{[4]}$-оснащенное) погружение
$(h,\Lambda,\zeta_L)$, $h: L^{n-4k} \looparrowright \R^n$ ($h:
L^{n-8k} \looparrowright \R^n$) представляет элемент $z =
\delta^{\Z/2^{[3]},k} \in Imm^{\Z/2^{[3]}}(n-4k,4k)$ ($z =
\delta^{\Z/2^{[4]},k} \in Imm^{\Z/2^{[4]}}(n-8k,8k)$) и является
$\I_a \oplus \II_d$--оснащенным ($\I_a \oplus \II_a$--оснащенным)
погружением. Тогда $\D_4$--оснащенное ($\Z/2^{[3]}$--оснащенное)
погружение $(g,\Psi,\eta_N)$ допускает $\I_a \oplus
\II_d$--структуру ($\I_a \oplus \II_a$--структуру), заданную
редукцией $\zeta_{a \oplus \dd,L}$ ($\zeta_{a \oplus \aa,L}$)
структурного отображения $\zeta_L$.
\[  \]

Следующие теоремы аналогичны Теореме 7.



\begin{theorem}
Предположим, что $\D_4$--оснащенное погружение $(g,\Psi,\eta_N)$,
представляет элемент $y \in
Imm^{\D_4}(n-\frac{n-n_s}{16},\frac{n-n_s}{16})$, $n_s=2^s-2$,
$n>n_s$, $s \ge 6$. Предположим, что задано отображение
$\eta_{b,N}: N^{n-\frac{n-n_s}{16}} \to K(\J_b,1)$, при этом
выполнено уравнение:
\begin{eqnarray}\label{etaad'}
\Theta_{\D_4}^k(y) = \langle (\eta_{N}^{\frac{15(n-n_s)}{32}})
\eta_{b,N}^{\frac{n_s} {2}};[N] \rangle,
\end{eqnarray}
где $[N]$ -- фундаментальный класс многообразия $N^{n-
\frac{n-n_s}{16}}$, $\Theta_{\D_4}^k$ -- характеристическое число,
определенное по формуле ($\ref{44}$).
 Тогда элемент $J^{\Z/2^{[3]},\frac{n-n_s}{32}}(y)$ в группе
$Imm^{\Z/2^{[3]}}(n-\frac{n-n_s}{8},\frac{n-n_s}{8})$ представлен
$\Z/2^{[3]}$--оснащенным погружением $(h,\Lambda,\zeta_L)$,
которое допускает $\I_a \oplus \II_d$-структуру.
\end{theorem}

\begin{theorem}
Предположим, что $\Z/2^{[3]}$--оснащенное погружение
$(g,\Psi,\eta_N)$, представляет элемент $y \in
Imm^{\Z/2^{[3]}}(n-\frac{n-n_s}{8},\frac{n-n_s}{8})$, $n_s=2^s-2$,
$n>n_s$, $s \ge 6$. Предположим, что задано отображение $\eta_{a
\oplus \dd,N}:   N^{n-\frac{n-n_s}{8}} \to K(\I_a \oplus
\II_d,1)$, при этом выполнено уравнение:
\begin{eqnarray}\label{etaad''}
\Theta_{\Z/2^{[3]}}^k(y) = \langle (\pi_{b,a \oplus
\dd,N}^{\ast}\eta_{N})^{\frac{7(n-n_s)}{32}} \bar
\eta_{b,N}^{\frac{n_s} {2}};[\bar N_{b}] \rangle,
\end{eqnarray}
где накрывающее многообразие $\bar N_{b}^{n-\frac{n-n_s}{8}}$ и
отображение $\bar \eta_{b,N}: \bar N_{b}^{n-\frac{n-n_s}{8}} \to
K(\J_b,1)$ над отображением  $\eta_{[3],N}$  определяются
совершенно аналогичным понятиям в  $(\ref{etaad})$, через $[\bar
N_{b}]$ обозначен фундаментальный класс многообразия $\bar
N^{n-\frac{n-n_s}{8}}_b$.

 Тогда элемент $J^{\Z/2^{[3]},\frac{n-n_s}{32}}(y)$ в группе
$Imm^{\Z/2^{[3]}}(n-\frac{n-n_s}{8},\frac{n-n_s}{8})$ представлен
$\Z/2^{[3]}$--оснащенным погружением $(h,\Lambda,\zeta_L)$,
которое допускает $\I_a \oplus \II_d$-структуру.
\end{theorem}

\begin{corollary}
Предположим, что выполнены условия Теоремы 6 (т.е. $n_s$ является
натуральным числом вида  $2^s-2$, $n>n_s$, $s \ge 6$ и элемент $x
\in Imm^{sf}(n-\frac{n-n_s}{32},\frac{n-n_s}{32})$ допускает
ретракцию порядка  $q=\frac{n_s}{2}$). Тогда элемент
$\delta^k_{\Z^{[3]}} \circ \delta^k_{[2]}(x)$, определенный при
помощи композиции гомоморфизмов $(\ref{6})$, $k=\frac{n-n_s}{32}$,
представлен $\Z/2^{[3]}$--оснащенным погружением $(h,\Lambda,
\zeta_L)$, которое допускает бициклическую структуру.
\end{corollary}

\section{$\Q_a \oplus \QQ_a$--структура (бикватернионная структура)
 $\Z/2^{[5]}$--оснащенного погружения}

Вспомним определение кватернионной подгруппы
 $\Q_a \subset \Z/2^{[3]}$, которая содержит подгруппу $\I_a \subset \Q_a$, см. [A], раздел 2.

 Определим подгруппу
\begin{eqnarray} \label{iaaQ}
i_{\Q_a \oplus \QQ_a}: \Q_a \oplus \QQ_a \subset
 \Z/2^{[6]}.
\end{eqnarray}

 Рассмотрим базис $(\h_{1,+}, \h_{2,+}, \h_{1,+},
\h_{2,-}, \hh_{1,+}, \hh_{2,+}, \hh_{1,-}, \hh_{2,-})$
пространства $\R^8$, который был определен при построении
бициклической структуры.

Опишем базис пространства $\R^{32}$. Этот базис состоит из 32
векторов, разбитых на два подмножества из 16 векторов

\begin{eqnarray}\label{h}
\h_{1,\ast,\ast \ast},\h_{2,\ast,\ast \ast},\h_{3,\ast,\ast
\ast},\h_{4,\ast,\ast \ast},
\end{eqnarray}

\begin{eqnarray}\label{hh}
\hh_{1,\ast,\ast \ast}, \hh_{2,\ast,\ast \ast}, \hh_{3,\ast,\ast
\ast}, \hh_{4,\ast,\ast \ast},
\end{eqnarray}
где символы $\ast$, $\ast \ast$ независимо принимают значения
$+,-$. В каждом из подпространств, порожденных 4 векторами
($\ref{h}$), для которых символы $\ast$, $\ast \ast$ принимают
одинаковые значения, представление группы $\Q_a \oplus \QQ_a$
задано тривиально на втором слагаемом $\QQ_a$, при этом образующие
$\i$,$\j$,$\k$ слагаемого $\Q_a$ действуют преобразованиями,
определенными матрицами:

\begin{eqnarray}\label{Q a1}
\left(
\begin{array}{cccc}
0 & -1 & 0 & 0 \\
1 & 0 & 0 & 0 \\
0 & 0 & 0 & -1 \\
0 & 0 & 1 & 0 \\
\end{array}
\right),
\end{eqnarray}
\begin{eqnarray}\label{Q a2}
 \left(
\begin{array}{cccc}
0 & 0 & -1 & 0 \\
0 & 0 & 0 & 1 \\
1 & 0 & 0 & 0 \\
0 & -1 & 0 & 0 \\
\end{array}
\right),
\end{eqnarray}
\begin{eqnarray}\label{Q a3}
\left(
\begin{array}{cccc}
0 & 0 & 0 & -1 \\
0 & 0 & -1 & 0 \\
0 & 1 & 0 & 0 \\
1 & 0 & 0 & 0 \\
\end{array}
\right).
\end{eqnarray}

Образующая $\i$ ($\j$) первого слагаемого $\Q_a$ действует
центральной симметрией в каждом 4-мерном подпространстве
($\ref{hh}$), порожденном векторами, для которых индекс $\ast$
($\ast \ast$) принимает значение $-$, а индекс $\ast \ast$
($\ast$) принимает любое фиксированное значение, в оставшейся паре
пространств, порожденных векторами ($\ref{hh}$) действие
тождественно.

В каждом из подпространств, порожденных 4 векторами ($\ref{hh}$) ,
для которых символы $\ast$, $\ast \ast$ принимают  одинаковые
значения, представление группы $\Q_a \oplus \QQ_a$ задано
тривиально на первом слагаемом $\Q_a$, причем образующие
$\i$,$\j$,$\k$ слагаемого $\QQ_a$ представлены матрицами ($\ref{Q
a1}$),($\ref{Q a2}$),($\ref{Q a3}$).

Образующая $\i$ ($\j$) второго слагаемого $\QQ_a$ действует
центральной симметрией в каждом 4-мерном подпространстве
($\ref{h}$), порожденном векторами, для которых индекс $\ast$
($\ast \ast$) принимает значение $-$, а индекс $\ast \ast$
($\ast$) принимает любое фиксированное значение, в оставшейся паре
пространств, порожденных векторами ($\ref{h}$) действие
тождественно. Поскольку преобразование центральной симметрии лежит
в центре группы $\Q_a$ ($\QQ_a$),  преобразование любого элемента
из $\QQ_a$ ($\Q_a$) в неприводимом подпространстве представления
$\Q_a$ ($\QQ_a$) коммутирует с элементами указанного
представления. Подгруппа ($\ref{iaaQ}$)
 определена.

Рассмотрим подгруппу
 $i_{\Q_a \oplus \II_a, \Q_a \oplus \QQ_a}: \Q_a \oplus \II_a \subset \Q_a \oplus \QQ_a$, которая определена
 прямой суммой группы $\Q_a$ c циклической подгруппой $\II_a$ второго слагаемого.
 Определено вложение  $i_{\Q_a \oplus \II_a}: \Q_a \oplus \II_a \subset
 \Z/2^{[4]}$, индуцированное из вложения ($\ref{iaaQ}$).  При этом определена коммутативная диаграмма:

\begin{eqnarray}
\label{a,aaQ}
\begin{array}{ccc}
\qquad \J_b & \stackrel {i_b}{\longrightarrow}& \qquad \D_4 \\
i_{b,a \oplus \dd} \downarrow \qquad & & i_{[3]} \downarrow \\
\qquad \I_a \oplus \II_d &  \stackrel {i_{a \oplus \dd}}{\longrightarrow}& \qquad \Z/2^{[3]}\\
i_{a \oplus \dd,a \oplus \aa} \downarrow \qquad &  & i_{[4]} \downarrow \\
\qquad \I_a \oplus \II_a & \stackrel{i_{a \oplus
\aa}}{\longrightarrow} &
\qquad \Z/2^{[4]} \\
i_{\I_a \oplus \II_a, \Q_a \oplus \II_a} \downarrow \qquad & &  i_{[5]} \downarrow \\
\qquad \Q_a \oplus \II_a &  \stackrel {i_{\Q_a \oplus \II_a}}
{\longrightarrow}&
\qquad \Z/2^{[5]} \\
i_{\Q_a \oplus \II_a, \Q_a \oplus \QQ_a} \downarrow \qquad & & i_{[6]} \downarrow \\
\qquad \Q_a \oplus \QQ_a & \stackrel{i_{\Q_a \oplus \QQ_a}}{\longrightarrow} & \qquad \Z/2^{[6]},\\
\end{array}
\end{eqnarray}
которая включает в себя диаграмму ($\ref{a,aa}$) в качестве
поддиаграммы.

Следующее определение аналогично Определению 8.

\begin{definition}
 Пусть $\Z/2^{[5]}$--оcнащенное ($\Z/2^{[6]}$--оcнащенное)
погружение $(h,\Lambda,\zeta_L)$, $h: L^{n-16k} \looparrowright
\R^n$
 ($h: L^{n-32k} \looparrowright \R^n$)  представляет элемент
$z \in Imm^{\Z/2^{[5]}}(n-16k,16k)$ ($z \in
Imm^{\Z/2^{[6]}}(n-32k,32k)$). Скажем, что это
$\Z/2^{[5]}$--оcнащенное ($\Z/2^{[6]}$--оcнащенное) погружение
является $\Q_a \oplus \II_a$--оснащенным ($\Q_a \oplus
\QQ_a$--оснащенным) погружением, если структурное отображение
$\zeta_L: L^{n-16k} \to K(\Z/2^{[5]},1)$ ($\zeta_L: L^{n-32k} \to
K(\Z/2^{[6]},1)$) представлено в виде композиции $\zeta_{\Q_a
\oplus \II_a,L}: L^{n-16k} \to K(\Q_a \oplus \II_a,1)$
($\zeta_{\Q_a \oplus \QQ_a,L}: L^{n-32k} \to K(\Q_a \oplus
\QQ_a,1)$) и отображения $i_{\Q_a\oplus \II_a}: K(\Q_a \oplus
\II_a,1) \to K(\Z/2^{[5]},1)$ ($i_{\Q_a\oplus \QQ_a}: K(\Q_a
\oplus \QQ_a,1) \to K(\Z/2^{[6]},1)$).
\end{definition}

Рассмотрим аналоги соотношений ($\ref{ibeta}$), ($\ref{i3}$),
($\ref{i4}$) для групп $\Q_a \oplus \II_a$ и $\Q_a \oplus
\QQ_a$--оснащенных погружений соответственно.

Группа когомологий $H^{16}(K(\Q_a \oplus \II_a,1);\Z/2)$
 ($H^{32}(K(\Q_a \oplus \Q_a,1);\Z/2)$) содержит элемент $\tau_{\Q_a
\oplus \II_a}$, ($\tau_{\Q_a \oplus \QQ_a}$, который определяется
нижеследующим уравнением ($\ref{i5}$) (($\ref{i6}$)).

Рассмотрим отображение $i_{\Q_a \oplus \II_a}: K(\Q_a \oplus
\II_a,1) \to K(\Z/2^{[5]},1)$ $\quad$ ($i_{\Q_a \oplus \QQ_a}:
K(\Q_a \oplus \QQ_a,1) \to K(\Z/2^{[5]},1)$) и рассмотрим обратный
образ $i_{\Q_a \oplus \II_a}^{\ast}(\tau_{[5]})$  ($i_{\Q_a \oplus
\QQ_a}^{\ast}(\tau_{[6]})$) характеристического эйлерового класса
$\tau_{[5]} \in H^{16}(K(\Z/2^{[5]}),1);\Z/2)$ ($\tau_{[6]} \in
H^{32}(K(\Z/2^{[6]}),1);\Z/2)$) универсального расслоения.
 Определим
\begin{eqnarray}\label{i5}
 i_{\Q_a \oplus \II_a}^{\ast} (\tau_{[5]}) = \tau_{\Q_a \oplus \II_a},
\end{eqnarray}
\begin{eqnarray}\label{i6}
 i_{\Q_a \oplus \QQ_a}^{\ast} (\tau_{[6]}) = \tau_{\Q_a \oplus \QQ_a}.
\end{eqnarray}

 Для отображения $\zeta_{\Q_a \oplus \II_a}:
L^{n-16k} \to K(\Q_a \oplus \II_a,1)$ $\zeta_{\Q_a \oplus \QQ_a}:
L^{n-32k} \to K(\I_a \oplus \II_d,1)$ аналогом характеристического
класса $\bar \zeta_{[2],L}$  служит характеристический класс $\bar
\zeta_{b,L} \in H^2(\bar L_{b}^{n-16k};\Z/2)$, при $d=5$ ($\bar
\zeta_{b,L} \in H^2(\bar L_{b}^{n-32k};\Z/2)$, при $d=6$).
Определим этот 2-мерный характеристический класс.

  Отображение $\bar \zeta_{b,L}$
определено как 8-листное накрытие над отображением $\zeta_{[5],L}$
 относительно подгруппы $i_{b,\Q_a \oplus \II_a}: \J_b \subset
\Q_a \oplus \II_a$ ( как 16-листное накрытие над отображением
 $\zeta_{[6],L}$ относительно подгруппы $i_{b, \Q_a \oplus \QQ_a}: \J_b
\subset \Q_a \oplus \QQ_a$). Над многообразием $L^{n-16k}$
($L^{n-32k}$) указанное накрытие обозначается через $\pi_{b,\Q_a
\oplus \II_a,L}$ ($\pi_{b,\Q_a \oplus \QQ_a,L}$).

\begin{definition}
Пусть $\Z/2^{[4]}$--оснащенное погружение   $(g,\Psi,\eta_N)$, $g:
N^{n-8k} \looparrowright \R^n$ представляет элемент $y \in
Imm^{\Z/2^{[4]}}(n-8k,8k)$, причем $n > 32k$. Пусть
$\Z/2^{[5]}$--оcнащенное погружение $(h,\Lambda,\zeta_L)$, $h:
L^{n-16k} \looparrowright \R^n$ является погружением точек
самопересечения погружения $g$ и представляет элемент
$\delta_{\Z/2^{[5]}}^k(y) \in Imm^{ \Z/2^{[5]}}(n-16k,16k)$.
 Скажем, что $\Z/2^{[4]}$--оснащенное  погружение
$(g,\Psi,\eta_N)$ допускает $\Q_a \oplus \II_a$--структуру, если
существует отображение $\zeta_{\Q_a \oplus \II_a,L}: L^{n-16k} \to
K(\Q_a \oplus \II_a,1)$, удовлетворяющее уравнению:

\begin{eqnarray}\label{etaQaIa}
\Theta_{\Z/2^{[4]}}^k(y) = \langle \pi^{\ast}_{b,\Q_a \oplus
\II_a, L}(\zeta^{k}) \bar
  \zeta_{b}^{\frac{n-32k} {2}};[\bar L_{b}] \rangle,
\end{eqnarray}
где $[\bar L_{b}]$-- фундаментальный класс многообразия $\bar
L^{n-16k}_{b}$, характеристическое число
$\Theta_{\Z/2^{[4]}}^k(y)$ определено по формуле
($\ref{thetaz2^3}$) при $d=4$.
\end{definition}

\begin{definition}
Пусть $\Z/2^{[5]}$--оснащенное погружение   $(g,\Psi,\eta_N)$, $g:
N^{n-16k} \looparrowright \R^n$ представляет элемент $y \in
Imm^{\Z/2^{[5]}}(n-16k,16k)$, причем $n > 32k$. Пусть
$\Z/2^{[6]}$--оcнащенное погружение $(h,\Lambda,\zeta_L)$, $h:
L^{n-32k} \looparrowright \R^n$ является погружением точек
самопересечения погружения $g$ и представляет элемент
$\delta_{\Z/2^{[6]}}^k(y) \in Imm^{\Z/2^{[6]}}(n-32k,32k)$.
 Скажем, что $\Z/2^{[5]}$--оснащенное  погружение
$(g,\Psi,\eta_N)$ допускает бикватернионную структуру ($\Q_a
\oplus \QQ_a$--структуру), если существует отображение
$\zeta_{\Q_a \oplus \QQ_a,L}: L^{n-32k} \to K(\Q_a \oplus
\QQ_a,1)$, удовлетворяющее уравнению:

\begin{eqnarray}\label{zetQaQa}
\Theta_{\Z/2^{[5]}}^k(y) = \langle \bar
  \zeta_{b}^{\frac{n-32k} {2}};[\bar L_{b}] \rangle,
\end{eqnarray}
где $[\bar L_{b}]$-- фундаментальный класс многообразия $\bar
L^{n-32k}_{b}$, характеристическое число
$\Theta_{\Z/2^{[5]}}^k(y)$ определено по формуле
($\ref{thetaz2^3}$) при $d=5$.
\end{definition}

\subsubsection*{Пример}

Пусть $\Z/2^{[4]}$-оснащенное      ($\Z/2^{[5]}$-оснащенное)
погружение $(g,\Psi,\eta_N)$, $g: N^{n-8k} \looparrowright \R^n$
($g: N^{n-16k} \looparrowright \R^n$) представляет элемент $y \in
Imm^{\Z/2^{[4]}}(n-8k,8k)$ ($y \in Imm^{\Z/2^{[5]}}(n-16k,16k)$) и
является $\Q_a \oplus \II_a$--оснащенным ($\Q_a \oplus
\QQ_a$--оснащенным) погружением, причем $n
> 32k$.
Пусть $\Z/2^{[5]}$-оснащенное ($\Z/2^{[6]}$-оснащенное) погружение
$(h,\Lambda,\zeta_L)$, $h: L^{n-16k} \looparrowright \R^n$ ($h:
L^{n-32k} \looparrowright \R^n$) представляет элемент $
\delta_{\Z/2^{[5]}}^k(y) \in Imm^{\Z/2^{[5]}}(n-16k,16k)$ ($
\delta_{\Z/2^{[6]}}^k(y) \in Imm^{\Z/2^{[6]}}(n-32k,32k)$) и
является $\Q_a \oplus \II_a$--оснащенным ($\Q_a \oplus
\QQ_a$--оснащенным) погружением. Тогда $\Z/2^{[4]}$--оснащенное
($\Z/2^{[5]}$--оснащенное) погружение $(g,\Psi,\eta_N)$ допускает
$\Q_a \oplus \II_a$--структуру ($\Q_a \oplus \QQ_a$--структуру),
заданную редукцией $\zeta_{\Q_a \oplus \II_a,L}$ ($\zeta_{\Q_a
\oplus \II_a,L}$) структурного отображения $\zeta_{[4],L}$
($\zeta_{[5],L}$).
\[  \]

Следующие теоремы аналогичны Теоремам 6,11,12.

\begin{theorem}
Предположим, что $\Z/2^{[4]}$--оснащенное погружение
$(g,\Psi,\eta_N)$, представляет элемент $y \in
Imm^{\Z/2^{[4]}}(n-\frac{n-n_s}{4},\frac{n-n_s}{4})$, $n_s=2^s-2$,
$n>5n_s+168$, $s \ge 6$. Предположим, что задано отображение
$\eta_{a \oplus \aa,N}: N^{n-\frac{n-n_s}{4}} \to K(\I_a \oplus
\II_a,1)$, при этом выполнено уравнение:
\begin{eqnarray}\label{etaQaIa'}
\Theta_{\Z/2^{[4]}}^k(y) = \langle \pi_{b,a \oplus
\dd,N}^{\ast}(\eta_{[4],N}^{\frac{3(n-n_s)}{32}}) \bar
\eta_{b,N}^{\frac{n_s} {2}};[\bar N_b] \rangle,
\end{eqnarray}
где накрывающее многообразие $\bar N_{b}^{n-\frac{n-n_s}{4}}$ и
отображение $\bar \eta_{b,N}: \bar N_{b}^{n-\frac{n-n_s}{4}} \to
K(\J_b,1)$ над отображением  $\eta_{[4],N}$  определяются
совершенно аналогичным понятиям в  $(\ref{etaQaIa})$. Тогда
элемент $J_{\Z/2^{[5]}}^{\frac{n-n_s}{32}}(y)$ в группе
$Imm^{\Z/2^{[5]}}(n-\frac{n-n_s}{2},\frac{n-n_s}{2})$ представлен
$\Z/2^{[5]}$--оснащенным погружением $(h,\Lambda,\zeta_L)$,
которое допускает $\Q_a \oplus \II_a$-структуру.
\end{theorem}

\begin{theorem}
Предположим, что $\Z/2^{[5]}$--оснащенное погружение
$(g,\Psi,\eta_N)$, представляет элемент $y \in
Imm^{\Z/2^{[5]}}(n-\frac{n-n_s}{2},\frac{n-n_s}{2})$, $n_s=2^s-2$,
$n > \frac{5n_s + 84}{3}$, $s \ge 6$. Предположим, что задано
отображение $\eta_{\Q_a \oplus \II_a,N}: N^{n-\frac{n-n_s}{2}} \to
K(\Q_a \oplus \II_a,1)$, при этом выполнено уравнение:
\begin{eqnarray}\label{etaQaQa'}
\Theta_{\Z/2^{[5]}}^k(y) = \langle \pi_{b,\Q_a \oplus
\II_a,N}^{\ast}(\eta_{N}^{\frac{n-n_s}{32}}) \bar
\eta_{b,N}^{\frac{n_s} {2}};[\bar N_{b}] \rangle,
\end{eqnarray}
где накрывающее многообразие $\bar N_{b}^{n-\frac{n-n_s}{2}}$ и
отображение $\bar \eta_{b,N}: \bar N_{b}^{n-\frac{n-n_s}{2}} \to
K(\J_b,1)$ над отображением  $\eta_{[4],N}$  определяются
аналогично понятиям в  $(\ref{etaQaIa})$, $(\ref{etaQaIa'})$.
Тогда элемент $J_{\Z/2^{[6]}}^{\frac{n-n_s}{32}}(y)$ в группе
$Imm^{\Z/2^{[6]}}(n_s,n-n_s)$ представлен $\Z/2^{[6]}$--оснащенным
погружением $(h,\Lambda,\zeta_L)$, которое допускает
бикватернионную структуру.
\end{theorem}

\begin{corollary}
Предположим, что выполнены условия Теоремы 6 (т.е. задано
натуральное число $n_s$ вида  $2^s-2$, $n>n_s$, $s \ge 6$ и задан
элемент $x \in Imm^{sf}(n-\frac{n-n_s}{32},\frac{n-n_s}{32})$,
который допускает ретракцию порядка  $q=\frac{n_s}{2}$). Тогда
элемент
\begin{eqnarray} \label{ddddx}
\delta^k_{\Z^{[5]}} \circ \delta^k_{\Z^{[4]}} \circ
\delta^k_{\Z^{[3]}} \circ \delta^k_{\D_4}(x),
\end{eqnarray}
определенный при помощи последовательной композиции гомоморфизмов
($\ref{6}$) $k=\frac{n-n_s}{32}$, представлен
$\Z/2^{[5]}$--оснащенным погружением $(g,\Psi, \eta_N)$, которое
допускает бикватернионную структуру.
\end{corollary}

\section{Решение Проблемы Кервера}

В этом разделе мы докажем следующий результат.
\subsubsection*{Основная Теорема}

Существует натуральное $l_0$ такое, что для произвольного
натурального $l \ge l_0$,  $n=2^l-2$, инвариант Кервера,
определенный формулой (1) тривиальный.
\[  \]

\subsubsection*{Доказательство Основной Теоремы}

Определим $n_s=2^{12}-2$ и обозначим $\frac{n-n_s}{32}$ через $k$.
По Теореме 29 (Tеорема о ретракции) существует натуральное $l_0$
такое, что для произвольного натурального  $l \ge l_0$
произвольный элемент $x \in Imm^{sf}(n-k,k)$ допускает ретракцию
порядка $\frac{n_s}{2}=2^{11}-1$. Не ограничивая общности будем
считать, что $l_0 \ge 13$. Поскольку выполнены размерностные
условия, сформулированные в Теоремах 6,17,18
 по Cледствию
19, в классе кобордизма элемента ($\ref{ddddx}$) существует
$\Z/2^{[5]}$-оснащенное погружение $(g,\Psi,\eta_N)$, допускающее
бикватернионную структуру.

Рассмотрим многообразие самопересечения $L^{n_s}$  погружения $g$,
$\dim(L^{n_s})=n_s$. Многообразие $L^{n_s}$ снабжено отображением
\begin{eqnarray}
\zeta_{\Q_a \oplus \QQ_a}: L^{n_s} \to K(\Q_a \oplus \QQ_a,1).
\end{eqnarray}
Это отображение представляется прямым произведением двух
отображений
\begin{eqnarray}
\zeta_{\Q_a \oplus \QQ_a} = \zeta_{\Q_a} \times \zeta_{\QQ_a}:
L^{n_s} \to K(\Q_a ,1) \times K(\QQ_a,1).
\end{eqnarray}

Нам потребуется следующая лемма. Рассмотрим образующий класс
когомологий $\tau_{\Q_a}  \times \tau_{\QQ_a} \in H^8(K(\Q_a
\oplus \QQ_a,1);\Z/2)$. Обозначим $\zeta_{\Q_a \oplus
\QQ_a}^{\ast}(\rho_{\Q_a} \times \rho_{\QQ_a}) \in
H^8(L^{n_s};\Z/2)$ через $\rho_{\Q_a \oplus \QQ_a, L}$. Рассмотрим
16-листное накрытие $\pi_{tot,L}: \bar L^{n_s}_{[2]} \to L^{n_s}$.
Индуцируем класс $\rho_{\Q_a \oplus \QQ_a}$ на накрывающее в класс
$\pi_{tot,L}^{\ast}(\rho_{\Q_a \oplus \QQ_a, L}) \in H^8(\bar
L^{n_s}_{[2]};\Z/2)$.

\begin{lemma}
В группе $H^8(\bar L_{[2]}^{n_s};\Z/2)$ выполнено равенство:
\begin{eqnarray} \label{totQQ}
\pi_{tot,L}^{\ast}(\rho_{\Q_a \oplus \QQ_a,L}) = \bar
\zeta_{b,L}^{4}.
\end{eqnarray}
\end{lemma}

\subsubsection*{Доказательство Леммы 20}

Доказательство проводится прямым вычислением и опускается.
\[  \]

 Рассмотрим подмногообразие $i_K: K^{14} \subset L^{n_s}$, двойственное
когомологическому классу  $\rho_{\Q_a \oplus \QQ_a}^{2^9-2} \in
H^{n_s-14}(L^{n_s};\Z/2)$.  Многообразие $K^{14}$ снабжено
отображением $\zeta_{\Q_a \oplus  \QQ_a,K}: K^{14} \to K(\Q_a
\oplus \QQ_a,1)$, которое определено как ограничение отображения
$\zeta_{\Q_a \oplus \QQ_a, L}$ на подмногообразие $K^{14} \subset
L^{n_s}$. Определено 16-листное накрытие  $\pi_{tot,K}: \bar
K_{b}^{14} \to K^{14}$. Определен индуцированный
характиристический класс $\bar \zeta_{b,K} \in H^2(\bar
K_{b};\Z/2)= \bar
  i_{b,\Q_a \oplus \QQ_a}^{\ast}(\bar \zeta_{b,L})$.

Формула $(\ref{zetQaQa})$ приобретает вид:
\begin{eqnarray} \label{arfQQ}
 \Theta_{sf}^k(y) = \langle \bar
\zeta_{b,L}^{\frac{n_s}{2}};[\bar L_{[2]}] \rangle.
\end{eqnarray}
По Лемме 20 характеристическое число ($\ref{arfQQ}$) равно
характеристическому числу
\begin{eqnarray} \label{arQQ}
\langle \bar \zeta_{b,L}^{7};[\bar L_{b}] \rangle.
\end{eqnarray}

Докажем, что характеристическое число $(\ref{arQQ})$ обращается в
нуль. Над пространством $K(\Q_a,1)$  ($K(\Q_a,1)$) определено
4-мерное расслоение $\rho_{\Q_a}$ ($\rho_{\QQ_a}$) со структурной
группой $\Q_a$ ($\QQ_a$), см. [A], раздел 2. Следовательно, над
пространством $K(\Q_a \oplus \QQ_a,1)= K(\Q_a,1) \times
K(\QQ_a,1)$ также определены 4-мерные расслоения $\chi_{\Q_a}$ и
$\chi_{\QQ_a}$) по формуле $\chi_{\Q_a} = p^{\ast}_{\Q_a \oplus
\QQ_a,\Q_a}(\rho_{\Q_a})$,   $\chi_{\QQ_a} = p^{\ast}_{\Q_a \oplus
\QQ_a,\Q_a}(\rho_{\QQ_a})$.

Нормальное расслоение к многообразию $K^{14}$ представлено суммой
Уитни
\begin{eqnarray} \label{nuK}
\nu_K = (2^9-2)(\zeta_ {\Q_a,K} \oplus \zeta_{\QQ_a,K}) \oplus
\Omega_K,
\end{eqnarray}
где $\zeta_{\Q_a,K} = \zeta_{\Q_a \oplus \QQ_a}(\chi_{\Q_a})$,
$\zeta_{\QQ_a,K} = \zeta_{\Q_a \oplus \QQ_a} (\chi_{\QQ_a})$,
$\Omega_K$--ограничение нормального расслоения многообразя
$L^{n_s}$ на подмногообразие $K^{14}\subset L^{n_s}$.

Расслоение $\Omega_K$ изоморфно сумме Уитни $k$ копий
$32$--мерного расслоения $\zeta_{[6]}$ со структурной группой
$\Z/2^{[6]}$. Поскольку $k \equiv 0 \pmod{2^{7}}$, а база $K^{14}$
имеет размерность $14$, то рассуждая по аналогии с Предложением 34
из [A], заключаем, что расслоение $\Omega_K$ тривиально.

Рассмотрим многообразие $-K^{14}$, которое получено из $K^{14}$
изменением ориентации. Рассмотрим отображение $F= id \cup -id:
K^{14} \cup -K^{14} \to K^{14}$.  Вычислим характеристические
циклы $p_1(\nu_K)$, $p_1(\nu_{-K})$ нормального расслоения
многообразий $K^{14}$, $-K^{14}$, которое определено на каждой
копии многообразия по формуле ($\ref{nuK}$) (см. аналогичное
вычисление в Предложении 33 из [A])  Согласно вычислениям,
получим, что цикл
\begin{eqnarray} \label{doc}
F_{\ast}([p_1 (\nu_K)]^{op} \cup [p_1(\nu_{-K})]^{op} ]) \in
H_{10}(K^{14};\Z)
\end{eqnarray}
двойственен в смысле Пуанкаре коциклу $\zeta_{\Q_a \oplus
\QQ_a,K}^{\ast}(\gamma)$, где $\gamma \in H^{10}(K(\Q_a \oplus
\QQ_a,1;\Z)$ ненулевой когомологический класс, индуцированный из
элемента $4t \in Im(H^4(K(\Q_a),1;\Z)$ при гомоморфизме
$H^4(K(\Q_a ,1);\Z) \stackrel{p^{\ast}_{\Q_a,\Q_a \oplus
\QQ_a}}{\longrightarrow} K(\Q_a \oplus \QQ_a,1)$.

Условие нетривиальности характеристического числа ($\ref{arfQQ}$)
эквивалентно тому, что ориентированный цикл $\zeta_{\Q_a \oplus
\QQ_a,K,\ast}([K]) \in H_{14}(K(\Q_a \oplus \QQ_a,1);\Z)$, при
разложении по стандартному базису, содержит моном $a_{\Q_a}
\otimes a_{\QQ_a}$, где $a_{\Q_a} \in H_7(K(\Q_a,1);\Z)$,
$a_{\QQ_a} \in H_7(K(\QQ_a,1);\Z)$ -- образущие. Если число
($\ref{arfQQ}$) не равно нулю, то цикл ($\ref{doc}$) не обращается
в нуль, поскольку частное класса гомологий $a_{\Q_a} \otimes
a_{\QQ_a}$ на класс когомологий $\gamma$ не обращается в нуль. С
другой стороны, по построению отображение $F$ кобордантно пустому
и характеристический цикл ($\ref{doc}$) обращается в нуль.

Доказано, что на классе оснащенного кобордизма $x$
характеристическое число (6) равно нулю. Проблема Кервера решена.

\section{Доказательство Теорем 11 и 12}

Начнем доказательство со следующей конструкции. Определим число
$n_s$ из условия Теоремы 6. Рассмотрим многообразие
  $ZZ= S^{n-\frac{n-n_s}{16}+3}/\i \times S^{n-\frac{n-n_s}{16}+3}/\i$.
 Это многообразие является прямым произведением стандартных линзовых пространств
 $\pmod{4}$. Заметим, что $\dim(ZZ) > n$. Выберем внутри
 многообразия $ZZ$ подмногообразие $Z$ с особенностями в коразмерности
 2, такое, что $Z$ вкладывается в $\R^n$, в частности, $\dim(Z)<n$. По поводу
 многообразий с особенностями см. [B-R-S].

Рассмотрим в многообразии
 $ZZ$ семейство подмногообразий
$$Z_j, \quad j=0, \dots, j_{max}, \qquad j_{max}=\frac{15n+n_s+64}{32}$$
размерности $n-\frac{n-n_s}{16}+4$ и коразмерности
$n-\frac{n-n_s}{16}+2$, определенное по формуле
$$Z_0 = S^{n-\frac{n-n_s}{16}+3}/\i  \times S^1/\i, \quad Z_1= S^{n-\frac{n-n_s}{16}+1}/\i
\times S^3/\i, \quad \dots,$$
$$ Z_j = S^{n-\frac{n-n_s}{16}+3-2j}/\i \times S^{2j+1}/\i, \quad
\dots, \quad Z_{j_{max}} = S^1/\i \times
S^{n-\frac{n-n_s}{16}+3}/\i.$$ Вложение соответствующего
подмногообразия семейства в многообразие $ZZ$ определено как
прямое произведение двух стандартных вложений.

Объединение $\cup_{j=0} ^{j_{max}} Z_j$ семейства подмногообразий
$\{Z_j\}$ многообразия $ZZ$ является полиэдром
(стратифицированноым подмногообразием с особенностями в
коразмерности 2) размерности $n-\frac{n-n_s}{16}+4$, которое
обозначим $Z_{a \oplus \aa} \subset ZZ$.

Рассмотрим цепочку подгрупп
\begin{eqnarray}\label{zep}
\J_b \stackrel{i_{b,a\oplus \dd}}{\longrightarrow} \I_a \oplus
\II_d \stackrel {i_{a \oplus \dd, a \oplus \aa}}{\longrightarrow}
\I_a \oplus \II_a,
\end{eqnarray}
которая индуцирует башню 2-листных накрытий:
\begin{eqnarray} \label{zz}
XX \stackrel{p_{XX,YY}}{\longrightarrow} YY
\stackrel{p_{YY,ZZ}}{\longrightarrow} ZZ.
\end{eqnarray}

 Определим башню накрытий
\begin{eqnarray} \label{z}
X_{b}   \stackrel{p_{Y_{a \oplus \dd}}}{\longrightarrow} Y_{a
\oplus \dd}  \stackrel{ p_{Z_{a \oplus \aa}}}{\longrightarrow}
Z_{a \oplus \aa}.
\end{eqnarray}
 Эта башня накрытий индуцированна из цепочки подгрупп
 ($\ref{zep}$) при помощи вложения
$Z_{a \oplus \aa} \subset ZZ$. Таким образом, накрывающее
пространство $X_{b}$ в башне ($\ref{z}$) также является
стратифицированным многообразием с особенностями в коразмерности
2. Это пространство определено явным образом как объединение
семейства подмногообразий в $XX = \RP^{n-\frac{n-n_s}{16}+3}
\times \RP^{n-\frac{n-n_s}{16}+3}$, определенных по формуле:
$$X_0
= \RP^{n-\frac{n-n_s}{16}+3} \times \RP^{1}, \quad \dots, \quad
X_j = \RP^{n-\frac{n-n_s}{16}+3-2j} \times \RP^{2j+1}, \dots $$
$$ X_{j_{max}} = \RP^{1} \times
\RP^{n-\frac{n-n_s}{16}+3}.$$

Среднее накрывающее пространство $Y_{a \oplus \dd}$ в башне
($\ref{z}$) также является стратифицированным многообразием с
особенностями в коразмерности 2. Это пространство является
объединением семейства подмногообразий в $YY =
S^{n-\frac{n-n_s}{16}+3}/\i \times \RP^{n-\frac{n-n_s}{16}+3}$,
определенных по формуле:
$$Y_0
= S^{n-\frac{n-n_s}{16}+3}/\i \times \RP^{1}, \quad \dots, \quad
Y_j = S^{n-\frac{n-n_s}{16}+3-2j}/\i \times \RP^{2j+1}, \dots $$
$$ X_{j_{max}} = S^{1}/\i \times
\RP^{n-\frac{n-n_s}{16}+3}.$$

Определены отображения $\eta_{X}: X_{b} \to K(\J_b,1)$, $\eta_{Y}:
Y_{a \oplus \dd} \to K(\I_a \oplus \II_d,1)$, $\eta_{Z}: Z_{a
\oplus \aa} \to K(\I_a \oplus \II_a,1)$, согласованные с
включением подгрупп ($\ref{zep}$) и башней накрытий ($\ref{z}$).
Отображение $\eta_X$ представлено прямым произведением отображений
$\eta_{X,d} \times \eta_{X,\dd}:  X_{b} \to K(\I_d,1) \times
K(\II_d,1)$. Отображение $\eta_Y$ представлено прямым
произведением отображений $\eta_{Y,a} \times \eta_{Y,\dd}: Y_{a
\oplus \dd} \to K(\I_a,1) \times K(\II_d,1)$. Отображение
$\zeta_Z$ представлено прямым произведением отображений
$\zeta_{Z,a} \times \zeta_{Z,\aa}: Z_{a \oplus \aa} \to K(\I_a,1)
\times K(\II_a,1)$.

Определим многообразие с особенностями $J$. Для произвольного
 $j= 0, \dots, j_{max}$ определим пространство
$J_j = S^{n-\frac{n-n_s}{16}-2j+3} \times S^{2j+1}$. Сферы
$S^{n-\frac{n-n_s}{16}-2j+3}$, $S^{2j+1}$ переобозначим для
краткости через $J_{j,1}$, $J_{j,2}$ соответственно. Справедлива
формула $J_j = J_{j,1} \times J_{j,2}$

Определено стандартное включение $i_{J_j}: J_{j,1} \times J_{j,2}
\subset S^{\frac{n-n_s}{16}+3} \times S^{\frac{n-n_s}{16}+3}$, где
каждый сомножитель включается в сферу-образ как стандартная
подсфера, лежащая в подпространстве с соответствующим числом
первых ненулевых координат. Объединение $\cup_{j=0}^{j_{max}}
Im(i_{J_j})$ образов всех этих вложений является искомым
пространством, которое обозначим через $J  \subset
S^{\frac{n-n_s}{4}+3} \times S^{\frac{n-n_s}{4}+3}$.

Определим разветвленное накрытие
\begin{eqnarray}\label{pz}
\varphi_Z: Z_{a \oplus \aa} \to J.
\end{eqnarray}
Для каждого
 $j= 0, \dots, j_{max}$ определим
 отображение $\varphi_j: Z_j \to J_j$. Определены стандартные разветвленные накрытия (см.
[A], глава 3, Определение отображения $d$):
\begin{eqnarray}\label{p1}
p_{Z_j,1}: S^{\frac{n-n_s}{16}-2j+3}/\i \to
S^{\frac{n-n_s}{16}-2j+3},
\end{eqnarray}
\begin{eqnarray}\label{p2}
 p_{Z_j,2}: S^{2j+3}/\i \to  S^{2j+3}.
\end{eqnarray}
Отображение $\varphi_{Z_j}$ определено как декартово произведение
 разветвленных накрытий
$$ \varphi_{Z_j} = p_{Z_j,1} \times p_{Z_j,2}: Z_j=S^{\frac{n-n_s}{16}-2j+3}/\i  \times S^{2j+1}/\i
\to J_{j,1} \times J_{j,2}. $$

Определено разветвленное накрытие ($\ref{pz}$) в результате
склейки разветвленных накрытий $\varphi_{Z_j}$ по подпространствам
попарных пересечений семейства подмногообразий $Z_j$ в
многообразии $ZZ$.

Определим разветвленное накрытие
\begin{eqnarray}\label{py}
\varphi_Y: Y_{a \oplus \dd} \to J.
\end{eqnarray}
Для каждого
 $j= 0, \dots, j_{max}$ определим
 отображение $\varphi_j: Y_j \to J_j$. Определены стандартные разветвленные
 накрытия
\begin{eqnarray}\label{pp1}
p_{Y_j,1}: S^{n-\frac{n-n_s}{16}-2j+3}/\i \to
S^{n-\frac{n-n_s}{16}-2j+3},
\end{eqnarray}
\begin{eqnarray}\label{pp2}
 p_{Y_j,2}: \RP^{2j+3} \to   S^{2j+3}.
\end{eqnarray}
Отображение $\varphi_{Y_j}$ определено как декартово произведение
 разветвленных накрытий
$$ \varphi_{Y_j} = p_{Y_j,1} \times p_{Y_j,2}: Y_j=S^{\frac{n-n_s}{16}-2j+3}/\i  \times \RP^{2j+1}
\to J_{j,1} \times J_{j,2}. $$

Определено разветвленное накрытие ($\ref{py}$) в результате
склейки разветвленных накрытий $\varphi_{Y_j}$ по подпространствам
попарных пересечений семейства подмногообразий $Z_j$ в
многообразии $ZZ$.

Определим вложение $i_J: J \subset \R^n$. Это вложение строится в
результате склейки стандартных вложений торов в семейство
$j_{max}+1$ евклидовых подпространств в $\R^n$ размерности
$n-\frac{n-n_s}{4} + 6$, проходящих через начало координат.
Подпространство с номером $j$ семейства содержит пару пространств
дополнительных размерностей $\R^{n-\frac{n-n_s}{16}- 2j+4}$,
$\R^{2j+2}$, пересекающихся в начале координат. Пара
подпространств с соседними номерами пересекаются по
подпространству коразмерности 2. Пересечение в пространстве с
меньшим (большим) соседним номером теряет коразмерность 2 вдоль
первого (второго) подпространства выбранной пары. Внутри
пространства с номером $j$ семейства рассматривается стандартный
тор $J_j$ с образующими вдоль выбранных дополнительных
подпространств. Объединяя семейство вложенний торов, получим
вложение $i_J$.

\subsubsection*{Конструкция отображения
$d_Y: Y_{a \oplus \dd} \to \R^n$}

Сначала определим вспомогательное отображение $\hat d_Z: Z_{a
\oplus \aa} \to \R^n$. Это отображение определено в результате
малой регулярной $PL$--деформации композиции $i_J \circ  \varphi_Z
: Z_{a \oplus \aa} \to \R^n$, причем сама деформация и ее калибр
$\varepsilon$ выбираются в процессе доказательства (см.
аналогичное построение в [A], Лемма 24). Рассмотрим композицию
$\hat d_Z \circ \varphi_Y: Y_{a \oplus \dd} \to \R^n$ и определим
искомое отображение $d_Y: Y_{a \oplus \dd} \to \R^n$ в результате
малой регулярной $PL$--деформации этой композиции калибра
$\varepsilon'$, причем $\varepsilon' << \varepsilon$.

\subsubsection*{Конструкция отображения
$d_X: X_{a \oplus \dd} \to \R^n$}

Определим вспомогательное отображение $\hat d_Y: Y_{a \oplus \dd}
\to \R^n$. Это отображение определено в результате малой
регулярной $PL$--деформации композиции $i_J \circ  \varphi_Y :
Y_{a \oplus \dd} \to \R^n$, причем сама деформация и ее калибр
$\varepsilon$ выбираются в процессе доказательства (см.
аналогичное построение в [A], Лемма 24). Рассмотрим композицию
$\hat d_Y \circ \varphi_X: X_{b} \to \R^n$ и определим искомое
отображение $d_X: Y_{b} \to \R^n$ в результате малой регулярной
$PL$--деформации этой композиции калибра $\varepsilon'$, причем
$\varepsilon' << \varepsilon$.

\subsubsection*{$\I_a \oplus \II_d$--структура для
отображения $d_X: X_b \to \R^n$}

Рассмотрим полиэдр самопересечения отображения $\hat d_Y: Y_{a
\oplus \aa} \to \R^n$ и обозначим его через $\hat N(d_Y)$. Полиэдр
$\hat N(d_Y)$ является многообразием с особенностями (в
коразмерности 2) с краем, этот край обозначим через $\partial \hat
N(d_Y)$.

Рассмотрим полиэдр самопересечения отображения $d_X: X_{b} \to
\R^n$ и обозначим его через $N(d_X)$. Полиэдр $N(d_X)$ является
многообразием с особенностями (в коразмерности 2) с краем, этот
край обозначим через $\partial N(d_X)$. Многообразие с
особенностями с краем $N(d_X)$ представляется в объединение двух
многообразий с особенностями c краем $N(d_X)=N_{X,antigiag} \cup
N_{X,\Gamma}$ по общему краю, таким образом, что:

1. Многообразие с особенностями с краем $ N_{X,\Gamma}$ является
накрывающим пространством при регулярном 4-листном накрытии
$p_{N_{X,\Gamma}}: N_{X,\Gamma} \to \hat N(d_Y)$.

2. Многообразие с особенностями с краем $ N_{X,antidiag}$
возникает при деформации двулистного накрытия $p_{Y_{a \oplus
\dd}}: X_b \to Y_{a \oplus \dd}$ внутри регулярной (погруженной)
окрестности неособых точек полиэдра $\hat d_Y(Y_{a \oplus \dd})$.

Рассмотрим полиэдр самопересечения отображения $d_X: X_{b} \to
\R^n$ и обозначим его через $N(d_X)$. Полиэдр $N(d_X)$ является
многообразием с особенностями (в коразмерности 2) с краем, этот
край обозначим через $\partial N(d_X)$. Многообразие с
особенностями с краем $N(d_X)$ представляется в объединение двух
многообразий с особенностями c краем $N(d_X)=N_{X,antigiag} \cup
N_{X,\Gamma}$ по общему краю, таким образом, что:

1. Многообразие с особенностями с краем $ N_{X,\Gamma}$ является
накрывающим пространством при регулярном 4-листном накрытии
$p_{N_{X,\Gamma}}: N_{X,\Gamma} \to N(d_Y)$.

2. Многообразие с особенностями с краем $ N_{X,antidiag}$
возникает при деформации двулистного накрытия $p_X$ внутри
регулярной (погруженной) окрестности неособых точек полиэдра
$d_Y(Y_{a \oplus \dd})$.

Повторяя рассуждения из  Леммы 24 [A], определим отображение
\begin{eqnarray}\label{yN}
\zeta_{a \oplus \dd, N(d_X)}: (N(d_X),\partial N(d_X)) \to (K(\I_a
\oplus \II_d,1), K(\J_b,1)).
\end{eqnarray}
Это отображение определено в результате склейки
характеристического отображения на $N_{X,antidiag}$ с
предварительно построенным отображением на $N_{Y,\Gamma}$, по
антидиагональной части границы, где указанный отображения
гомотопны. Граничные условия на $\partial N(d_X)$ определяются
композицией $\partial N(d_X) \subset X_{b}
\stackrel{\eta_{X_b}}{\longrightarrow} K(\J_b,1)$. Отображение
($\ref{yN}$) определяет $\I_a \oplus \II_d$ структуру для
отображения $d_X$.

\subsubsection*{$\I_a \oplus \II_a$--структура для
отображения $d_Y: Y_{a \oplus \dd} \to \R^n$}

Аналогично предыдущему построению, определим многообразие с
особенностями с краем $N(d_Y)$ и отображение
\begin{eqnarray}\label{zN}
\zeta_{a\oplus \aa,N(d_Y)}: (N(d_Y),\partial N(d_Y)) \to (K(\I_a
\oplus \II_a,1), K(\I_a \oplus \II_d,1)).
\end{eqnarray}

Граничные условия на $\partial N(d_Y)$ определяются композицией
$\partial N(d_Y) \subset Y_{a \oplus \dd} \stackrel{\eta_{Y_{a
\oplus \dd}}}{\longrightarrow} K(\I_a \oplus \II_d,1)$.
Отображение ($\ref{zN}$) определяет $\I_a \oplus \II_d$ структуру
для отображения $d_Y$.

\subsubsection*{Конструкция $\D_4$--оснащенного погружения с
$\I_a \oplus \II_d$--структурой в Теореме 11}

Пусть задано $\D_4$--оснащенное погружение $(g,\Psi,\eta_N)$, $g:
N^{n-\frac{n-n_s}{16}} \looparrowright \R^n$. По условию теоремы
задано отображение $\eta_{b,N}: N^{n-\frac{n-n_s}{16}} \to
K(\J_b,1)$, при этом выполнено уравнение ($\ref{etaad'}$).
Отображение $\eta_{b,N}$ определяет однозначно (с точностью до
гомотопии) отображение $\eta_{b,X}: N^{n-\frac{n-n_s}{16}} \to
X_{b}$, поскольку $X_b$ вкладывается в $K(\J_b,1)$ как остов
стандартного клеточного разбиения, который содержит меньший остов
стандартного клеточного разбиения размерности
$n-\frac{n-n_s}{16}+1 = \dim(N)+1$.

Рассмотрим композицию $d_X \circ \eta_{b,X}:
N^{n-\frac{n-n_s}{16}} \to \R^n$ и рассмотрим малую деформацию
этого отображения в погружение $g_1$ в классе регулярной гомотопии
данного погружения $g: N^{n-\frac{n-n_s}{16}} \looparrowright
\R^n$. Калибр $\delta$ деформации $d_X \circ \eta_{b,X} \mapsto
g_1$ выбирается много меньшим $\varepsilon'$. Погружение $g_1$
снабжено $\D_4$--оснащением $\Psi_1$ c тем же характеристическим
классом $\eta_N$, при этом тройка $(g_1, \Psi_1, \eta_N)$
определяет в группе
$Imm^{\D_4}(n-\frac{n-n_s}{16},\frac{n-n_s}{16})$ элемент $y$.

Обозначим через $L^{n-\frac{n-n_s}{8}}$ многообразие
самопересечения погружения $g_1$. Определено разбиение
\begin{eqnarray}\label{LL}
L^{n-\frac{n-n_s}{8}} =  L^{n- \frac{n-n_s}{8}}_{cycl} \cup
L^{n-\frac{n-n_s}{8}}_{b}
\end{eqnarray}
по общей границе. При этом многообразие
$L^{n-\frac{n-n_s}{8}}_{cycl}$ погружено в регулярную
(погруженную) окрестность полиэдра $N(d_X)$. Многообразие
$L^{n-\frac{n-n_s}{8}}_{b}$ погружено в регулярную (погруженную)
регулярных точек полиэдра $d_X(X_b)$ так, что определено
отображение
\begin{eqnarray}\label{projektLb}
\pi_{L_{b}}: L^{n-\frac{n-n_s}{8}}_{b} \to X_b
\end{eqnarray}
проекции этой части многообразия ($\ref{LL}$) на центральный
полиэдр в рассматриваемой окрестности.
 Общая граница
этих многообразий погружена в регулярную (погруженную) окрестность
критических значений отображения  $d_X: X_b \to \R^n$ (cм.
аналогичное Предложение 20 из [A]).

Определим искомое отображение
\begin{eqnarray}\label{Ladd}
\zeta_{a \oplus \dd,L}:  L^{n- \frac{n-n_s}{8}} \to K(\I_a \oplus
\II_d),1).
\end{eqnarray}
 Отображение $\zeta_{a \oplus \dd}$ зададим отдельно на
компонентах разбиения ($\ref{LL}$). На компоненте
$L^{n-\frac{n-n_s}{8}}_{cycl}$ отображение $\zeta_{a \oplus
\dd,L}$ определено отображением  ($\ref{yN}$), которое
продолжается на всю погруженную регулярную окрестность
многообразия с особенностями $N(d_X)$ вне критических значений. На
компоненте $L^{n-\frac{n-n_s}{8}}_{b}$ отображение $\zeta_{a
\oplus \dd,L}$ определено композицией
$$ L^{n-\frac{n-n_s}{8}}_{b}
\stackrel{\pi_{L_{b}}}{\longrightarrow}
 X_b \stackrel{\eta_{X}}{\longrightarrow} K(\J_b,1)
 \stackrel{i_{b,a \oplus \dd}}{\longrightarrow}
  K(\I_a \oplus \II_d,1).$$
На общей границе указанные отображения можно склеить, что следует
из выполнения граничных условий для отображения ($\ref{yN}$).
Отображение ($\ref{Ladd}$), определяющее $\I_a \oplus
\II_d$--структуру $\D_4$--оснащенного погружения
$(g_1,\Psi_1,\eta_N)$ определено.

\subsubsection*{Конструкция $\Z/2^{[3]}$--оснащенного погружения с
бициклической структурой в Теореме  12}

Пусть задано $\Z/2^{[3]}$--оснащенное погружение
$(g,\Psi,\eta_N)$, $g: N^{n-\frac{n-n_s}{8}} \looparrowright
\R^n$, определяющее элемент $y \in
Imm^{\Z/2^{[3]}}(n-\frac{n-n_s}{8},\frac{n-n_s}{8})$. По условию
теоремы задано отображение $\eta_{a \oplus \dd}:
N^{n-\frac{n-n_s}{8}} \to K(\I_a \oplus \II_d,1)$, при этом
выполнено уравнение ($\ref{etaad''}$). Отображение $\eta_{a \oplus
\dd,N}$ определяет однозначно (с точностью до гомотопии)
отображение $\eta_{b,X}: N^{n-\frac{n-n_s}{8}} \to X_{b}$,
поскольку $X_b$ вкладывается в $K(\J_b,1)$ как остов стандартного
клеточного разбиения, который заведомо содержит меньший остов
стандартного клеточного разбиения размерности $n-\frac{n-n_s}{8}+1
= \dim(N)+1$.

Рассмотрим композицию $d_Y \circ \eta_{a \oplus \dd,Y}:
N^{n-\frac{n-n_s}{8}} \to \R^n$ и рассмотрим малую деформацию
этого отображения в погружение $g_1$ в классе регулярной гомотопии
данного погружения $g: N^{n-\frac{n-n_s}{8}} \looparrowright
\R^n$. Калибр $\delta$ деформации $d_Y \circ \eta_{a \oplus \dd,Y}
\mapsto g_1$ выбирается много меньшим $\varepsilon'$. Погружение
$g_1$ снабжено $\Z/3^{[3]}$--оснащением $\Psi_1$ c тем же
характеристическим классом $\eta_N$, при этом тройка $(g_1,
\Psi_1, \eta_N)$ определяет в группе
$Imm^{\Z/2^{[3]}}(n-\frac{n-n_s}{8},\frac{n-n_s}{8})$ элемент $y$.

Обозначим через $L^{n-\frac{n-n_s}{4}}$ многообразие
самопересечения погружения $g_1$. Определено разбиение
\begin{eqnarray}\label{LL'}
L^{n-\frac{n-n_s}{4}} =  L^{n- \frac{n-n_s}{4}}_{cycl} \cup
 L^{n-\frac{n-n_s}{4}}_ {a \oplus \dd}
\end{eqnarray}
по общей границе. При этом многообразие
$L^{n-\frac{n-n_s}{4}}_{cycl}$ погружено в регулярную
(погруженную) окрестность полиэдра $N(d_Y)$. Многообразие
$L^{n-\frac{n-n_s}{4}}_{a \oplus \dd}$ погружено в регулярную
(погруженную) регулярных точек полиэдра $d_Y(Y_{a \oplus \dd})$
так, что определено отображение
\begin{eqnarray}\label{projektLadd}
\pi_{L_{a \oplus \dd}}: L^{n-\frac{n-n_s}{4}}_{a \oplus \dd} \to
Y_{a \oplus \dd}
\end{eqnarray}
проекции этой части многообразия ($\ref{LL'}$) на центральный
полиэдр в рассматриваемой окрестности.
 Общая граница
этих многообразий погружена в регулярную (погруженную) окрестность
критических значений отображения $d_Y: Y_{a \oplus \dd} \to \R^n$
(cм. аналогичное Предложение 20 из [A]).

Определим искомое отображение
\begin{eqnarray}\label{Laaa}
\zeta_{a \oplus \aa,L}: L^{n- \frac{n-n_s}{4}} \to K(\I_a \oplus
\II_a,1).
\end{eqnarray}
 Отображение $\zeta_{a \oplus \aa}$ зададим отдельно на
компонентах разбиения ($\ref{LL'}$). На компоненте
$L^{n-\frac{n-n_s}{4}}_{cycl}$ отображение $\zeta_{a \oplus
\aa,L}$ определено отображением  ($\ref{zN}$), которое
продолжается на всю регулярную окрестность многообразия с
особенностями $N(d_Y)$. На компоненте $L^{n-\frac{n-n_s}{4}}_{a
\oplus \aa}$ отображение $\zeta_{a \oplus \aa, L}$ определено
композицией
$$ L^{n-\frac{n-n_s}{4}}_{a \oplus \dd}
\stackrel{\pi_{L_{a \oplus \dd}}}{\longrightarrow}
 Y_{a \oplus \dd} \stackrel {\eta_{Y}}{\longrightarrow} K(\I_a \oplus \II_d,1)
 \stackrel{i_{a \oplus \dd,a \oplus \aa}}{\longrightarrow}
  K(\I_a \oplus \II_a,1).$$
На общей границе указанные отображения можно склеить, что следует
из выполнения граничных условий для отображения ($\ref{zN}$).
Отображение ($\ref{Laaa}$), определяющее $\I_a \oplus
\II_a$--структуру $\Z/2^{[3]}$--оснащенного погружения
$(g_1,\Psi_1,\eta_N)$ определено.
\[  \]

\subsubsection*{Проверка уравнения ($\ref{etaad'}$)}

Рассмотрим $\D_4$--оснащенное погружение $(g_1,\Psi_1,\eta_N)$,
построенное выше. Это $\D_4$--оснащенное погружение определяет тот
же элемент $y \in
Imm^{\Z/2^{[3]}}(n-\frac{n-n_s}{16},\frac{n-n_s}{16})$. По условию
многообразие $N^{n-\frac{n-n_s}{16}}$ снабжено отображением
$\eta_{b,N}: N^{n-\frac{n-n_s}{16}} \to K(\J_b,1)$. Рассмотрим
многообразие $L^{n-\frac{n-n_s}{8}}$ самопересечения погружения
$g_1$, снабженное отображением ($\ref{Ladd}$).

Определим подмногообразие
\begin{eqnarray}\label{Neta}
N_{\eta}^ {n-\frac{n-n_s}{2}} \subset N^{n-\frac{n-n_s}{16}},
\end{eqnarray}
двойственное в смысле Пуанкаре коциклу $\eta^{\frac{7(n-n_s)}{32}}
\in H^{\frac{7(n-n_s)}{16}}(N^{n-\frac{n-n_s}{16}};\Z/2)$.
Рассмотрим погружение (общего положения) $g_{N_{\eta}}:
N_{\eta}^{n-\frac{n-n_s}{2}} \looparrowright \R^n$, определенное
как ограничение погружения $g_1$ на подмногообразие
($\ref{Neta}$). Обозначим через  $L^{n_s}_{\eta}$ многообразие
самопересечения погружения $g_{N_{\eta}}$. Определено вложение
подмногообразий
\begin{eqnarray}\label{Leta}
L^{n_s}_{\eta} \subset L^{n-\frac{n-n_s}{8}}.
\end{eqnarray}
 Определено
отображение
$$\zeta_{a \oplus \dd, L_{\eta}}: L^{n_s}_{\eta} \to
K(\I_a \oplus \II_d,1)$$ как ограничение отображения
($\ref{Ladd}$) на подмногообразие ($\ref{Leta}$). Заметим, что
подмногообразие ($\ref{Leta}$) представляет гомологический класс,
двойственный в смысле Пуанкаре коциклу
$\zeta^{\frac{7(n-n_s)}{32}} \in
H^{\frac{7(n-n_s)}{8}}(L^{n-\frac{n-n_s}{8}};\Z/2)$.

Поэтому уравнение ($\ref{etaad}$) эквивалентно следующему:
\begin{eqnarray}\label{Larf}
\Theta_{\D_4}^{\frac{n-n_s}{32}} = \langle \bar
  \zeta_{b,L_{\eta}}^{\frac{n_s}{2}};[\bar L_{b, \eta}]
  \rangle.
\end{eqnarray}

Каноническое 2-листное накрывающее $\bar L^{n_s}_{\eta,b}$ над
многообразием ($\ref{Leta}$) естественно погружено  в исходное
многообразие          $N^{n-\frac{n-n_s}{16}}$ и его
фундаментальный цикл представляет (по Теореме Герберта)
гомологический класс, двойственный в смысле Пуанкаре коциклу
$\eta^{\frac{15(n-n_s)}{32}} \in  H^{\frac{15(n-n_s)}{16}}(N^{n-
\frac{n-n_s}{16}};\Z/2)$. Уравнение равнение ($\ref{etaad'}$)
эквивалентно следующему:
\begin{eqnarray}\label{Narf}
\Theta_{\D_4}^{\frac{n-n_s}{32}}(y) = \langle
\eta_{b,N_{\eta}}^{\frac{n_s}{2}};[N_{\eta}]
  \rangle,
\end{eqnarray}
где отображение $\eta_{b,N_{\eta}}:  N_{\eta}^{n_s} \to K(\J_b,1)$
определено в результате ограничения отображения $\eta_{b,N}$ на
подмногообразие ($\ref{Neta}$).

Для доказательства теоремы осталось заметить, что правые части
равенств ($\ref{Larf}$) ($\ref{Narf}$) равны. Это доказано в
следующей лемме. Теорема 11 доказана.

\begin{lemma}

Классы гомологий
\begin{eqnarray}\label{aoplusd}
\bar \zeta_{b, \ast}([\bar L_{b,\eta}] \in
H_{n_s}(K(\J_b,1);\Z/2),
\end{eqnarray}
\begin{eqnarray}\label{etabar}
\eta_{b,\ast}([\bar L_{b,\eta}]) \in H_{n_s}(K(\J_b,1);\Z/2)
\end{eqnarray}
равны.
\end{lemma}

Переформулируем лемму и докажем более общее утверждение. Класс
гомологий ($\ref{etabar}$) можно обобщить и определить в более
общей ситуации, без предположения о том, что многообразие
$L^{n_s}_{\eta}$ является замкнутым.

Пусть определено замкнутое ориентированное многообразие $R^{2r}$ с
особенностями в коразмерности 2 размерности $\dim(R)=2r$,
$\frac{n}{2} \le 2r \le n-\frac{n-n_s}{16}$. Предположим, что
задано отображение $\eta_{b,R}: R^{2r} \to X_b \stackrel{\eta_X,b}
{\longrightarrow} K(\J_b,1)$.

Рассмотрим отображение $g_{R}: R^{2r} \to \R^n$ общего положения,
которое определено в результате малой деформацией общего положения
отображения  $R^{2r}  \stackrel{g_{b,R}} {\longrightarrow} \to X_b
\stackrel{d_X} {\longrightarrow} \R^n$. Определено ориентированное
многообразие с особенностями $T^{n-4r}$ самопересечния отображения
$g_{R}$. Край $\partial T^{n-4r}$ многообразия $T^{n-4r}$ состоит
из критических точек отображения $g_{R}$. Определено каноническое
2-листное накрытие $\bar T^{n-4r} \to T^{n-4r}$, разветвленное
вдоль края.  Определено отображение $\zeta_{b,T}: \bar T^{n-4r}
\to K(\J_b,1)$ в результате композиции погружения $\bar T^{n-4r}
\looparrowright R^{2r}$ с отображением $\eta_{b,R}$.

Определен класс гомологий
\begin{eqnarray}\label{etabarT}
\eta_{b,T,\ast}([\bar T]) \in H_{n-4r}(K(\J_b,1);\Z/2),
\end{eqnarray}
обобщающий класс гомологий ($\ref{etabar}$). Класс гомологий
($\ref{aoplusd}$) также можно обобщить на рассматриваемый случай.

Определено разбиение
\begin{eqnarray}\label{TT}
T^{n-4r} =  T^{n-4r}_{cycl} \cup T^{n-4r}_ {b},
\end{eqnarray}
аналогичное разбиению ($\ref{LL}$).

Край $\partial T^{n-4r}$ многообразия с особенностями $T^{n-4r}$
целиком лежит в компоненте  $T^{n-4r}_ {b}$ и проекция
\begin{eqnarray}\label{projektTb}
\pi_{T_b}: T^{n-4r}_b \to X_b
\end{eqnarray}
аналогичная отображению $(\ref{projektLb})$, переводит край
$\partial T^{n-4r}$ в регулярную часть многообразия с
особенностями $X_b$.

Определено отображение
\begin{eqnarray}\label{Tadd}
\zeta_{a \oplus \dd,T}:  (T^{n-4r}, \partial T^{n-4r}) \to (K(\I_a
\oplus \II_d,1),(K(\J_b,1)),
\end{eqnarray}
аналогичное отображению ($\ref{Ladd}$).

Поскольку коразмерность отображения $g_{R}$ четна, и многообразие
с особенностями $R^{2r}$ является ориентированным, то и
многообразие с особенностями с краем $T^{n-4r}$ также является
ориентированным. Следовательно относительный класс гомологий
$$\zeta_{a \oplus
\dd,T,\ast}: ([T^{n-4r}, \partial T^{n-4r}]) \in H_{n-4r}(K(\I_a
\oplus \II_d,1),(K(\J_b,1));\Z/2)$$ является приведением по модулю
2 соответствующего целочисленного класса из группы
$H_{n-4r}(K(\I_a \oplus \II_d,1),K(\J_b,1));\Z)$.

Ограничение отображения $\zeta_{a \oplus \dd,T}$ на край $\partial
T^{n-4r}$ принимает значения в подпространстве $K(\J_b,1) \subset
K(\I_a \oplus \II_d,1)$. Образ фундаментального цикла $\zeta_{a
\oplus \dd,T,\ast}([\partial T^{n-4r}])$ многообразия с
особенностями, краем, определяет нулевой цикл в группе
$H_{n_s-1}(K(\J_b,1);\Z/2)$, поскольку лежит в образе нулевой
группы $H_{n-4r-1}(K(\J_b,1);\Z)$ при гомоморфизме приведения по
модулю 2 (см. аналогичное утверждение в Лемме 16 из [A]).
Следовательно, класс гомологий
\begin{eqnarray}\label{aoplusdT}
\bar \zeta_{b, \ast}([\bar T_{b}] \in H_{n-4r}(K(\J_b,1);\Z/2)
\end{eqnarray}
определен.

В частности, характеристические классы ($\ref{etabarT}$),
($\ref{aoplusdT}$) определены, если в качестве многообразия
$R^{2r}$ выбрать многообразие $R^{2r}_i=\RP^{2r-2i-1} \times
\RP^{2i+1}$,
 а в качестве
отображения $\eta_{b,R}$ выбрать произвольное стандартное
отображение, определенное как декартово произведение
покоординатных вложениий
\begin{eqnarray}\label{coord}
\eta_{b,R_i}: \RP^{2r-2i-1} \times \RP^{2i+1} \subset X_j
\stackrel{\varphi_{Х_j}}{\longrightarrow} X_b
\stackrel{\eta_{b,X}}{\longrightarrow} K(\J_b,1),
\end{eqnarray}
при произвольных значениях $i,j$, $0 \le 2i \le 2r-2$, $0 \le j
\le j_{max}$.

\begin{lemma}
Обозначим образ фундаментального класса $\eta_{b,R_i,\ast}([R_i])
\in  H_{2r}(K(\J_b,1);\Z)$ через $(2r-2i-1 \times 2i+1)$.

--1.  Для класса гомологий $(2r-2i-1 \times  2i+1) \in
H_{2r}(K(\J_b,1);\Z)$, заданного отображением $(\ref{coord})$
класс гомологий $(\ref{etabarT})$, лежащий в той же группе, равен
$(4r-4i-2-\frac{n}{2} \times 4i+2-\frac{n}{2})$, если оба числа в
скобках строго положительные, и равен нулю, если хотябы одно из
чисел в скобках отрицательно.

--2.  Характеристический класс $(\ref{aoplusdT})$, обобщающий
класс $(\ref{aoplusd})$ для ориентированных многообразий с
особенностями при $2r=n-\frac{n+n_s}{2}$, совпадает с классом
$(\ref{etabarT})$, обобщающий характеристический класс
$(\ref{etabar})$.
\end{lemma}

\subsubsection*{Доказательство Леммы 22}
Утверждение 1 доказывается прямым вычислением, которое опускается.
Утверждение 2 вытекает из построения отображения $d_X$, см.
аналогичную формулу (16) в [A]. Лемма 22 доказана.

\subsubsection*{Доказательство Леммы 21}
Набор стандартных отображений ($\ref{coord}$) при
$2r=n-\frac{n+n_s}{2}$ реализует все элементы в группе
$H_{n-\frac{n-n_s}{2}}(X_b;\Z)$. Поскольку многообразие
$N^{n-\frac{n-n_s}{2}}_{\eta}$ оказывается ориентированным,
отображение $g_{\eta}: N^{n-\frac{n-n_s}{2}}_{\eta} \to \R^n$
кобордантно дизъюнктному набору стандартных отображений
$(\ref{coord})$ в классе отображений ориентированных многообразий
с особенностями в коразмерности 2.  Характеристические классы
отображения при кобордизме отображения сохраняются. Следовательно,
характеристические числа ($\ref{aoplusd}$) ($\ref{etabar}$) равны.
Лемма 21 доказана.

\subsubsection*{Проверка уравнения ($\ref{etaad''}$)}

Рассмотрим $\Z/2^{[3]}$--оснащенное погружение
$(g_1,\Psi_1,\eta_N)$, построенное выше. Это погружение определяет
элемент $y \in Imm^{\Z/2^{[3]}}(n-\frac{n-n_s}{8})$. Многообразие
$N^{n-\frac{n-n_s}{8}}$ снабжено отображением $\eta_{a \oplus
\dd,N}: N^{n-\frac{n-n_s}{8}} \to K(\I_a \oplus \II_d,1)$.
Рассмотрим многообразие $L^{n-\frac{n-n_s}{4}}$ самопересечения
погружения $g_1$, снабженное отображением ($\ref{Laaa}$).

Определим подмногообразие
\begin{eqnarray}\label{Neta'}
N_{\eta}^ {n-\frac{n-n_s}{2}} \subset N^{n-\frac{n-n_s}{8}},
\end{eqnarray}
двойственное в смысле Пуанкаре коциклу $\eta^{\frac{3(n-n_s)}{32}}
\in H^{\frac{3(n-n_s)}{8}}(N^{n-\frac{n-n_s}{8}};\Z/2)$.
Рассмотрим погружение (общего положения) $g_{N_{\eta}}:
N_{\eta}^{n-\frac{n-n_s}{2}} \looparrowright \R^n$, определенное
как ограничение погружения $g_1$ на подмногообразие
($\ref{Neta'}$). Обозначим через  $L^{n_s}_{\eta}$ многообразие
самопересечения погружения $g_{N_{\eta}}$. Определено вложение
подмногообразий
\begin{eqnarray}\label{Leta'}
L^{n_s}_{\eta} \subset L^{n-\frac{n-n_s}{4}}.
\end{eqnarray}
 Определено
отображение
$$\zeta_{a \oplus \aa, L_{\eta}}: L^{n_s}_{\eta} \to
K(\I_a \oplus \II_a,1)$$ как ограничение отображения
($\ref{Laaa}$) на подмногообразие ($\ref{Leta'}$). Заметим, что
подмногообразие ($\ref{Leta'}$) представляет гомологический класс,
двойственный в смысле Пуанкаре коциклу
$\zeta^{\frac{3(n-n_s)}{16}} \in
H^{\frac{3(n-n_s)}{4}}(L^{n-\frac{n-n_s}{4}};\Z/2)$. Каноническое
2-листное накрывающее $\bar L^{n_s}_{\eta}$ над многообразием
($\ref{Leta'}$) естественно погружено  в исходное многообразие
$N^{n-\frac{n-n_s}{8}}$ и его фундаментальный цикл представляет
(по Теореме Герберта) гомологический класс, двойственный в смысле
Пуанкаре коциклу $\eta^{\frac{7(n-n_s)}{32}} \in
H^{\frac{7(n-n_s)}{8}}(N^{n-\frac{n-n_s}{8}};\Z/2)$. Уравнение
($\ref{zetaad}$) эквивалентно следующему:
\begin{eqnarray}\label{Larf'}
\Theta_{\Z/2^{[3]}}^{\frac{n-n_s}{32}}(y) = \langle \bar
  \zeta_{b,L_{\eta}}^{\frac{n_s}{2}};[\bar L_{b, \eta}]
  \rangle.
\end{eqnarray}
Уравнение равнение ($\ref{etaad''}$) эквивалентно следующему:
\begin{eqnarray}\label{Narf'}
\Theta_{\Z/2^{[3]}}^{\frac{n-n_s}{32}}(y) = \langle \bar
\eta_{b,N_{\eta}}^{\frac{n_s}{2}};[N_{b,\eta}]
  \rangle,
\end{eqnarray}
где многообразие $\bar N_{b,\eta}$ и отображение $\bar
\eta_{b,N_{\eta}}: N_{b,\eta} \to K(\J_b,1)$ определены как
2-листное накрывающее над многообразием $N^{\frac{n+n_s}{2}}$ и
отображением $\eta_{a \oplus d,N_{\eta}}: N_{\eta} \to K(\I_a
\oplus \II_d,1)$.

Для доказательства теоремы осталось заметить, что правые части
равенств ($\ref{Larf'}$) ($\ref{Narf'}$) равны. Это доказано в
следующей лемме. Теорема 12 доказана.

\begin{lemma}

Классы гомологий
\begin{eqnarray}\label{aoplusa}
\bar \zeta_{b, L_{\eta},\ast}([\bar L_{b,\eta}] \in
H_{n_s}(K(\J_b,1);\Z/2),
\end{eqnarray}
\begin{eqnarray}\label{etabar'}
\bar \eta_{b,\ast}([\bar L_{b,\eta}]) \in H_{n_s}(K(\J_b,1);\Z/2)
\end{eqnarray}
равны.
\end{lemma}

 Потребуется лемма, аналогичная Лемме 22.

Определим  характеристические классы, аналогичные
($\ref{etabarT}$), ($\ref{aoplusdT}$). Для заданного натурального
$r$, $r=1 \pmod{2}$, определим семейство ориентированных
многообразий $R^{2r}_i=S^{2r-4i+1}/\i \times \RP^{4i-1}$, $1 \le i
\le \frac{r-2}{2}$, $2r \le n-\frac{n-n_s}{32} = \dim(Y_{a \oplus
\dd})-4$. В качестве отображения $\eta_{a \oplus \dd,R}: R^{2r}
\to Y_{a \oplus \dd}$ выберем произвольное стандартное
отображение, определенное как декартово произведение
покоординатных вложениий
\begin{eqnarray}\label{coord'}
\eta_{a \oplus \dd,R_i}: S^{2r-4i+1}/i \times \RP^{4i-1} \subset
Y_j \stackrel{\varphi_{Y_j}}{\longrightarrow} Y_{a \oplus \dd}
\stackrel{\eta_{a \oplus \dd,Y}}{\longrightarrow} K(\I_a \oplus
\II_d,1),
\end{eqnarray}
при некоторых произвольных значениях $i,j$, $0 \le i \le \frac{r-2}{2}$, $0
\le j \le j_{max}$.

Определено отображение $g_R: R^{2r} \to \R^n$, как результат
деформации общего положения отображения $d_Y \circ \eta_{a \oplus
\dd,R}$. Определено ориентированное многообразие с особенностями с
краем $T^{4r-n}$ размерности $\dim(T)=4r-n$, как многообразие
самопересечения отображения $g_{b,R}$. Многообразие с
особенностями  $\bar T^{4r-n}$ служит каноническим 2-листным
разветвленным накрывающим
 над многообразием $T^{4r-n}$. Поскольку многообразие с особенностями
$\bar T^{4r-n}$ погружено в многообразие $R^{2r}$, поэтому
  пределен класс гомологий
\begin{eqnarray}\label{54}
\bar \zeta_{b, T, \ast}([\bar T_{b}] \in
H_{4r-n}(K(\J_b,1);\Z/2),
\end{eqnarray}
обобщающий класс гомологий ($\ref{aoplusa}$).

Определен класс гомологий
\begin{eqnarray}\label{etabarT'}
\eta_{b,T,\ast}([\bar T_b]) \in H_{4r-n}(K(\J_b,1);\Z/2),
\end{eqnarray}
обобщающий класс гомологий ($\ref{etabar'}$).

\begin{lemma}
Обозначим образ фундаментального класса $\eta_{a \oplus
\dd,R_i,\ast}([R_i]) \in  H_{2r}(K(\I_a \oplus \II_d,1);\Z)$
через $(2r-4i+1 \times 4i-1)$.

--1.  Для класса гомологий $(2r-4i+1 \times  4i-1) \in
H_{2r}(K(\I_a \oplus \II_d,1);\Z)$, заданного отображением
$(\ref{coord'})$, класс гомологий $(\ref{etabarT'})$, лежащий в
той же группе, равен   $(4r-8i+2-\frac{n}{2} \times
8i-2-\frac{n}{2})$, если оба числа в скобках строго положительные,
и равен нулю, если хотябы одно из чисел отрицательно.

--2.  Характеристический класс $(\ref{54})$, обобщающий класс
$(\ref{aoplusa})$ для ориентированных многообразий с особенностями
при $2r=n-\frac{n+n_s}{2}$, совпадает с классом
$(\ref{etabarT'})$, обобщающий характеристический класс
$(\ref{etabar'})$.
\end{lemma}

\subsubsection*{Доказательство Леммы 24}
Доказательство аналогично доказательству Леммы 22.
\[  \]

\subsubsection*{Доказательство Леммы 23}

Рассмотрим пару многообразий, определенных формулой
$(\ref{Neta'})$ и переобозначим эту пару многообразий через
\begin{eqnarray}\label{Neta'3}
N_{\eta,[3]}^ {n-\frac{n-n_s}{2}} \subset
N^{n-\frac{n-n_s}{8}}_{[3]}.
\end{eqnarray}
Многообразие $N^{n-\frac{n-n_s}{8}}_{[3]}$ снабжено отображением
\begin{eqnarray}\label{444}
\eta_{a \oplus \dd,[3]}: N^{n-\frac{n-n_s}{8}}_{[3]} \to K(\I_a
\oplus \II_d,1).
\end{eqnarray}
Многообразие $N^{n-\frac{n-n_s}{8}}_{[3]}$ с точностью до
нормального кобордизма совпадает с многообразием
$L^{n-\frac{n-n_s}{8}}$ из Теоремы 11, т.е. является
$\Z/2^{[3]}$--оснащенным многообразием самопересечения
$\D_4$--оснащенного погружения $(g_1, \Psi_1, \eta_{N})$
многообразия $N^{n-\frac{n-n_s}{16}}$, определеного в условии Теоремы
11.

Переобозначим многообразие $N^{n-\frac{n-n_s}{16}}$ через
$N^{n-\frac{n-n_s}{16}}_{[2]}$. Рассмотрим подмногообразие
\begin{eqnarray}\label{45}
N^{n-\frac{n-n_s}{4}}_{\eta,[2]} \subset
N^{n-\frac{n-n_s}{16}}_{[2]},
\end{eqnarray}
двойственное в смысле Пуанкаре когомологическому классу
$\eta_N^{\frac{3(n-n_s)}{32}}$, $\eta_N \in
H^2(N^{n-\frac{n-n_s}{16}}_{[2]};\Z/2)$. Определено погружение
\begin{eqnarray}\label{46}
g_{N_{\eta,[2]}}: N^{n-\frac{n-n_s}{4}}_{\eta,[2]} \looparrowright
\R^n,
\end{eqnarray}
как ограничение погружения $g_1$ на подмногообразие $(\ref{45})$.

Многообразие $N^{n-\frac{n-n_s}{8}}_{[3]}$, определенное в
($\ref{Neta'3}$) является многообразием самопересечения погружения
$g_1$. Подмногообразие   $N_{\eta,[3]}^ {n-\frac{n-n_s}{2}}$,
определенное в ($\ref{Neta'3}$), является многообразием
самопересечения погружения $(\ref{46})$.

Применим Лемму 22 к гомологическому классу
\begin{eqnarray}\label{47}
 \eta_{b,N_{\eta,[2]}}: N_{\eta,[2]} \to X_b \subset K(\J_b,1)
\end{eqnarray}
 и вычислим класс гомологий $(\ref{444})$.  Класс гомологий
 $(\ref{444})$ удовлетворяет условиям Леммы 24, быть может, с точностью до прибавления четного
 гомологического класса. Поскольку
 классы гомологий  $(\ref{aoplusa})$,  $(\ref{etabarT'})$
 лежат в группе, в которой четный элемент равен нулю, не
 ограничивая общности, можно предположить, что гомологический
 класс $(\ref{47})$ определяется линейной комбинацией стандартных
 образующих, описанных в Лемме 24. Лемма 23 доказана.

\section{Доказательство Теорем 17 и 18}

Построения в этом разделе  аналогичны построениям раздела 6.
Начнем доказательство со следующей конструкции. Определим число
$n_s$ из условия Теоремы 6.

Рассмотрим многообразие
  $YY= S^{n-\frac{n-n_s}{4}+5}/\Q_a \times S^{n-\frac{n-n_s}{4}+5}/\i$.
 Это многообразие является прямым произведением однородного кватернионного и линзового пространств.
 Заметим, что $\dim(YY) > n$. Выберем внутри
 многообразия $YY$ подмногообразие $Y$ с особенностями в коразмерности
 4, такое, что многообразие $Y$ вкладывается в
 $\R^n$. В частности, $\dim(Y)<n$.

Рассмотрим в многообразии
 $YY$ семейство подмногообразий
$$Y_j, \quad j=0, \dots, j_{max}, \qquad j_{max}=\frac{3n+n_s+24}{16}$$
размерности  $(n-\frac{n-n_s}{4}+8)$ и коразмерности
$(n-\frac{n-n_s}{4}+2)$, определенное по формуле
$$Y_0 = S^{n-\frac{n-n_s}{4}+5}/\Q_a  \times S^3/\i, \quad Z_1=
S^{n-\frac{n-n_s}{4}+1}/\Q_a \times S^7/\i, \quad \dots,$$
$$ Y_j = S^{n-\frac{n-n_s}{4}+5-4j}/\Q_a \times S^{4j+3}/\i, \quad
\dots, \quad Y_{j_{max}} = S^3/\Q_a \times
S^{n-\frac{n-n_s}{4}+3}/\i.$$ Вложение соответствующего
подмногообразия семейства в многообразие $YY$ определено как
прямое произведение стандартных вложений.

Объединение  $\cup_{j=0} ^{j_{max}} Y_j$ семейства подмногообразий
$\{Y_j\}$ многообразия $YY$ является полиэдром (стратифицированным
подмногообразием с особенностями в коразмерности 4) размерности
$n-\frac{n-n_s}{4}+8$, которое обозначим $Y_{\Q_a \oplus \II_a}
\subset YY$.

Рассмотрим подгруппу
\begin{eqnarray}\label{zepQ}
 \I_a \oplus
\II_a \stackrel {i_{a \oplus \aa, \Q_a \oplus
\II_a}}{\longrightarrow} \Q_a   \oplus \II_a,
\end{eqnarray}
которая индуцирует 2-листное накрытие:
\begin{eqnarray} \label{zzQ}
XX \stackrel{p_{XX,YY}}{\longrightarrow} YY.
\end{eqnarray}

 Определим накрытие
\begin{eqnarray}\label{zQ}
X_{\I_a \oplus \II_a} \stackrel{p_{Y_{\Q_a \oplus
\II_a}}}{\longrightarrow} Y_{\Q_a \oplus \II_a}.
\end{eqnarray}
 Это накрытие индуцировано подгруппой
 ($\ref{zepQ}$) при помощи вложения
$Y_{\Q_a \oplus \II_a} \subset YY$. Таким образом, накрывающее
пространство $X_{\I_a \oplus \II_a}$ в накрытии ($\ref{zQ}$) также
является стратифицированным многообразием с особенностями в
коразмерности 4. Это пространство определено явным образом как
объединение семейства подмногообразий в $XX =
S^{n-\frac{n-n_s}{4}+5}/\i \times S^{n-\frac{n-n_s}{4}+5}/\i$,
определенных по формуле:
$$X_0
= S^{n-\frac{n-n_s}{4}+5}/\i \times S^{3}/\i, \quad \dots, \quad
X_j = S^{n-\frac{n-n_s}{4}+5-4j}/\i \times S^{4j+3}/\i, \dots $$
$$ X_{j_{max}} = S^{3}/\i \times
S^{n-\frac{n-n_s}{4}+5}/\i.$$

Рассмотрим многообразие
  $ZZ= S^{n-\frac{n-n_s}{2}+5}/\Q_a \times S^{n-\frac{n-n_s}{2}+5}/\Q_a$.
 Это многообразие является прямым произведением двух однородных кватернионных пространств.
 Заметим, что $\dim(ZZ) > n$. Выберем внутри
 многообразия $ZZ$ подмногообразие  $Z$ с особенностями в коразмерности
 4, такое, что многообразие $Z$ вкладывается в
 $\R^n$. В частности, $\dim(Z)<n$.

Рассмотрим в многообразии
 $ZZ$ семейство подмногообразий
$$Z_j, \quad j=0, \dots, j'_{max}, \qquad j'_{max}=\frac{n+n_s+6}{2}$$
размерности  $(n-\frac{n-n_s}{2}+8)$ и коразмерности
$(n-\frac{n-n_s}{2}+2)$, определенное по формуле
$$Z_0 = S^{n-\frac{n-n_s}{2}+5}/\Q_a  \times S^3/\Q_a, \quad Z_1=
S^{n-\frac{n-n_s}{2}+1}/\Q_a \times S^7/\Q_a, \quad \dots,$$
$$ Z_j = S^{n-\frac{n-n_s}{2}+5-4j}/\Q_a \times S^{4j+3}/\Q_a, \quad
\dots, \quad Z_{j'_{max}} = S^3/\Q_a \times
S^{n-\frac{n-n_s}{2}+3}/\Q_a.$$ Вложение соответствующего
подмногообразия семейства в многообразие $ZZ$ определено как
прямое произведение стандартных вложений.

Объединение  $\cup_{j=0} ^{j'_{max}} Z_j$ семейства
подмногообразий $\{Z_j\}$ многообразия $ZZ$ является полиэдром
(стратифицированным подмногообразием с особенностями в
коразмерности 4) размерности $n-\frac{n-n_s}{4}+8$, которое
обозначим $Z_{\Q_a \oplus \QQ_a} \subset ZZ$.

Рассмотрим подгруппу
\begin{eqnarray}\label{yQ}
 \Q_a \oplus
\II_a \stackrel {i_{\Q_a \oplus \II_a, \Q_a \oplus
\QQ_a}}{\longrightarrow} \Q_a \oplus \QQ_a,
\end{eqnarray}
которая индуцирует 2-листное накрытие:
\begin{eqnarray} \label{yyQ}
YY' \stackrel{p'_{YY',ZZ}}{\longrightarrow} ZZ.
\end{eqnarray}

 Определим накрытие
\begin{eqnarray}\label{yyyQ}
Y'_{\Q_a \oplus \II_a} \stackrel{p_{Z_{\Q_a \oplus \QQ_a}}}{\longrightarrow} Z_{\Q_a \oplus \QQ_a}.
\end{eqnarray}
 Это накрытие индуцировано подгруппой
 ($\ref{yQ}$) при помощи вложения
$Z_{\Q_a \oplus \II_a} \subset ZZ$. Таким образом, накрывающее
пространство $Y'_{\Q_a \oplus \II_a}$ в накрытии ($\ref{yyyQ}$)
также является стратифицированным многообразием с особенностями в
коразмерности 4. Это пространство определено явным образом как
объединение семейства подмногообразий в $YY' =
S^{n-\frac{n-n_s}{2}+5}/\Q_a \times S^{n-\frac{n-n_s}{2}+5}/\i$,
определенных по формуле:
$$Y'_0
= S^{n-\frac{n-n_s}{2}+5}/\Q_a \times S^{3}/\i, \quad \dots, \quad
Y'_j = S^{n-\frac{n-n_s}{2}+5-4j}/\Q_a \times S^{4j+3}/\i, \dots
$$
$$Y' _{j_{max}} = S^{3}/\Q_a \times
S^{n-\frac{n-n_s}{2}+5}/\i.$$

Определены отображения $\eta_{X}: X_{\I_a \oplus \II_a} \to K(\I_a
\oplus \II_a,1)$, $\eta_{Y}: Y_{\Q_a \oplus \II_a} \to K(\Q_a
\oplus \II_a,1)$, $\eta_{Y'}: Y'_{\Q_a \oplus \II_a} \to K(\Q_a
\oplus \II_a,1)$, $\eta_{Z}: Z_{\Q_a \oplus \QQ_a} \to K(\Q_a
\oplus \QQ_a,1)$, согласованные с включениями подгрупп
($\ref{yQ}$),($\ref{zepQ}$) и накрытиями
($\ref{zQ}$),($\ref{yyyQ}$). Отображение $\eta_X$ представлено
прямым произведением отображений $\eta_{X, a} \times \eta_{X,\aa}:
X_{\I_a \oplus \II_a} \to K(\I_a,1)  \times K(\II_a,1)$.
Отображение $\eta_Y$ представлено прямым произведением отображений
$\eta_{Y,\Q_a} \times \eta_{Y,\aa}: Y_{\Q_a \oplus \II_a} \to
K(\Q_a,1) \times K(\II_a,1)$. Отображение $\eta_{Y'}$ представлено
прямым произведением отображений $\eta_{Y',\Q_a} \times
\eta_{Y',\aa}: Y'_{\Q_a \oplus \II_a} \to K(\Q_a,1) \times
K(\II_a,1)$. Отображение $\zeta_Z$ представлено прямым
произведением отображений $\zeta_{Z,\Q_a} \times \zeta_{Z,\QQ_a}:
Z_{\Q_a \oplus \QQ_a} \to K(\Q_a,1) \times K(\QQ_a,1)$.

Определим многообразие с особенностями $J_{\Q_a \oplus \II_a}$.
Для произвольного
 $j= 0, \dots, j_{max}$ определим пространство
$J_{\Q_a,j}$ как джойн $\frac{3n+n_s-16j+24}{16}$ копий
кватернионного однородного пространства $S^{3}/\Q_a$. При тех же
значениях $j$ определим пространство $J_{\II_a,j}=S^{4j+3}$,
диффеоморфное стандартной $4j+3$--мерной сфере, которое удобно
себе представлять как пространство джойна $j+1$ экземпляров
стандартной сферы $S^3$. Пространство $J_{\Q_a,j} \times
J_{\II_a,j}$ переобозначим для краткости через $J_{\Q_a \oplus
\II_a,j}$.

Определено стандартное включение $i_{J_{\Q_a \oplus \II_a,j}}:
J_{\Q_a \oplus \II_a,j} \subset J_{\Q_a,0} \times J_{\II_a,
j_{max}}$, где сомножитель $J_{\Q_a,j}$ ($J_{\II_a,j}$) включается
в джойн-образ как стандартный подджойн, построенный для
соответствующего числа однородных кватернионных пространств
(трехмерных сфер) при стандартной нумерации. Объединение
$\cup_{j=0}^{j_{max}} Im(i_{J_{\Q_a \oplus \II_a,j}})$ образов
всех этих вложений является искомым пространством, которое
обозначим через $J_{\Q_a \oplus \II_a}$, $J_{\Q_a \oplus \II_a}
\subset J_{\Q_a,j} \times J_{\II_a,j}$ .

Определим многообразие с особенностями $J_{\Q_a \oplus \QQ_a}$.
Для произвольного
 $j= 0, \dots, j'_{max}$ определим пространство
$J'_{\Q_a,j}$ как джойн $\frac{n+n_s-8j+12}{8}$ копий
кватернионного однородного пространства $S^{3}/\Q_a$. При тех же
значениях $j$ определим пространство $J'_{\QQ_a,j}$, диффеоморфное
пространству джойна $j+1$ экземпляра  кватернионного однородного
пространства $S^3/\Q_a$. Пространство $J'_{\Q_a,j} \times
J'_{\QQ_a,j}$ переобозначим для краткости через $J_{\Q_a \oplus
\QQ_a,j}$.

Определено стандартное включение $i_{J_{\Q_a \oplus \QQ_a,j}}:
J_{\Q_a \oplus \Q_a,j} \subset  J'_{\Q_a,0} \times J'_{\QQ_a,
j'_{max}}$, где сомножители $J'_{\Q_a,j}$, $J'_{\QQ_a,j}$
включаются в джойн-образ как стандартные подджойны, построенные
для соответствующего числа однородных кватернионных пространств
при стандартной нумерации. Объединение $\cup_{j=0}^{j'_{max}}
Im(i_{J_{\Q_a \oplus \QQ_a,j}})$ образов всех этих вложений
является искомым пространством, которое обозначим через $J_{\Q_a
\oplus \QQ_a}$, $J_{\Q_a \oplus \QQ_a} \subset J'_{\Q_a,0} \times
J'_{\QQ_a,j'_{max}}$ .

Определим разветвленное накрытие
\begin{eqnarray}\label{varphiyQ}
\varphi_Y: Y_{\Q_a \oplus \II_a} \to J_{\Q_a \oplus \II_a}.
\end{eqnarray}
Для каждого
 $j= 0, \dots, j_{max}$ определим
 отображение
 \begin{eqnarray}\label{varphiYj}
 \varphi_{Y,j}: Y_{\Q_a \oplus \II_a,j} \to J_{\Q_a \oplus \II_a,j}.
 \end{eqnarray}
 Определено стандартное разветвленное накрытие (см.
[A], глава 3, Определение отображения $c$):
\begin{eqnarray}\label{p1Q}
p_{Y_{\Q_a,j},1}: S^{n-\frac{n-n_s}{4}-4j+5}/\Q_a \to
J_{\Q_a,j},
\end{eqnarray}

 Определено стандартное разветвленное накрытие (см.
[A], глава 3, Определение отображения $d$):
\begin{eqnarray}\label{p2Q}
p_{Y_{\II_a,j},1}: S^{4j+3}/\i \to
J_{\II_a,j},
\end{eqnarray}

Отображение ($\ref{varphiYj}$) определено как декартово произведение
описанных выше разветвленных накрытий ($\ref{p1Q}$), ($\ref{p2Q}$) по первому сомножителю и
по второму сомножителю:
$$ \varphi_{Y_{\Q_a \oplus \II_a,j}} = p_{Y_{\Q_a,j},1} \times p_{Y_{\II_a,j},2}:
S^{n-\frac{n-n_s}{4}-4j+5}/\i  \times S^{4j+3}/\i \to $$
$$J_{\Q_a,j} \times J_{\II_a,j} = J_{\Q_a \oplus \II_a,j}. $$

Определено разветвленное накрытие ($\ref{varphiyQ}$) в результате
склейки разветвленных накрытий $\varphi_{Y_{\Q_a \oplus \II_a,j}}$ по подпространствам
попарных пересечений семейства подмногообразий $Y_{\Q_a \oplus \II_a,j}$ в
многообразии $YY$.

Определим разветвленное накрытие
\begin{eqnarray}\label{varphizQ}
\varphi_Z: Z_{\Q_a \oplus \QQ_a} \to J_{\Q_a \oplus \QQ_a}.
\end{eqnarray}
Для каждого
 $j= 0, \dots, j'_{max}$ определим
 отображение
 \begin{eqnarray}\label{varphiZj}
 \varphi_{Z_{\Q_a \oplus \QQ_a,j}}: Z_{\Q_a \oplus \QQ_a,j} \to J_{\Q_a \oplus \QQ_a,j}.
 \end{eqnarray}
 Определены стандартные разветвленные накрытия (см.
[A], глава 3, Определение отображения $c$):
\begin{eqnarray}\label{pz1Q}
p_{Z_{\Q_a,j},1}: S^{n-\frac{n-n_s}{2}-4j+5}/\Q_a \to
J_{\Q_a,j},
\end{eqnarray}

\begin{eqnarray}\label{pz2Q}
p_{Z_{\QQ_a,j},1}: S^{4j+3}/\Q_a \to J_{\QQ_a,j}.
\end{eqnarray}

Отображение ($\ref{varphiZj}$) определено как декартово произведение
описанных выше разветвленных накрытий ($\ref{pz1Q}$), ($\ref{pz2Q}$) по первому сомножителю и
по второму сомножителю:
$$ \varphi_{Z_{\Q_a \oplus \QQ_a},j} = p_{Z_{\Q_a,j},1} \times p_{Z_{\QQ_a,j},2}:
S^{n-\frac{n-n_s}{2}-4j+5}/\Q_a  \times S^{4j+3}/\Q_a \to $$
$$ J_{\Q_a,j} \times J_{\QQ_a,j} = J_{\Q_a \oplus \QQ_a,j}. $$

Определено разветвленное накрытие ($\ref{varphizQ}$) в результате
склейки разветвленных накрытий $\varphi_{Z_{\Q_a \oplus \QQ_a,j}}$ по подпространствам
попарных пересечений семейства подмногообразий $Y_{\Q_a \oplus \II_a,j}$ в
многообразии $ZZ$.

Определим вложение $i_{J_{\Q_a \oplus \II_a}}: J_{\Q_a \oplus \II_a} \subset \R^n$. Это вложение строится в
результате склейки семейства стандартных вложений пространств $J_{\Q_a \oplus \II_a,j}$, $0 \le j \le j_{max}$ в семейство
$j_{max}+1$ евклидовых подпространств в $\R^n$ размерности
$\frac{5(3n+n_s+24)}{16}+3$,  проходящих через начало координат.
Здесь мы воспользовались неравенством $n > 5n_s + 168$, которое выполнено по условию Теоремы 17.
При указанном неравенстве размерность евклидовых пространств рассматриваемого семейства меньше $n$.

Подпространство с номером $j$ семейства содержит пару пространств
дополнительных размерностей: $\R^{\frac{5(3n+n_s+24)}{16}- 5j-1}$,
$\R^{5j+4}$, пересекающихся в начале координат. Пара
подпространств с соседними номерами пересекаются по
подпространству коразмерности 4. Пересечение в пространстве с
меньшим (большим) соседним номером теряет коразмерность 4 вдоль
первого (второго) подпространства выбранной пары. Внутри
пространства с номером $j$ семейства рассматривается декартово произведение
вложений
\begin{eqnarray}\label{ii11}
i_{J_{\Q_a,j}}: J_{\Q_a,j} \subset \R^{\frac{5(3n+n_s+24)}{16}-
5j-1},
\end{eqnarray}
\begin{eqnarray}\label{ii22}
i_{J_{\II_a,j}}: J_{\II_a,j} \subset \R^{5j+4}.
\end{eqnarray}
Вложение ($\ref{ii11}$) определено как джойн
$\frac{3n+n_s+24-16j}{16}$ копий вложений Масси $S^3/\Q_a \subset
\R^4$ (см. [A], раздел 3, Определение отображения $c$). Вложение
($\ref{ii22}$) определено как стандартное вложение
$(4j+3)$--мерной сферы в евклидово пространство.
 Объединяя семейство построенных вложенний, получим
вложение $i_{J_{\Q_a \oplus \II_a}}: J_{\Q_a \oplus \II_a} \subset \R^n$.

Определим вложение $i_{J_{\Q_a \oplus \QQ_a}}: J_{\Q_a \oplus \QQ_a} \subset \R^n$. Это вложение строится в
результате склейки семейства стандартных вложений пространств $J_{\Q_a \oplus \QQ_a,j}$, $0 \le j \le j'_{max}$
в семейство
$j'_{max}+1$ евклидовых подпространств в $\R^n$ размерности
$\frac{5(n+n_s+12)}{8}+3$,  проходящих через начало координат.
Здесь мы воспользовались неравенством $n > \frac{5n_s + 84}{3}$, которое выполнено по условию Теоремы 18.
При указанном неравенстве размерность евклидовых пространств рассматриваемого семейства меньше $n$.

Подпространство с номером $j$ семейства содержит пару пространств
дополнительных размерностей $\R^{\frac{5(n+n_s+12)}{8}- 5j-1}$,
$\R^{5j+4}$, пересекающихся в начале координат. Пара
подпространств с соседними номерами пересекаются по
подпространству коразмерности 4. Пересечение в пространстве с
меньшим (большим) соседним номером теряет коразмерность 4 вдоль
первого (второго) подпространства выбранной пары. Внутри
пространства с номером $j$ семейства рассматривается декартово
произведение вложений
\begin{eqnarray}\label{i11}
i_{J'_{\Q_a,j}}: J'_{\Q_a,j} \subset \R^{\frac{5(n+n_s+12)}{8}-
5j-1},
\end{eqnarray}
\begin{eqnarray}\label{i22}
i_{J'_{\QQ_a,j}}: J'_{\QQ_a,j} \subset \R^{5j+4}.
\end{eqnarray}
Вложение ($\ref{i11}$) определено как джойн
$\frac{n+n_s+12-16j}{8}$ копий вложений Масси $S^3/\Q_a \subset
\R^4$ (см. [A], раздел 3, Определение отображения $c$). Вложение
($\ref{i22}$) определено как джойн $j+1$ копии вложения Масси
$S^3/\Q_a \subset \R^4$.
 Объединяя семейство построенных вложенний, получим
вложение $i_{J_{\Q_a \oplus \QQ_a}}: J_{\Q_a \oplus \QQ_a} \subset
\R^n$.

\subsubsection*{Конструкция отображения
$c_{Y'_{\Q_a \oplus \II_a}}: Y'_{\Q_a \oplus \II_a} \to \R^n$}

Сначала определим вспомогательное отображение $\hat d_{Z_{\Q_a
\oplus \QQ_a}}: Z_{\Q_a  \oplus \QQ_a} \to \R^n$. Это отображение
определено в результате малой регулярной $PL$--деформации
композиции  $i_{J_{\Q_a \oplus \QQ_a}} \circ \varphi_{Z_{\Q_a
\oplus \QQ_a}} : Z_{\Q_a \oplus \QQ_a} \to \R^n$, причем сама
деформация и ее калибр $\varepsilon$ выбираются в процессе
доказательства (см. аналогичное построение в [A], Лемма 25).
Рассмотрим композицию $\hat d_{Z_{\Q_a \oplus \QQ_a}} \circ
p_{Y'_{\Q_a \oplus \II_a}}: Y'_{\Q_a \oplus \II_a} \to \R^n$ и
определим искомое отображение $ c_{Y'_{\Q_a \oplus \II_a}}:
Y'_{\Q_a \oplus \II_a} \to \R^n$ в результате малой регулярной
$PL$--деформации этой композиции калибра $\varepsilon'$, причем
$\varepsilon' << \varepsilon$.

\subsubsection*{Конструкция отображения
$c_{X_{\I_a \oplus \II_a}}: X_{\I_a \oplus \II_a} \to \R^n$}

Определим вспомогательное отображение $\hat d_Y: Y_{\Q_a \oplus
\II_a} \to \R^n$. Это отображение определено в результате малой
регулярной $PL$--деформации композиции $i_{J_{\Q_a \oplus \II_a}}
\circ  \varphi_{Y_{\Q_a \oplus \II_a}} : Y_{\Q_a \oplus \II_a} \to
\R^n$, причем сама деформация и ее калибр $\varepsilon$ выбираются
в процессе доказательства (см. аналогичное построение в [A], Лемма
25). Рассмотрим композицию  $\hat d_{Y_{\Q_a \oplus \II_a}}  \circ
\varphi_{X_{\I_a \oplus \II_a}}:  X_{\I_a \oplus \II_a} \to \R^n$
и определим искомое отображение $d_{X_{\I_a \oplus \II_a}}:
X_{\I_a \oplus \II_a} \to \R^n$ в результате малой регулярной
$PL$--деформации этой композиции калибра $\varepsilon'$, причем
$\varepsilon' << \varepsilon$.

\subsubsection*{$\Q_a \oplus \II_a$--структура для
отображения $c_{X_{\I_a \oplus \II_a}}: X_{\I_a \oplus \II_a} \to
\R^n$}

Рассмотрим полиэдр самопересечения отображения $\hat c_{Y_{\Q_a
\oplus \II_a}}: Y_{\Q_a \oplus \II_a} \to \R^n$ и обозначим его
через $\hat N(d_{Y_{\Q_a \oplus \II_a}})$. Полиэдр $\hat
N(d_{Y_{\Q_a \oplus \II_a}})$ является многообразием с
особенностями (в коразмерности 4) с краем, этот край обозначим
через $\partial \hat N(d_{Y_{\Q_a \oplus \II_a}})$.

Рассмотрим полиэдр самопересечения отображения $c_{X_{\I_a \oplus
\II_a}}: X_{\I_a \oplus \II_a} \to \R^n$ и обозначим его через
$N(c_{X_{\I_a \oplus \II_a}})$. Полиэдр $N(c_{X_{\I_a \oplus
\II_a}})$ является многообразием с особенностями (в коразмерности
4) с краем, этот край обозначим через $\partial N(c_{X_{\I_a
\oplus \II_a}})$. Многообразие с особенностями с краем
$N(c_{X_{\I_a \oplus \II_a}})$ представляется в объединение двух
многообразий с особенностями c краем $N(c_{X_{\I_a \oplus
\II_a}})=N_{X_{\I_a \oplus \II_a},antigiag} \cup N_{X_{\I_a \oplus
\II_a},\Gamma}$ по общему краю, таким образом, что:

1. Многообразие с особенностями с краем $ N_{X_{\I_a \oplus
\II_a,\Gamma}}$ является накрывающим пространством при регулярном
4-листном накрытии $p_{N_{X_{\I_a \oplus \II_a,\Gamma}}}:
N_{X_{\I_a \oplus \II_a,\Gamma}} \to \hat N(c_{Y_{\Q_a \oplus
\II_a}})$.

2. Многообразие с особенностями с краем $ N_{X_{\I_a \oplus
\II_a},antidiag}$ возникает при деформации двулистного накрытия
$p_{Y_{\I_a \oplus \II_a}}: X_{\I_a \oplus \II_a} \to Y_{\Q_a
\oplus \II_a}$ внутри регулярной (погруженной) окрестности
неособых точек полиэдра  $\hat c_{Y_{\Q_a \oplus \II_a}}(Y_{\Q_a
\oplus \II_a})$.

Рассмотрим полиэдр самопересечения отображения $c_{X_{\I_a \oplus
\II_a}}: X_{\I_a \oplus \II_a} \to \R^n$ и обозначим его через
$N(c_{X_{\I_a \oplus \II_a}})$. Полиэдр $N(c_{X_{\I_a \oplus
\II_a}})$ является многообразием с особенностями (в коразмерности
2) с краем, этот край обозначим через $\partial N(c_{X_{\I_a
\oplus \II_a}})$. Многообразие с особенностями с краем
$N(c_{X_{\I_a \oplus \II_a}})$ представляется в объединение двух
многообразий с особенностями c краем $N(c_{X_{\I_a \oplus
\II_a}})=N_{X_{\I_a \oplus \II_a,antigiag}} \cup N_{X_{\I_a \oplus
\II_a,\Gamma}}$ по общему краю, таким образом, что:

1. Многообразие с особенностями с краем $ N_{X_{\I_a \oplus
\II_a},\Gamma}$ является накрывающим пространством при регулярном
4-листном накрытии $p_{N_{X_{\I_a \oplus \II_a},\Gamma}}:
N_{X_{\I_a \oplus \II_a},\Gamma} \to N(c_{Y_{\Q_a \oplus
\II_a}})$.

2. Многообразие с особенностями с краем $ N_{X_{\I_a \oplus
\II_a},antidiag}$ возникает при деформации двулистного накрытия
$p_{X_{\I_a \oplus \II_a}}$ внутри регулярной (погруженной)
окрестности неособых точек полиэдра $c_{Y_{\Q_a \oplus
\II_a}}(Y_{\Q_a \oplus \II_a})$.

Повторяя рассуждения из  Леммы 25 [A], определим отображение
$$\zeta_{\Q_a \oplus \II_a, N(c_{X_{\I_a \oplus \II_a}})}:
(N(c_{X_{\I_a \oplus \II_a}}),\partial N(c_{X_{\I_a \oplus
\II_a}})) \to$$
\begin{eqnarray}\label{yNQ}
(K(\Q_a \oplus \II_a,1), K(\I_a \oplus \II_a,1)).
\end{eqnarray}
Это отображение определено в результате склейки
характеристического отображения на $N_{X_{\I_a \oplus
\II_a},antidiag}$ с предварительно построенным отображением на
$N_{Y_{\Q_a \oplus \II_a},\Gamma}$, по антидиагональной части
границы, где указанный отображения гомотопны. Граничные условия на
$\partial N(c_{X_{\I_a \oplus   \II_a}})$ определяются композицией
$\partial N(c_{X_{\I_a \oplus \II_a}}) \subset X_{\I_a \oplus
\II_a} \stackrel{\eta_{X_{\I_a \oplus \II_a}}}{\longrightarrow}
K(\I_a \oplus \II_a,1)$. Отображение ($\ref{yNQ}$) определяет
$\I_a \oplus \II_a$ структуру для отображения $c_{X_{\I_a \oplus
\II_a}}$.

\subsubsection*{Бикватернионная структура для
отображения $c_{Y'_{\Q_a \oplus \II_a}}: Y'_{\Q_a \oplus \II_a} \to \R^n$}

Аналогично предыдущему построению, определим многообразие с
особенностями с краем $N(c_{Y'_{\Q_a \oplus \II_a}})$ и отображение
$$\zeta_{N(c_{Y'_{\Q_a \oplus \II_a}})}: (N(c_{Y'_{\Q_a \oplus
\II_a}}),\partial N(c_{Y'_{\Q_a \oplus \II_a}}) \to$$
\begin{eqnarray}\label{zNQ}
(K(\Q_a \oplus \QQ_a,1), K(\Q_a \oplus \II_a,1)).
\end{eqnarray}

Граничные условия на $\partial N(c_{Y'_{\Q_a \oplus \II_a}})$ определяются композицией
$\partial N(c_{Y'_{\Q_a \oplus \II_a}}) \subset Y'_{\Q_a \oplus \II_a} \stackrel{\eta_{Y'_{\Q_a
\oplus \II_a}}}{\longrightarrow} K(\Q_a \oplus \II_a,1)$.
Отображение ($\ref{zNQ}$) определяет $\Q_a \oplus \QQ_a$ структуру
для отображения $c_{Y'_{\Q_a \oplus \II_a}}$.

\subsubsection*{Конструкция $\Z^{[4]}$--оснащенного погружения с
$\Q_a \oplus \II_a$--структурой в Теореме 17}

Пусть задано $\Z^{[4]}$--оснащенное погружение $(g,\Psi,\eta_N)$,
$g: N^{n-\frac{n-n_s}{4}} \looparrowright \R^n$. По условию
теоремы задано отображение $\eta_{a \oplus \aa,N}:
N^{n-\frac{n-n_s}{4}} \to K(\I_a \oplus \II_a,1)$, при этом
выполнено уравнение ($\ref{etaQaIa'}$). Отображение $\eta_{a
\oplus \aa,N}$ определяет однозначно (с точностью до гомотопии)
отображение $\eta_{a \oplus \aa,X}: N^{n-\frac{n-n_s}{4}} \to
X_{\I_a \oplus \II_a}$, поскольку $X_{\I_a \oplus \II_a}$
вкладывается в $K(\I_a \oplus \II_a,1)$ как остов стандартного
клеточного разбиения, который содержит меньший остов стандартного
клеточного разбиения размерности $n-\frac{n-n_s}{4}+1 =
\dim(N)+1$.

Рассмотрим композицию $c_{X_{\I_a \oplus \II_a}} \circ \eta_{X_{\I_a \oplus \II_a}}:
N^{n-\frac{n-n_s}{4}} \to \R^n$ и рассмотрим малую деформацию
этого отображения в погружение $g_1$ в классе регулярной гомотопии
данного погружения $g: N^{n-\frac{n-n_s}{4}} \looparrowright
\R^n$. Калибр $\delta$ деформации $c_{X_{\I_a \oplus \II_a}} \circ \eta_{X_{\I_a \oplus \II_a}} \mapsto
g_1$ выбирается много меньшим $\varepsilon'$. Погружение $g_1$
снабжено $\D_4$--оснащением $\Psi_1$ c тем же характеристическим
классом $\eta_N$, при этом тройка $(g_1, \Psi_1, \eta_N)$
определяет в группе
$Imm^{\Z^{[4]}}(n-\frac{n-n_s}{4},\frac{n-n_s}{4})$ элемент $y$.

Обозначим через $L^{n-\frac{n-n_s}{2}}$ многообразие
самопересечения погружения $g_1$. Определено разбиение
\begin{eqnarray}\label{LLQ}
L^{n-\frac{n-n_s}{2}} =  L^{n- \frac{n-n_s}{2}}_{cycl} \cup
L^{n-\frac{n-n_s}{2}}_{\I_a \oplus \II_a}
\end{eqnarray}
по общей границе. При этом многообразие
$L^{n-\frac{n-n_s}{2}}_{cycl}$ погружено в регулярную
(погруженную) окрестность полиэдра $N(c_{X_{\I_a \oplus \II_a}})$. Многообразие
$L^{n-\frac{n-n_s}{2}}_{\I_a \oplus \II_a}$ погружено в регулярную (погруженную)
некритических точек полиэдра $d_{X_{\I_a \oplus \II_a}}(X_{\I_a \oplus \II_a})$ так, что определено
отображение
\begin{eqnarray}\label{projektLbQ}
\pi_{L_{\I_a \oplus \II_a}}: L^{n-\frac{n-n_s}{2}}_{\I_a \oplus \II_a} \to X_{\I_a \oplus \II_a}
\end{eqnarray}
проекции этой части многообразия ($\ref{LLQ}$) на центральный
полиэдр в рассматриваемой окрестности.
 Общая граница
этих многообразий погружена в регулярную (погруженную) окрестность
критических значений отображения  $c_{X_{\I_a \oplus \II_a}}: X_{\I_a \oplus \II_a} \to \R^n$ (cм.
аналогичное Предложение 21 из [A]).

Определим искомое отображение
\begin{eqnarray}\label{LaddQ}
\zeta_{\Q_a \oplus \II_a,L}:  L^{n- \frac{n-n_s}{2}} \to K(\Q_a \oplus
\II_a),1).
\end{eqnarray}
 Отображение $\zeta_{\Q_a \oplus \II_a,L}$ зададим отдельно на
компонентах разбиения ($\ref{LLQ}$). На компоненте
$L^{n-\frac{n-n_s}{2}}_{cycl}$ отображение $\zeta_{\Q_a \oplus
\II_a,L}$ определено отображением  ($\ref{yNQ}$), которое
продолжается на всю погруженную регулярную окрестность
многообразия с особенностями $N(c_{X_{\I_a \oplus \II_a}})$ вне критических значений. На
компоненте $L^{n-\frac{n-n_s}{2}}_{\Q_a \oplus \II_a}$ отображение $\zeta_{\Q_a
\oplus \II_a,L}$ определено композицией
$$ L^{n-\frac{n-n_s}{2}}_{\I_a \oplus \II_a}
\stackrel{\pi_{X_{ \I_a \oplus \II_a}}}{\longrightarrow}
 X_{\I_a \oplus \II_a} \stackrel{\eta_{X_{\I_a \oplus \II_a}}}{\longrightarrow} K(\I_a \oplus \II_a,1)
 \stackrel{i_{\I_a \oplus \II_a,\Q_a \oplus \II_a}}{\longrightarrow}
  K(\Q_a \oplus \II_d,1).$$
На общей границе указанные отображения можно склеить, что следует
из выполнения граничных условий для отображения ($\ref{yNQ}$).
Отображение ($\ref{LaddQ}$), определяющее $\Q_a \oplus
\II_a$--структуру $\Z^{[4]}$--оснащенного погружения
$(g_1,\Psi_1,\eta_N)$ определено.

\subsubsection*{Конструкция $\Z/2^{[5]}$--оснащенного погружения с
бикватернионной структурой в Теореме 18}

Пусть задано $\Z/2^{[5]}$--оснащенное погружение
$(g,\Psi,\eta_N)$, $g: N^{n-\frac{n-n_s}{2}} \looparrowright
\R^n$, определяющее элемент $y \in
Imm^{\Z/2^{[5]}}(n-\frac{n-n_s}{2}, \frac{n-n_s}{2})$. По условию
теоремы задано отображение $\eta_{\Q_a \oplus \II_a}:
N^{n-\frac{n-n_s}{2}} \to K(\Q_a \oplus \II_a,1)$, при этом
выполнено уравнение ($\ref{etaQaQa'}$). Отображение $\eta_{\Q_a
\oplus \II_a,N}$ определяет однозначно (с точностью до гомотопии)
отображение $\eta_{Y'_{\Q_a \oplus \II_a}}: N^{n-\frac{n-n_s}{2}}
\to Y'_{\Q_a \oplus \II_a}$, поскольку $Y'_{\Q_a \oplus \II_a}$
вкладывается в $K(\Q_a \oplus \II_a,1)$ как остов стандартного
клеточного разбиения, который заведомо содержит меньший остов
стандартного клеточного разбиения размерности $n-\frac{n-n_s}{2}+1
= \dim(N)+1$.

Рассмотрим композицию $c_{Y'_{\Q_a \oplus \II_a}} \circ
\eta_{Y'_{\Q_a \oplus \II_a}}: N^{n-\frac{n-n_s}{2}} \to \R^n$ и
рассмотрим малую деформацию этого отображения в погружение $g_1$ в
классе регулярной гомотопии данного погружения $g:
N^{n-\frac{n-n_s}{2}} \looparrowright \R^n$. Калибр $\delta$
деформации $d_{Y'_{\Q_a \oplus \II_a}} \circ \eta_{Y'_{\Q_a \oplus
\II_a}} \mapsto g_1$ выбирается много меньшим $\varepsilon'$.
Погружение $g_1$ снабжено $\Z/3^{[5]}$--оснащением $\Psi_1$ c тем
же характеристическим классом $\eta_N$, при этом тройка $(g_1,
\Psi_1, \eta_N)$ определяет в группе
$Imm^{\Z/2^{[5]}}(n-\frac{n-n_s}{2},\frac{n-n_s}{2})$ элемент $y$.

Обозначим через $L^{n_s}$ многообразие самопересечения погружения
$g_1$. Определено разбиение
\begin{eqnarray}\label{LL'Q}
L^{n_s} =  L^{n_s}_{cycl} \cup
 L^{n_s}_{\Q_a \oplus \II_a}
\end{eqnarray}
по общей границе. При этом многообразие $L^{n_s}_{cycl}$ погружено
в регулярную (погруженную) окрестность полиэдра $N(d_{Y_{\Q_a
\oplus \II_a}})$. Многообразие $L^{n_s}_{\Q_a \oplus \II_a}$
погружено в регулярную (погруженную) регулярных точек полиэдра
$c_{Y'_{\Q_a \oplus \II_a}}(Y'_{\Q_a \oplus \II_a})$ так, что
определено отображение
\begin{eqnarray}\label{projektLaddQ}
\pi_{L_{\Q_a \oplus \II_a}}:   L^{n_s}_{\Q_a \oplus \II_a} \to
Y'_{\Q_a \oplus \II_a}
\end{eqnarray}
проекции этой части многообразия ($\ref{LL'Q}$) на центральный
полиэдр в рассматриваемой окрестности.
 Общая граница
этих многообразий погружена в регулярную (погруженную) окрестность
критических значений отображения $c_{Y'_{\Q_a \oplus \II_a}}:
Y'_{\Q_a \oplus \II_a} \to \R^n$ (cм. аналогичное Предложение 21
из [A]).

Определим искомое отображение
\begin{eqnarray}\label{LaaaQ}
\zeta_{\Q_a \oplus \II_a,L}:  L^{n_s} \to K(\Q_a \oplus \QQ_a,1).
\end{eqnarray}
 Отображение $\zeta_{\Q_a \oplus \QQ_a}$ зададим отдельно на
компонентах разбиения ($\ref{LL'Q}$). На компоненте
$L^{n_s}_{cycl}$ отображение $\zeta_{\Q_a \oplus \QQ_a,L}$
определено отображением  ($\ref{zNQ}$), которое продолжается на
всю регулярную окрестность многообразия с особенностями
$N(d_{Y_{\Q_a \oplus \II_a}})$. На компоненте $L^{n_s}_{\Q_a
\oplus \II_a}$ отображение $\zeta_{\Q_a \oplus \QQ_a, L}$
определено композицией
$$ L^{n_s}_{\Q_a \oplus \II_a}
 \stackrel{\pi_{L_{\Q_a \oplus \II_a}}}{\longrightarrow}
 Y'_{\Q_a \oplus \II_a} \stackrel {\eta_{Y_{\Q_a \oplus \II_a}}}{\longrightarrow} K(\Q_a \oplus \II_a,1)
 \stackrel{i_{\Q_a \oplus \II_a,\Q_a \oplus \QQ_a}}{\longrightarrow}
  K(\Q_a \oplus \QQ_a,1).$$

На общей границе указанные отображения можно склеить, что следует
из выполнения граничных условий для отображения ($\ref{zNQ}$).
Отображение ($\ref{LaaaQ}$), определяющее $\Q_a \oplus
\QQ_a$--структуру $\Z/2^{[5]}$--оснащенного погружения
$(g_1,\Psi_1,\eta_N)$ определено.
\[  \]

\subsubsection*{Проверка уравнения ($\ref{etaQaIa'}$)}

Рассмотрим $\Z/4^{[4]}$--оснащенное погружение
$(g_1,\Psi_1,\eta_N)$, построенное выше. Это погружение определяет
то же элемент $y \in
Imm^{\Z/2^{[4]}}(n-\frac{n-n_s}{4},\frac{n-n_s}{4})$. По условию
многообразие $N^{n-\frac{n-n_s}{4}}$ снабжено отображением
$\eta_{\Q_a \oplus \II_a,N}: N^{n-\frac{n-n_s}{4}} \to K(\Q_a
\oplus \II_a,1)$. Рассмотрим многообразие $L^{n-\frac{n-n_s}{2}}$
самопересечения погружения $g_1$, снабженное отображением
($\ref{LaddQ}$).

Определим подмногообразие
\begin{eqnarray}\label{NetaQ}
N_{\eta}^{n-\frac{n-n_s}{2}} \subset N^{n-\frac{n-n_s}{4}},
\end{eqnarray}
двойственное в смысле Пуанкаре коциклу
$\eta_{[4]}^{\frac{n-n_s}{32}} \in
H^{\frac{n-n_s}{4}}(N^{n-\frac{n-n_s}{4}};\Z/2)$. Рассмотрим
погружение (общего положения) $g_{N_{\eta}}:
N_{\eta}^{n-\frac{n-n_s}{2}} \looparrowright \R^n$, определенное
как ограничение погружения $g_1$ на подмногообразие
($\ref{NetaQ}$). Обозначим через  $L^{n_s}_{\eta}$ многообразие
самопересечения погружения $g_{N_{\eta}}$. Определено вложение
подмногообразий
\begin{eqnarray}\label{LetaQ}
L^{n_s}_{\eta} \subset L^{n-\frac{n-n_s}{2}}.
\end{eqnarray}
 Определено
отображение
$$\zeta_{\Q_a \oplus \II_a, L_{\eta}}: L^{n_s}_{\eta} \to
K(\Q_a \oplus \II_a,1)$$ как ограничение отображения
($\ref{LaddQ}$) на подмногообразие ($\ref{LetaQ}$). Заметим, что
подмногообразие ($\ref{LetaQ}$) представляет гомологический класс,
двойственный в смысле Пуанкаре коциклу $\zeta^{\frac{n-n_s)}{32}}
\in H^{\frac{n-n_s}{2}}(L^{n-\frac{n-n_s}{2}};\Z/2)$.

Поэтому уравнение ($\ref{etaQaIa}$) эквивалентно следующему:
\begin{eqnarray}\label{LarfQ}
\Theta_{\D_4}^{\frac{n-n_s}{32}}(y) = \langle \bar
  \zeta_{b,L_{\eta}}^{\frac{n_s}{2}};[\bar L_{b, \eta}]
  \rangle.
\end{eqnarray}

Каноническое 2-листное накрывающее $\bar L^{n_s}_{\eta,b}$ над
многообразием ($\ref{LetaQ}$) естественно погружено  в исходное
многообразие          $N^{n-\frac{n-n_s}{4}}$ и его
фундаментальный цикл представляет (по Теореме Герберта)
гомологический класс, двойственный в смысле Пуанкаре коциклу
$\eta^{\frac{n-n_s}{32}} \in  H^{\frac{3(n-n_s)}{4}}(N^{n-
\frac{n-n_s}{4}};\Z/2)$. Уравнение равнение ($\ref{etaQaIa'}$)
эквивалентно следующему:
\begin{eqnarray}\label{NarfQ}
\Theta_{\D_4}^{\frac{n-n_s}{32}}(y) = \langle
\eta_{b,N_{\eta}}^{\frac{n_s}{2}};[\bar N_{\eta,b}]
  \rangle,
\end{eqnarray}
где отображение $\eta_{b,N_{\eta}}:  \bar N_{\eta,b}^{n_s} \to
K(\J_b,1)$ определено в результате накрытия над ограничением
отображения $\eta_{b,N}$ на подмногообразие ($\ref{NetaQ}$).

Для доказательства теоремы осталось заметить, что правые части
равенств ($\ref{LarfQ}$) ($\ref{NarfQ}$) равны. Это доказано в
следующей лемме. Теорема 17 доказана.

\begin{lemma}

Классы гомологий
\begin{eqnarray}\label{aoplusdQ}
\bar \zeta_{b, \ast}([\bar L_{b,\eta}] \in
H_{n_s}(K(\J_b,1);\Z/2),
\end{eqnarray}
\begin{eqnarray}\label{etabarQ}
\bar \eta_{b,\ast}([\bar N_{b,\eta}]) \in H_{n_s}(K(\J_b,1);\Z/2),
\end{eqnarray}
аналогичные гомологическим классам $(\ref{aoplusd})$,
$(\ref{etabar})$, равны.
\end{lemma}

Переформулируем лемму и докажем более общее утверждение. Класс
гомологий ($\ref{etabarQ}$) можно обобщить и определить в более
общей ситуации, без предположения о том, что многообразие
$L^{n_s}_{\eta}$ является замкнутым.

Пусть определено замкнутое ориентированное многообразие $R^{2r}$ с
особенностями в коразмерности 4 размерности $\dim(R)=2r$, $2r=-2
\pmod{8}$, $\frac{n}{2} \le 2r \le n-\frac{n-n_s}{4}$.
Предположим, что задано отображение $\eta_{\I_a \oplus \II_a,R}:
R^{2r} \to X_{\I_a \oplus \II_a} \stackrel{\eta_{X_{\I_a \oplus
\II_a}}} {\longrightarrow} K(\I_a \oplus \II_a,1)$.

Рассмотрим отображение $g_{R}: R^{2r} \to \R^n$ общего положения,
которое определено в результате малой деформацией общего положения
отображения  $R^{2r}  \stackrel{g_{\I_a \oplus \II_a,R}}
{\longrightarrow} \to X_{\I_a \oplus \II_a} \stackrel{d_{X_{\I_a
\oplus \II_a}}} {\longrightarrow} \R^n$. Определено
ориентированное многообразие с особенностями $T^{n-4r}$
самопересечния отображения $g_{R}$. Край $\partial T^{n-4r}$
многообразия $T^{n-4r}$ состоит из критических точек отображения
$g_{R}$. Определено каноническое 2-листное накрытие $\bar T^{n-4r}
\to T^{n-4r}$, разветвленное вдоль края. Определено отображение
$\zeta_{\Q_a \oplus \II_a,T}: \bar T^{n-4r} \to K(\Q_a \oplus
\II_a,1)$ в результате композиции погружения $\bar T^{n-4r}
\looparrowright R^{2r}$ с отображением $\eta_{\I_a \oplus
\II_a,R}$.

Определен класс гомологий
\begin{eqnarray}\label{etabarTQ}
\eta_{\I_a \oplus \II_a,T,\ast}([\bar T]) \in H_{n-4r}(K(\I_a
\oplus \II_a,1);\Z/2),
\end{eqnarray}
обобщающий класс гомологий ($\ref{etabarQ}$). Класс гомологий
($\ref{aoplusdQ}$) также можно обобщить на рассматриваемый случай.

Определено разбиение
\begin{eqnarray}\label{TTQ}
T^{n-4r} =  T^{n-4r}_{cycl} \cup T^{n-4r}_ {\I_a \oplus \II_a},
\end{eqnarray}
аналогичное разбиению ($\ref{LLQ}$).

Край $\partial T^{n-4r}$ многообразия с особенностями $T^{n-4r}$
целиком лежит в компоненте  $T^{n-4r}_ {\I_a \oplus \II_a}$ и
проекция
\begin{eqnarray}\label{projektTbQ}
\pi_{T_{\I_a \oplus \II_a}}: T^{n-4r}_{\I_a \oplus \II_a} \to
X_{\I_a \oplus \II_a},
\end{eqnarray}
аналогичная отображению $(\ref{projektLbQ})$, переводит край
$\partial T^{n-4r}$ в регулярную часть многообразия с
особенностями $X_{\I_a \oplus \II_a}$.

Определено отображение
\begin{eqnarray}\label{TaddQ}
\zeta_{\Q_a \oplus \II_a,T}:  (T^{n-4r}, \partial T^{n-4r}) \to
(K(\Q_a \oplus \II_a,1),(K(\I_a \oplus \II_a,1)),
\end{eqnarray}
аналогичное отображению ($\ref{LaddQ}$).

Поскольку коразмерность отображения $g_{R}$ четна, и многообразие
с особенностями $R^{2r}$ является ориентированным, то и
многообразие с особенностями с краем $T^{n-4r}$ также является
ориентированным. Следовательно относительный класс гомологий
$$\zeta_{\Q_a \oplus
\II_a,T,\ast}: ([T^{n-4r}, \partial T^{n-4r}]) \in H_{n-4r}(K(\Q_a
\oplus \II_a,1),(K(\I_a \oplus \II_a,1));\Z/2)$$ является
приведением по модулю 2 соответствующего целочисленного класса из
группы $H_{n-4r}(K(\Q_a \oplus \II_a,1),K(\I_a \oplus
\II_a,1));\Z)$.

Ограничение отображения $\zeta_{\Q_a \oplus \II_a,T}$ на край
$\partial T^{n-4r}$ принимает значения в подпространстве $K(\I_a
\oplus \II_a,1) \subset K(\Q_a \oplus \II_a,1)$. Образ
фундаментального цикла $\zeta_{\Q_a \oplus \II_a,T,\ast}([\partial
T^{n-4r}])$ многообразия с особенностями, краем, определяет
нулевой цикл в группе $H_{n_s-1}(K(\I_a \oplus \II_a,1);\Z/2)$,
поскольку лежит в образе нулевой группы $H_{n-4r-1}(K(\I_a \oplus
\II_a,1);\Z)$ при гомоморфизме приведения по модулю 2 (см.
аналогичное утверждение в Лемме 16 из [A]). Следовательно, класс
гомологий
\begin{eqnarray}\label{aoplusdTQ}
\bar \zeta_{\I_a \oplus \II_a, \ast}([\bar T_{b}] \in
H_{n-4r}(K(\I_a \oplus \II_a,1);\Z/2)
\end{eqnarray}
определен.

В частности, характеристические классы ($\ref{etabarTQ}$),
($\ref{aoplusdTQ}$) определены, если в качестве многообразия
$R^{2r}$ выбрать многообразие $R^{2r}_i=S^{2r-8i+1}/\i \times
S^{8i-1}/\i$,
 а в качестве
отображения $\eta_{\I_a \oplus \II_a,R}$ выбрать произвольное
стандартное отображение, определенное как декартово произведение
покоординатных вложениий
$$\eta_{\I_a \oplus \II_a,R_i}:
S^{2r-8i+1} \times \RP^{8i-1} \subset X_{\I_a \oplus \II_a,j}$$
\begin{eqnarray}\label{coordQ}
\stackrel{\varphi_{Х_{\I_a \oplus \II_a},j}}{\longrightarrow}
X_{\I_a \oplus \II_a} \stackrel{\eta_{X_{\I_a \oplus
\II_a}}}{\longrightarrow} K(\I_a \oplus \II_a,1),
\end{eqnarray}
при произвольных значениях $i,j$, $0 \le 8i \le 2r+1$, $0 \le j
\le j_{max}$.

\begin{lemma}
Обозначим образ фундаментального класса $\eta_{\I_a \oplus
\II_a,R_i,\ast}([R_i]) \in  H_{2r}(K(\I_a \oplus \II_a,1);\Z)$
через $(2r-8i+1 \times 8i-1)$.

-- 1.  Для класса гомологий  $(2r-8i+1 \times  8i-1) \in
H_{2r}(K(\I_a \oplus \II_a,1);\Z)$, заданного отображением
$(\ref{coordQ})$, класс гомологий $(\ref{etabarTQ})$, лежащий в
той же группе, равен $(4r-16i+2-\frac{n}{2} \times
16i-2-\frac{n}{2})$, если оба числа в скобках строго
положительные, и равен нулю, если хотябы одно из чисел в скобках
отрицательно.

-- 2.  Характеристический класс $(\ref{aoplusdTQ})$, обобщающий
класс $(\ref{aoplusdQ})$ для ориентированных многообразий с
особенностями при  $2r=n-\frac{n+n_s}{2}$, совпадает с классом
$(\ref{etabarTQ})$, обобщающий характеристический класс
$(\ref{etabarQ})$.
\end{lemma}

\subsubsection*{Доказательство Леммы 26}
Утверждение 1 доказывается прямым вычислением, которое опускается.
Утверждение 2 вытекает из построения отображения $d_{X_{\I_a
\oplus \II_a}}$, см. аналогичную формулу (20) в [A]. Лемма 26
доказана.

\subsubsection*{Доказательство Леммы 25}

Рассмотрим пару многообразий, определенных формулой
$(\ref{NetaQ})$ и переобозначим эту пару многообразий через
\begin{eqnarray}\label{Neta'3Q}
N_{\eta,[4]}^ {n-\frac{n-n_s}{2}} \subset
N^{n-\frac{n-n_s}{4}}_{[4]}.
\end{eqnarray}
Многообразие $N^{n-\frac{n-n_s}{4}}_{[4]}$ снабжено отображением
\begin{eqnarray}\label{444Q}
\eta_{\I_a \oplus \II_a,[4]}: N^{n-\frac{n-n_s}{4}}_{[4]} \to K(\I_a
\oplus \II_a,1).
\end{eqnarray}
Многообразие $N^{n-\frac{n-n_s}{4}}_{[4]}$ с точностью до
нормального кобордизма совпадает с многообразием
$L^{n-\frac{n-n_s}{4}}$, которое  является
$\Z/2^{[4]}$--оснащенным многообразием самопересечения
$\Z/2^{[3]}$--оснащенного погружения $(g_1, \Psi_1, \eta_{N})$
многообразия $N^{n-\frac{n-n_s}{8}}$, определеного в условии
Теоремы 12.

Переобозначим многообразие $N^{n-\frac{n-n_s}{8}}$ через
$N^{n-\frac{n-n_s}{8}}_{[3]}$. Рассмотрим подмногообразие
\begin{eqnarray}\label{45Q}
N^{n-\frac{n-n_s}{4}}_{\eta,[3]} \subset
N^{n-\frac{n-n_s}{8}}_{[3]},
\end{eqnarray}
двойственное в смысле Пуанкаре когомологическому классу
$\eta_N^{\frac{n-n_s}{16}}$, $\eta_N \in
H^4(N^{n-\frac{n-n_s}{8}}_{[3]};\Z/2)$. Определено погружение
\begin{eqnarray}\label{46Q}
g_{N_{\eta,[3]}}: N^{n-\frac{n-n_s}{4}}_{\eta,[3]} \looparrowright
\R^n,
\end{eqnarray}
как ограничение погружения $g_1$ на подмногообразие $(\ref{45Q})$.

Многообразие $N^{n-\frac{n-n_s}{4}}_{[4]}$, определенное в
($\ref{Neta'3Q}$) является многообразием самопересечения
погружения $g_1$. Подмногообразие   $N_{\eta,[4]}^
{n-\frac{n-n_s}{2}}$, определенное в ($\ref{Neta'3Q}$), является
многообразием самопересечения погружения $(\ref{46Q})$.

Применим Лемму 24 к гомологическому классу
\begin{eqnarray}\label{47Q}
 \eta_{a \oplus \dd,N_{\eta,[3]}}: N_{\eta,[3]}^{n-\frac{n-n_s}{2}} \to X_{\I_a \oplus \II_d} \subset K(\I_a \oplus \II_d,1)
\end{eqnarray}
 и вычислим класс гомологий $(\ref{444Q})$.  Класс гомологий
 $(\ref{444Q})$ удовлетворяет условиям Леммы 26, быть может, с точностью до прибавления четного
 гомологического класса. Поскольку
 классы гомологий    $(\ref{aoplusdQ})$,  $(\ref{etabarQ})$
 лежат в группе, в которой четный элемент равен нулю, не
 ограничивая общности, можно предположить, что гомологический
 класс $(\ref{47Q})$ определяется линейной комбинацией стандартных
 образующих, описанных в Лемме 26. Лемма 25 доказана.

Воспользуемся Леммой 25. Набор стандартных отображений
($\ref{coordQ}$) при $2r=n-\frac{n-n_s}{2}$ реализует все элементы
в группе $H_{n-\frac{n-n_s}{2}}(X_{\I_a \oplus \II_a};\Z)$.
Поскольку многообразие $N^{n-\frac{n-n_s}{2}}_{\eta}$ оказывается
ориентированным, отображение $g_{\eta}:
N^{n-\frac{n-n_s}{2}}_{\eta} \to \R^n$ кобордантно дизъюнктному
набору стандартных отображений $(\ref{coordQ})$ в классе
отображений ориентированных многообразий с особенностями в
коразмерности 4.  Характеристические классы отображения при
кобордизме отображения сохраняются. Следовательно,
характеристические числа ($\ref{aoplusdQ}$) ($\ref{etabarQ}$)
равны. Лемма 25 доказана.

\subsubsection*{Проверка уравнения ($\ref{etaQaQa'}$)}

Рассмотрим $\Z/2^{[5]}$--оснащенное погружение
$(g_1,\Psi_1,\eta_N)$, построенное выше. Это погружение определяет
элемент $y \in
Imm^{\Z/2^{[5]}}(n-\frac{n-n_s}{2},\frac{n-n_s}{2})$, тот же, что
и исходное $\Z/2^{[5]}$--оснащенное погружение $(g, \Psi,
\eta_N)$.

Многообразие $N^{n-\frac{n-n_s}{2}}$ снабжено отображением
$\eta_{\Q_a \oplus \II_a,N}: N^{n-\frac{n-n_s}{2}} \to K(\Q_a \oplus
\II_a,1)$. Рассмотрим многообразие $L^{n_s}$
самопересечения погружения $g_1$, снабженное отображением
($\ref{LaaaQ}$).

Определим подмногообразие
\begin{eqnarray}\label{Neta'Q}
N_{\eta}^ {n_s} \subset N^{n-\frac{n-n_s}{2}},
\end{eqnarray}
двойственное в смысле Пуанкаре коциклу $\eta^{\frac{n-n_s}{8}}
\in H^{\frac{n-n_s}{2}}(N^{n-\frac{n-n_s}{2}};\Z/2)$.
Обозначим через
\begin{eqnarray}\label{LLL}
L^{n_s}_{\eta}
\end{eqnarray}
 многообразие самопересечения погружения
$g_{N_{\eta}}$.
 Определено
отображение
$$\zeta_{\Q_a \oplus \Q_a, L_{\eta}}: L^{n_s}_{\eta} \to
K(\Q_a \oplus \QQ_a,1)$$ как ограничение отображения
($\ref{LaaaQ}$) на подмногообразие ($\ref{LLL}$).

Каноническое 2-листное накрывающее $\bar L^{n_s}_{\eta}$ над
многообразием ($\ref{LLL}$) естественно погружено  в исходное
многообразие $N^{n-\frac{n-n_s}{2}}$ и его фундаментальный цикл
представляет (по Теореме Герберта) гомологический класс,
двойственный в смысле Пуанкаре коциклу $\eta^{\frac{n-n_s}{32}}
\in H^{\frac{n-n_s}{2}}(N^{n-\frac{n-n_s}{2}};\Z/2)$. Уравнение
($\ref{zetQaQa}$) эквивалентно следующему:
\begin{eqnarray}\label{Larf'Q}
\Theta_{\Z/2^{[5]}}^{\frac{n-n_s}{32}}(y) = \langle \bar
  \zeta_{b,L}^{\frac{n_s}{2}};[\bar L_{b}]
  \rangle.
\end{eqnarray}
Уравнение ($\ref{etaQaQa'}$) эквивалентно следующему:
\begin{eqnarray}\label{Narf'Q}
\Theta_{\Z/2^{[5]}}^{\frac{n-n_s}{32}} = \langle \bar
\eta_{b,N}^{\frac{n_s}{2}};[\bar N_{b,\eta}]
  \rangle,
\end{eqnarray}
где многообразие $\bar N_{b,\eta}$ и отображение $\bar
\eta_{b,N_{\eta}}: \bar N_{b,\eta} \to K(\J_b,1)$ определены как
8-листное накрывающее над многообразием $N^{n_s}$ и отображением
$\eta_{\Q_a \oplus \II_a,N_{\eta}}: N_{\eta}^{n_s} \to K(\Q_a
\oplus \II_a,1)$.

Для доказательства теоремы осталось заметить, что правые части
равенств ($\ref{Larf'Q}$) ($\ref{Narf'Q}$) равны. Это доказано в
следующей Лемме. Теорема 18 доказана.

\begin{lemma}

Классы гомологий
\begin{eqnarray}\label{aoplusaQ}
\bar \zeta_{b, L,\ast}([\bar L_{b}] \in
H_{n_s}(K(\J_b,1);\Z/2),
\end{eqnarray}
\begin{eqnarray}\label{etabar'Q}
\bar \eta_{b,\ast}([\bar N_{b,\eta}]) \in H_{n_s}(K(\J_b,1);\Z/2)
\end{eqnarray}
равны.
\end{lemma}

 Потребуется лемма, аналогичная Лемме 22.

Определим  характеристические классы, аналогичные
($\ref{etabarTQ}$), ($\ref{aoplusdTQ}$). Для заданного
натурального $r$, $r=-2 \pmod{16}$ определим семейство
ориентированных многообразий $R^{2r}_i=S^{2r-16i+1}/\Q_a \times
S^{16i-1}/\i$, $0 \le i \le \frac{r-6}{8}$, $2r \le
n-\frac{n-n_s}{2} = \dim(Y'_{\Q_a \oplus \II_a})-8$. В качестве
отображения $\eta_{\Q_a \oplus \II_a,R}: R^{2r} \to Y_{\Q_a \oplus
\II_a}$ выберем произвольное стандартное отображение, определенное
как декартово произведение покоординатных вложениий
$$\eta_{\Q_a \oplus \II_a,R_i}: S^{2r-16i+1}/i \times \RP^{16i-1}
\subset Y_j \stackrel{\varphi_{Y_j}}{\longrightarrow}$$
\begin{eqnarray}\label{coord'Q}
Y_{\Q_a \oplus \II_a} \stackrel{\eta_{Y_{\Q_a \oplus
\II_a}}}{\longrightarrow} K(\Q_a \oplus \II_a,1),
\end{eqnarray}
при некоторых произвольных значениях $i,j$, $0 \le i \le \frac{r-6}{8}$, $0
\le j \le j'_{max}$.

Определено отображение $g_R: R^{2r} \to \R^n$, как результат
деформации общего положения отображения $d_{Y'_{\Q_a \oplus \II_a}} \circ \eta_{\Q_a \oplus
\II_a,R}$. Определено ориентированное многообразие с особенностями с
краем $T^{4r-n}$ размерности $\dim(T)=4r-n$, как многообразие
самопересечения отображения $g_R$. Многообразие с
особенностями  $\bar T^{4r-n}$ служит каноническим 2-листным
разветвленным накрывающим
 над многообразием $T^{4r-n}$. Многообразие с особенностями
$\bar T^{4r-n}$ погружено в многообразие $R^{2r}$, поэтому
  пределен класс гомологий
\begin{eqnarray}\label{54Q}
\bar \zeta_{b, T, \ast}([\bar T_{b}] \in
H_{4r-n}(K(\J_b,1);\Z/2),
\end{eqnarray}
обобщающий класс гомологий ($\ref{aoplusaQ}$).

Определен класс гомологий
\begin{eqnarray}\label{etabarT'Q}
\bar \eta_{b,T,\ast}([\bar T_{b}]) \in H_{4r-n}(K(\J_b,1);\Z/2),
\end{eqnarray}
обобщающий класс гомологий $(\ref{etabar'Q})$.

\begin{lemma}
Обозначим образ фундаментального класса $\eta_{\Q_a \oplus
\II_a,R_i,\ast}([R_i]) \in  H_{2r}(K(\I_a \oplus \II_d,1);\Z)$
через $(2r-16i+1 \times 16i-1)$.

--1.  Для класса гомологий $(2r-16i+1 \times  16i-1) \in
H_{r}(K(\Q_a \oplus \II_a,1);\Z)$, заданного отображением
$(\ref{coord'Q})$, класс гомологий $(\ref{etabarT'Q})$, лежащий в
той же группе, равен   $(4r-32i+2-\frac{n}{2} \times
32i+2-\frac{n}{2})$, если оба числа в скобках строго положительные,
и равен нулю, если хотябы одно из чисел отрицательно.

--2.  Характеристический класс $(\ref{54Q})$, обобщающий класс
$(\ref{aoplusaQ})$ для ориентированных многообразий с особенностями
при $2r=n-\frac{n+n_s}{2}$, совпадает с классом
$(\ref{etabarT'Q})$, обобщающий характеристический класс
$(\ref{etabar'Q})$.
\end{lemma}

\subsubsection*{Доказательство Леммы 28}
Доказательство аналогично доказательству Лемм 23,24.
\[  \]

\subsubsection*{Доказательство Леммы 27}

Рассмотрим многообразие  $(\ref{Neta'Q})$ и переобозначим его
через
\begin{eqnarray}\label{Neta'3QQ}
N_{[5]}^ {n-\frac{n-n_s}{2}}.
\end{eqnarray}
Многообразие ($\ref{Neta'3QQ}$) снабжено отображением
\begin{eqnarray}\label{444QQ}
\eta_{\Q_a \oplus \II_a,[5]}: N^{n-\frac{n-n_s}{2}}_{[5]} \to K(\Q_a
\oplus \II_a,1).
\end{eqnarray}
Многообразие ($\ref{Neta'3QQ}$) с точностью до нормального
кобордизма совпадает с многообразием $L^{n-\frac{n-n_s}{2}}$ из
Теоремы 17, т.е. является  $\Z/2^{[5]}$--оснащенным многообразием
самопересечения $\Z/2^{[4]}$--оснащенного погружения $(g_1,
\Psi_1, \eta_{N})$, где многообразие $N^{n-\frac{n-n_s}{4}}$
определено в условии Теоремы 17. Переобозначим многообразие
$N^{n-\frac{n-n_s}{4}}$ через $N^{n- \frac{n-n_s}{4}}_{[4]}$.

Применим Лемму 28 к гомологическому классу
\begin{eqnarray}\label{47Q}
 \eta_{\Q_a \oplus \II_a,N_{[5]}}: N^{n-\frac{n-n_s}{2}}_{[5]} \to Y'_{\Q_a \oplus \II_a} \subset K(\Q_a \oplus \II_a,1)
\end{eqnarray}
 и вычислим класс гомологий $(\ref{444Q})$.  Класс гомологий
 $(\ref{444Q})$ удовлетворяет условиям Леммы 24, быть может, с точностью до прибавления четного
 гомологического класса. Поскольку
 классы гомологий  $(\ref{54Q})$,  $(\ref{etabarT'Q})$
 лежат в группе, в которой четный элемент равен нулю, не
 ограничивая общности, можно предположить, что гомологический
 класс $(\ref{47Q})$ определяется линейной комбинацией стандартных
 образующих, описанных в Лемме 28. Лемма 27 доказана.

\section{Теорема о ретракции}

Мы докажем следующую теорему о ретракции, которая использовалась в
формулировке Теоремы 6.

 \begin{theorem}{Теорема о ретракции.}
Для произвольного натурального $d$ существует натуральное $l=l(d)$
такое, что произвольный элемент в группе кобордизма
$Imm^{sf}(2^{l'}-3,1)$, при $l' \ge l$, допускает ретракцию
порядка $d-1$ (см. Определение 5).
\end{theorem}

\subsubsection*{Замечание}
Необходимое и достаточное условие существования ретракции
переформулировано ниже в  Следствии 34. Это позволяет
переформулировать Теорему 29 без использования понятия ретракции в
следующей форме. Для произвольного натурального $d$ существует
натуральное $l=l(d)$ такое, что при $l' \ge l$ гомоморфизм $
J_{sf}^d: Imm^{sf}(2^{l'}-3,1) \to Imm^{sf}(d,2^{l'}-2-d)$ равен
нулю.

Пусть $M^{n-k}$--замкнутое $(n-k)$-мерное многообразие,
погруженное в $\R^n$ c коразмерностью $k$ со скошенным оснащением
$\Xi_M$ (характеристический класс этого скошенного оснащения
обозначен через $\kappa_M$). Предположим, кроме того, что
многообразие $M^{n-k}$ снабжено семейством 1-когомологических
классов $A_M=\{\kappa_M, \kappa_1, \dots, \kappa_j\}$, которые
представлены отображениями $\kappa_M$, $\kappa_i: M^{n-k} \to
\RP^{\infty}$, $i=1, \dots, j$. Обозначим набор $\kappa_1, \dots,
\kappa_j$ через $A'_M$.

\begin{definition} Определим группу кобордизма $Imm^{sf;\kappa_1,
\dots, \kappa_j}(n-k;k)$. Представителем этой группы являются
тройки $(M^{n-k},\Xi_M,A_M)$, где
 $M^{n-k}$--замкнутое $n-k$-мерное многообразие со скошенным оснащением, снабженное семейством
1-когомологических классов $A_M$.  Отношение кобордизма является
стандартным отношением эквивалентности на множестве
представителей, которое сохраняет описанную дополнительную
структуру.
\end{definition}

 В
случае $j=0$ группа совпадает с группой скошенно-оснащенных
погруженных многообразий $Imm^{sf}(n-k,k)$.  Если предположить,
что $\kappa_M=0$, причем аналогичное условие выполнено и для
кобордизма, задающего отношение эквивалентности  представителей,
то получим группу оснащенных кобордизмов многообразий, снабженных
семейством когомологических классов $A'_M$. Эта группа
обозначается через $Imm^{fr; \kappa_1, \dots, \kappa_j}(n-k,k)$.
По конструкции Понтрягина-Тома группа $Imm^{fr; \kappa_1, \dots,
\kappa_j}(n-k,k)$ изоморфна прямой сумме $\Pi_{n-k}(\prod_j
\RP^{\infty}) \oplus \Pi_{n-k}$, где первое слагаемое является
стабильной гомотопической группой пространства  $\prod_{i=1}^{i=j}
\RP_i^{\infty}$, а второе слагаемое является $(n-k)$-ой стабильной
гомотопической группой сфер $\Pi_{n-k}$.
\[  \]

Определен естественный (тавталогический) гомоморфизм
\begin{eqnarray}\label{101}
\delta: Imm^{fr; \kappa_1, \dots, \kappa_j}(n-k,k) \to Imm^{sf;
\kappa_1, \dots, \kappa_j}(n-k, k),
\end{eqnarray}
 при этом когомологический класс $\kappa_M$ набора $A_M$
определяется нулевым.

Определен другой естественный гомоморфизм
\begin{eqnarray}\label{102}
J^{k}: Imm^{sf;\kappa_1, \dots, \kappa_j}(n-1,1) \to
Imm^{sf;\kappa_1, \dots, \kappa_j}(n-k,k),
\end{eqnarray}
 посредством перехода от тройки $(M_0^{n-1},
\Xi_{M_0}, A_{M_0})$, представляющей элемент в исходной группе
$Imm^{sf;\kappa_1, \dots, \kappa_j}(n-1,1)$, к тройке $(M^{n-k},
\Xi_M, A_M)$, где подмногообразие $M^{n-k} \subset M_0^{n-1}$
двойственно в смысле Пуанкаре когомологическому классу
$\kappa_{M_0}^{k-1} \in H^{k-1}(M_0;\Z/2)$, скошенное оснащение
$\Xi_M$ определяется стандартным способом, а набор
когомологических классов $A_M$ получен ограничением набора
$A_{M_0}$ на подмногообразие $M^{n-k} \subset M_0^{n-1}$.
\[  \]

Опишем гомоморфизм трансфера $r_j:Imm^{sf;\kappa_1, \dots,
\kappa_j}(n-k,k) \to Imm^{sf;\hat \kappa_1 \dots, \hat
\kappa_{j-1}}(n-k,k)$ относительно когомологического класса
$\kappa_j$. Пусть $\alpha \in Imm^{sf;\kappa_1, \dots,
\kappa_j}(n-k,k)$ -- произвольный элемент, представленный тройкой
$(M^{n-k}, \Xi_M, A_M)$. Рассмотрим двулистное накрытие $p_j:\hat
M^{n-k}_j \to M^{n-k}$, $p_j=\kappa_j + \kappa_M$. Ниже через
$p_j$ мы будем также обозначать линейное расслоение,
соответствующее характеристическому классу $p_j$. Опишем скошенное
оснащение $\Xi_{\hat M_j}$ над $\hat M_j$. Рассмотрим погружение
$\varphi: M^{n-k} \looparrowright \R^n$ коразмерности $k$, класс
регулярной гомотопии которого определен скошенным оснащением
$\Xi_M$.
 Обозначим через
$\nu_{M}$ нормальное расслоение к $\varphi$.
 Над двулистном накрывающем
$\hat M_j$ определено нормальное расслоение $\nu_{\hat M_j}$
погружения $\varphi$, индуцированное из $\nu_{M}$ при проекции
$p_s$. Определим скошенное оснащение $\hat \Xi$ нормального
расслоения $\nu_{\hat M_j}$ как оснащение, индуцированное из
скошенного оснащения $\Xi$ на слои двулистного накрытия.

 Многообразие $\hat M_j$ снабжено
набором когомологических классов $A_{\hat M_j}=\{ \kappa_{\hat
M_j}, \hat \kappa_1, \dots, \hat \kappa_{j-1}\}$, индуцированных
из заданного набора когомологических классов с базы $M^{n-k}$ при
помощи накрытия $p_j$. Построенное скошенно-оснащенное
многообразие $\hat M_j$ c указанным набором классов когомологий
представляет элемент в группе кобордизма $Imm^{sf;\kappa_1 \dots,
\kappa_{j-1}}(n-k,k)$, который определяет образ исходного
элемента.

\begin{example}
Предположим, что  $j=0$. Тогда гомоморфизм трансфера $r_{1}:
Imm^{sf}(n-k,k) \to Imm^{fr}(n-k,k)$ изучен в [A-E].
\end{example}

Обозначим через
\begin{eqnarray}\label{103}
r_{j_1+1, \dots, j_2}: Imm^{sf;\kappa_1, \dots,
\kappa_{j_2}}(n-k,k) \to Imm^{sf;\hat \kappa_1, \dots, \hat
\kappa_{j_1}}(n-k,k)
\end{eqnarray}
композицию последовательных трансферов относительно набора
когомологических классов $A_M=\{\kappa_{j_1}+1, \dots,
\kappa_{j_2} \}$, $0 \le j_1 < j_2$. Более подробно, это означает,
что первый гомоморфизм трансфера берем относительно класса
$\kappa_{j_2} \in H^1(M^{n-k};\Z/2)$, далее второй -- относительно
класса $\hat \kappa_{j_2-1} \in H^1(\hat M_{j_2}^{n-k};\Z/2)$,
который получен индуцированием  $\kappa_{j_2-1}$ на двулистное
накрывающее $\hat M^{n-k}_{j_2}$ накрытия $p_{j_2}$ и т.д.
Нетрудно доказать, что результат применения гомоморфизма
$r_{j_1+1, \dots, j_2}$ не зависит от способа упорядочивания
когомологических классов внутри набора с индексами $j_1+1, \dots,
j_2$. В случае $j_1=0$, гомоморфизм трансфера $r_{1, \dots, j_2}$
строится относительно набора когомологических классов $A'_M$.
Назовем такой гомоморфизмом гомоморфизмом полного трансфера и
обозначим его через $r_{tot}$. Гомоморфизм полного трансфера можно
сразу определить при помощи перехода к $2^{j_2}$--листному
накрытию $\bar M^{n-k}_{tot} \to M^{n-k}$ cо структурной группой
$\Z/2^{j_2}$ при помощи гомоморфизма
$$\oplus_i (\kappa_N +
\kappa_{i}): H^1(M^{n-k};\Z/2) \to \Z/2^{j_2}, \qquad 1 \le i \le
j_2.$$

Накрывающее многообразие $\bar M^{n-k}_{tot}$ в этой конструкции
можно определить как тотальное пространство $\Z^j$--накрытия при
помощи другого упорядоченного набора характеристических классов.
Поменяем местами класс $\kappa_M$ с произвольным классом
$\kappa_{i}$ и определим накрытие, описанным выше способом.
Результат не изменится. Это следует из того, что во-первых,
результат применения трансфера не зависит от того, каким способом
упорядочен набор классов когомологий, относительно которых
строится трансфер. Во-вторых,  на накрывающем пространстве $\bar
M^{n-k}$ над $M^{n-k}$ относительно когомологического класса
$\kappa_{i}+\kappa_M$ индуцированные классы $\bar \kappa_i, \bar
\kappa_M \in H^1(\bar M^{n-k};\Z/2)$ оказываются равными.

\begin{proposition}
Для любого натурального четного числа $k$, $k=0 (mod 2)$,
$2^l-2=n>k>0$, существует такое  натуральное число $\psi$,
зависящее только от $k$, что гомоморфизм полного трансфера
\begin{eqnarray}\label{104}
r_{tot}: Imm^{sf;\kappa_1, \dots, \kappa_{\psi}}(n-k,k) \to
Imm^{sf}(n-k,k)
\end{eqnarray}
тривиальный.
\end{proposition}

\subsubsection*{Доказательство Предложения 32}

Определим  последовательность натуральных чисел $s_{n-k}, \dots,
s_1$, по убыванию индексов от $n-k$ до 1:
$$ s_{n-k} =
ord(\Pi_{n-k}), \quad s_{n-k-1}= ord\Pi_{n-k-1}), \quad \dots,
\quad  s_1= ord(\Pi_1).$$
 Здесь через $ord(\Pi_i)$ обозначена
степень числа 2, которая равна порядку 2-компоненты  $i$-ой
стабильной гомотопической группы сфер. Далее определим
$\psi(i)=\sum_{j=i}^{n-k} s_j$ и, наконец, искомое число $\psi =
\psi(1)+ \sigma +1$, где $\sigma = ord(Imm^{fr;\kappa_1, \dots,
\kappa_{\psi(1)}}(n-k,k))$. Таким образом $\psi > \psi(1) \ge
\psi(2) \ge \dots \ge \psi(n-k)$.

Пусть задан элемент $\alpha \in Imm^{sf;\kappa_1, \dots,
\kappa_{\psi}}(n-k,k)$, представленный тройкой $(M^{n-k}, \Xi_M,
A_M)$. Рассмотрим элемент $r_{\psi}(\alpha)$, представленный
тройкой $(\hat M^{n-k}, \Xi_{\hat M}, A_{\hat M})$, полученный
применением трансфера относительно классов когомологий с номерами
$\psi(1)+1, \dots, \psi$. Докажем, что этот элемент допускает
ретракцию порядка 0.

Рассмотрим произведение $\RP^{\infty}_0 \times
\prod_{i=1}^{i=\psi(j)} \RP^{\infty}(i)$, которое обозначим через
$X(\psi(j))$ и рассмотрим отображение $\lambda: M^{n-k} \to
X(\psi(1))$, покоординатно определенное набором $A_M$
когомологических классов, за исключением классов
$\kappa_{\psi(1)}+1, \dots, \kappa_{\psi}$. Переобозначим
пространство $X(\psi(1))$ для краткости через $X$. Рассмотрим
естественную фильтрацию
\begin{eqnarray}\label{105}
\dots \subset X^{(n-k+1)} \subset X^{(n-k)} \subset \dots \subset
X^{(1)} \subset X,
\end{eqnarray}
которая определена как произведение стандартных фильтраций по
каждой координате. Каждый страт $X^{(i)} \setminus X^{(i+1)}$,
$i=1, \dots n-k$, распадается в объединение открытых клеток
коразмерности $i$. Клетки определяются мультииндексами $\mu=(m_1,
\dots, m_{\psi(n-k)})$, $m_1+ \dots + m_{\psi(n-k)}=i$, $m_i \ge
0$, указывающими коразмерность минимального остова
соответствующего координатного проективного пространства, в
котором содержится данная клетка.


Предположим, что отображение $\lambda$ находится в общем положении
по отношению к фильтрации $(\ref{105})$. Обозначим через $L^0
\subset M^{n-k}$-- $0$-мерное подмногообразие, определенное как
прообраз страта $X^{n-k}$ фильтрации коразмерности $(n-k) =
\dim(M^{n-k})$.

Рассмотрим $2^{\psi-\psi(1)}$--листное накрытие $p: \hat M^{n-k}
\to M^{n-k}$, индуцированное набором когомологических классов
$\kappa_M + \kappa_{\psi(1)+1}, \dots, \kappa_M + \kappa_{\psi}$.
На многообразии $\hat M^{n-k}$ определено скошенное оснащение
$\Xi_{\hat M}$ и набор когомологических классов $A_{\hat M} =
\kappa_{\hat M}, \dots, \hat \kappa_{\psi(1)})$, индуцированный из
поднабора $A_{M,1} = (\kappa_M, \dots \kappa_{\psi(1)})$ набора
$A_M$ на накрывающее многообразие $\hat M$ накрытии $p: \hat
M^{n-k} \to M^{n-k}$. Рассмотрим тройку $(\hat M^{n-k}, \Xi_{\hat
M}, A_{\hat M})$. Эта тройка представляет трансфер  тройки
$(M^{n-k}, \Xi_M, A_M)$ относительно набора когомологических
классов $\kappa_{\psi(1)+1}, \dots, \kappa_{\psi}$. Докажем, что
эта тройка $(\hat M^{n-k}, \Xi_{\hat M}, A_{\hat M})$ в классе
нормального кобордизма допускает отображение $\hat \lambda: \hat
M^{n-k} \to X$, определенное набором характеристических классов
$A_{\hat M}$, для которого прообраз $\hat
\lambda_1^{-1}(X^{(n-k)})$ является пустым.

По предположению $k$ четно. Следовательно, поскольку $n$ четно,
многообразие $M^{n-k}$ является ориентированным.  Обозначим
целочисленный коэффициент пересечения образа $ \lambda(M^{n-k})$ с
ориентированной клеткой из $X^{(n-k)} \setminus X^{(n-k+1)}$
номера $\mu$ через $lk(\mu)$. Аналогично определен также
целочисленный набор $\{\hat lk(\mu)\}$
 для отображения $\hat
\lambda$. Очевидно, что набор $\{\hat lk(\mu)\}$ получен из набора
$\{lk(\mu)\}$ умножением на $2^{\psi-\psi(1)}$.

Рассмотрим двулистное накрытие $\tilde p: \tilde M^{n-k} \to
M^{n-k}$, индуцированное классом когомологий $\kappa_M$ и
рассмотрим композицию $\tilde \lambda: \tilde M^{n-k}
\stackrel{\tilde p}{\longrightarrow} M^{n-k}
\stackrel{\lambda}{\longrightarrow} X$. Многообразие $\tilde
M^{n-k}$ является оснащенным. Обозначим через $\Psi_{\tilde M}$
этого многообразия, индуцированное из скошенного оснащения $\Xi_M$
при накрытии $\tilde p$. Обозначим через $A_{\tilde M}$ набор
когомологических классов $\{\kappa_{\tilde M}, \tilde \kappa_1,
\dots, \tilde \kappa_{\psi(1)}\}$, индуцированных на 2-листное
накрывающее $\tilde M^{n-k}$ при накрытии $\tilde p$. Тройка
$(\tilde M^{n-k}, \Psi_{\tilde M}, A_{\tilde M})$ определяет
элемент из группы $Imm^{fr;\kappa_1, \kappa_{\psi(1)}}(n-k,k)$. По
определению порядок рассматриваемой группы равен $\psi - \psi(1)
-1 $. Следовательно,  $\psi - \psi(1) -1$ экземпляров тройки
$(\tilde M^{n-k}, \Xi_{\tilde M}, A_{\tilde M})$ нормально
ограничивают.

Рассмотрим отображение $\tilde \lambda: \tilde M^{n-k} \to X$.
Набор коэффициентов $\{lk(\tilde \lambda)\}$ пересечения образа
фундаментального класса $\tilde \lambda_{\ast}([\tilde M^{n-k}])$
с клетками пространства $X^{(k)} \setminus X^{(k+1)}$ равен
удвоенному набору $lk(\mu)$, вычисленному по отображению
$\lambda$.

Рассмотрим тройку $(2^{\sigma})(-\tilde M^{n-k}, -\Psi_{\tilde M},
A_{\tilde M})$, которая определена в результате дизъюнктного
объединения $2^{\sigma}$, $(\sigma = \psi - \psi(1)-1)$
экземпляров тройки $(-\tilde M^{n-k},   -\Psi_{\tilde M},
A_{\tilde M})$. Здесь оснащение $\Psi_{\tilde M}$ заменено на
обратное $-\Psi_{\tilde M}$ путем изменения ориентации первого
базисного вектора оснащения (в частности, ориентация многообразия
$\tilde M$ при этом также меняется). Набор коэффициентов
пересечения для отображения $2^{\sigma} \tilde \lambda$,
построенного по набору когомологических классов
$(2^{\sigma})A_{\tilde M}$, обозначим через $\tilde lk(\mu)$.

Рассмотрим тройку $(\hat M'^{n-k}, \Xi_{\hat M'}, A_{\hat M'})$,
которая определена как дизъюнктное объединение тройки $(\hat
M^{n-k}, \Xi_{\hat M}, A_{\hat M})$ с тройкой $(2^{\sigma})(\tilde
M^{n-k}, -\Psi_{\tilde N}, A_{\tilde N})$. Новая тройка
представляет тот же элемент в группе $Imm^{sf;\kappa_1, \dots,
\kappa_{\psi(1)}}(k,n-k)$, что и исходная тройка $(\hat M^{n-k},
\Xi_{\hat M}, A_{\hat M})$. Определено отображение $\hat \lambda'
: \hat M'^k \to X$. При этом набор коэффициентов $\hat lk'(\mu)$,
построенный для новой тройки, равен нулю, поскольку наборы
$\{(2^{\sigma})\tilde lk(\mu)\}$, $\{\hat lk(\mu)\}$ отличаются
знаком.

Определена нормальная перестройка тройки $(\hat M'^{n-k},
\Xi_{\hat M'}, A_{\hat M'})$ по семейству ручек индекса 1 при
которой $\hat \lambda'^{-1}(X^{(n-k)}) = \emptyset$. Начальный
этап индуктивной конструкции определен.

Переобозначим тройку $(\hat M'^{n-k}, \Xi_{M'}, A_{\hat M'})$
снова через $(M^{n-k}, \Xi_M, A_M)$. Отображение, построенное по
набору когомологических классов $A_M$ обозначим снова через
$\lambda: \hat M^{n-k} \to X$. По построению
$\lambda^{-1}(X^{(n-k)})=\emptyset$. Применим к рассматриваемому
элементу гомоморфизм трансфера $r_{\kappa_{\psi(2)+1}, \dots,
\kappa_{\psi(1)}}$. В результате получим тройку, обозначаемую для
краткости снова через $(\hat M^{n-k}, \Xi_{\hat M}, A_{\hat M})$.

Многообразие  $\hat M^{n-k}$ определено как накрывающее
пространство $2^{\psi(1)-\psi(2)}$--листного накрытия $p: \hat
M^{n-k} \to M^{n-k}$ относительно классов когомологий
$(\kappa_{\psi(2)+1}, \dots, \kappa_{\psi(1)})$. Скошенное
оснащение $\Xi_{\hat M}$ индуцировано из оснащения $\Xi_M$ при
накрытии $p$. Набор когомологических классов  $A_{\hat
M}=\{\kappa_{\hat M}, \dots, \hat \kappa_{\psi(2)}\}$ на
накрывающем многообразии $\hat M^{n-k}$ индуцирован из набора
$(\kappa_M, \dots, \kappa_{\psi(2)})$ при накрытии $p$.

Рассмотрим пространство $X(\psi(2))$, определенное как
произведение проективных пространств с индексами $(0,1, \dots,
\psi(2))$. Определено вложение $i_{\psi(2)}: X(\psi(2)) \subset
X(\psi(1))$, соответствующее указанному набору координат.
Пространство $X(\psi(2))$ само снабжено стратификацией, которая
обозначается через
\begin{eqnarray}\label{15}
\dots \subset X^{i}(\psi(2)) \subset X^{i-1}(\psi(2)) \subset
\dots \subset X(\psi(2)).
\end{eqnarray}
 Вложение $i_{\psi(2)}$
согласовано со стратификациями $(\ref{105})$, $(\ref{15})$.

Набор $A_{\hat M}$ когомологических классов определяет отображение
$\hat \lambda: \hat M^{n-k} \to X(\psi(2))$. Поскольку $
\lambda^{-1}(X^{(n-k)}(\psi(1))) = \emptyset$, то и $\hat
\lambda^{-1}(X^{(n-k)}(\psi(2))) = \emptyset$.

 Обозначим через $\hat L^1
\subset \hat M^{n-k}$ одномерное подмногообразие, определенное по
формуле $\hat L^1= \hat \lambda^{-1}(X^{(n-k-1)}(\psi(2)))$.
Ограничение когомологических классов набора $A_{\hat M}$ на
подмногообразие $\hat L^1$ является тривиальным. В частности, это
многообразие $\hat L^1$ является оснащенным. Компоненты
многообразия $\hat L^1$ индексированы в соответствии с клетками
старшей размерности пространства $X^{n-k-1}(\psi(2))$. При этом
компонента, соответствующая фиксированному мультииндексу,
распадается в дизъюнктное объединение $2^{s_1}$ оснащенных кривых
(возможно, несвязных), которые диффеоморфны между собой как
оснащенные 1-мерные многообразия.

Определено оснащенное 2-мерное многообразие $\tilde K^2$ с
оснащенной границей $\partial (\tilde K^2)= (\tilde L^1)$.
Оснащенное многообразие $(\tilde K^2, \Psi_K)$,  определяет
образующую ручки для нормальной перестройки, при которой класс
кобордизма тройки $(\hat M^{n-k}, \Xi_{\hat M}, A_{\hat M})$
сохраняется. В результате нормального кобордизма получим другую
тройку, обозначаемую через  $(\hat M'^{n-k}, \Xi_{\hat M'},
A_{\hat M'})$. Набор когомологических классов $A_{\hat M'}$
определяет отображение $\hat \lambda': \hat M'^{n-k} \to
X(\psi(2))$ и при этом $\hat
\lambda'^{-1}(X^{(n-k-1)}(\psi(2)))=\emptyset$.

Рассуждаем по индукции и  переходим к следующему этапу, при помощи
$2^{s_2}$--листного накрытия относительно набора классов
когомологий $(\hat \kappa_{\psi(3)+1}, \dots, \hat
\kappa_{\psi(2)})$. Изменяем параметр $j$ от 2 до $n-k$. На $j$-ом
шаге получаем тройку $(\tilde M_{j}^{n-k}, \Xi_{\tilde M_{j}},
A_{\tilde M_{j}})$, при этом отображение $\hat \lambda_{j}: \hat
M^{n-k}_{j} \to X(\psi(j))$, построенное по набору $A_{\tilde
M_{j}}$, удовлетворяет условиям: $\hat \lambda_{j}^{-1}
(X^{(n-k-j+1)}(\psi(j)))=\emptyset$, а многообразие $\hat
\lambda_{j}^{-1}(X^{(n-k-j)} (\psi(j)))$ является оснащенным
$j$-мерным многообразием, которое оснащенно ограничивает.

Переход к последующему этапу осуществляется посредством оснащенной
перестройки и переходом к       $2^{s_{j}}$-листному накрытию
относительно классов когомологий $(\hat \kappa_{\psi(j-1)+1},
\dots, \kappa_{\psi(j)})$.

На заключительном шаге конструкции  получим элемент
$r_{tot}(\alpha)$ (см. $(\ref{104})$), который является границей в
группе $Imm^{sf}(n-k,k)$. Предложение 32 доказано.
\[  \]

Опишем полное алгебраическое препятствие для построения ретракции
заданного порядка.

\begin{lemma}
Произвольный элемент $x$ в группе кобордизмов $Imm^{sf,\kappa_1,
\dots, \kappa_j}(n-k,k)$ допускает ретракцию порядка $i$, тогда и
только тогда, когда $J_{sf}^{k'}(x) \in Imm^{sf,\kappa_1, \dots,
\kappa_j}(n-k',k')$  (в предположении $i \le n-k' \le n-k$ )
допускает ретракцию того же порядка $i$.
\end{lemma}

\begin{corollary}
Для произвольного элемента $x \in Imm^{sf;\kappa_1, \dots,
\kappa_j}(n-k,k)$ полное препятствие к построению ретракции
порядка $q$ (в предположении $0 \le q \le n-k $) определено
элементом $J_{sf}^{n-q-1}(x) \in Imm^{sf;\kappa_1, \dots,
\kappa_j}(q+1,n-q-1)$.
\end{corollary}

Для доказательства Леммы 33 докажем еще одну лемму. Предположим,
что скошенно-оснащенное  многообразие
 $(M^{n-k},
\Xi_M)$, снабженное семейством $A_M=(\kappa_M, \kappa_1, \dots ,
\kappa_j)$ когомологических классов, представляет элемент в группе
кобордизмов
 $Imm^{sf;\kappa_1,
\dots, \kappa_j}(n-k,k)$. Пусть $\kappa'_M : M^{n-k} \to
\RP^{n-k-i}$, $i < n-k$ -- ретракция порядка   $(i-1)$
относительно характеристического класса $\kappa_M$ скошенного
оснащения $\Xi_M$. Рассмотрим подмногообразие $Q^i \subset
M^{n-k}$, определенное по формуле $Q^i = \kappa_M'^{-1}(pt), pt
\in \RP^{k-i}$. Многообразие $Q^i$ снабжено естественным
оснащением $\Xi_Q$ с нулевым характеристическим классом, поскольку
подмногообразие $Q^i \subset M^{n-k}$ является оснащенным
подмногообразием и поскольку ограничение скошенного оснащения
$\Xi_M \vert_Q$ является оснащением. Кроме того, определен набор
когомологических классов $A_Q$, полученный ограничением
когомологических классов набора $A_M$ на подмногообразие $Q^i
\subset M^{n-k}$. Класс $\kappa_Q=\kappa_M \vert _Q$ при этом
оказывается тривиальным. Тройка $(Q^i,\Xi_Q, A_Q)$ представляет
элемент в группе оснащенных кобордизмов $Imm^{fr;\kappa_1, \dots,
\kappa_j}(i,n-i)$, который мы обозначим через $y$.

\begin{lemma}
Тройка  $(M^{n-k}, \Xi_M, A_M)$, определяющая элемент $x$ в группе
кобордизма $Imm^{sf,\kappa_1, \dots, \kappa_j} (n-k,k)$ при
условиях, описанных выше, допускает ретракцию порядка $i$
относительно когомологического класса
 $\kappa_M$
 тогда и только тогда, когда элемент $y \in Imm^{fr,\kappa_1, \dots, \kappa_j}(i,n-i)$
 лежит в ядре
гомоморфизма $ \delta: Imm^{fr;\kappa_1, \dots, \kappa_j}(i,n-i)
\to Imm^{sf;\kappa_1, \dots, \kappa_j}(i,n-i)$.
\end{lemma}

\subsubsection*{Доказательство Леммы 35}

В начале докажем полноту препятствия к построению ретракции. По
условию для элемента $y \in Imm^{sf;\kappa_1, \dots,
\kappa_j}(i,n-i)$ имеем $\delta(y)=0$. Рассмотрим
скошенно-оснащенное многообразие $(P^{i+1},\Xi_P)$, $\partial
P^{i+1} = Q^i$, снабженное набором когомологических классов $A_P$,
причем ограничение $A_P \vert_{Q}$ совпадает с $A_Q$. Построим
многообразие $T^{n-k}$, определенное в результате перестройки
многообразия $M^{n-k}$ по ручке, которая описывается ниже.  Ручка
определена как сферическая граница дискового подоснащения в
дисковом $(n-k-i)$-мерном оснащении $\Xi_P$ над многообразием
$P^{i+1}$. Центральное подмногообразие  ручки имеет размерность
$(n-k-i)$, при этом размерность самой ручки совпадает с
размерностью $n-k$ многообразия.

Опишем ручку подробнее. Рассмотрим подрасслоение в нормальном
(оснащенном) расслоении над $P^{i+1}$, порожденное первыми $(k-i)$
векторами оснащения $\Xi_P$ в каждом слое и обозначим через $U_P$
подпространство дискового расслоения, ассоциированного с данным
векторным. Пространство $U_P$ является телом рассматриваемой
ручки. Это пространство является скошенно-оснащенным в
коразмерности $(n-k)$ многообразием с краем.

Граница $\partial U_P$ содержит подмногообразие $Q^i \times
D^{n-k-i}$, которое представляет дисковое расслоение над
многообразием $Q^i$. Вторая копия многообразия $Q^i \times
D^{n-k-i}$ вложена в  компоненту $M^{n-k} \times \{1\}$ границы
многообразия $M^{n-k} \times I$ и представляет регулярную дисковую
окрестность подмногообразия $Q^i \times \{1\} \subset M^{n-k}
\times \{1\}$. Определим многообразие $T^{n-k}$ по формуле
$$  T^{n-k}= \partial^+((M^{n-k} \times I) \cup_{Q^i \times D^{n-k-i}}  U_P), $$
где через $\partial^+$ обозначена компонента границы,
соответствующая перестроенной компоненте $M^{n-k} \times \{1\}$
границы $\partial(M^{n-k} \times I)$.

Заметим, что после стандартной операции "сглаживания углов"
$\quad$ многообразие $T^{n-k}$ становится гладким (замкнутым). Это
многообразие допускает каноническое скошенное оснащение
 $\Xi_T$ с характеристическим классом
$\kappa_T$. Оснащение $\Xi_T$ размерности $k$ получено в
результате склейки соответствующих скошенных оснащений к ручке и к
многообразию $\Xi_M$ вдоль общей части $Q^i$. Многообразие
$T^{n-k}$ снабжено набором когомологических классов
$A_T=\{\kappa_T,\kappa_1, \dots, \kappa_j\} \in
H^1(T^{n-k};\Z/2)$, каждый класс из этого набора получен в
результате склейки соответствующих когомологических классов на
многообразии $M^{n-k} \times I$ и на ручке. Тройка
 $(T^{n-k},\Xi_T, A_T)$
 представляет класс кобордизма исходного элемента $(M^{n-k},\Xi_M, A_M)$ в группе
кобордизмов
 $Imm^{sf; \kappa_1, \dots, \kappa_j}(n-k,k)$.

Докажем, что многообразие $T^{n-k}$ допускает ретракцию порядка
$i$ относительно когомологического класса $\kappa_T$. Докажем, что
этот когомологический класс представлен отображением $\kappa_T:
T^{n-k} \to \RP^{n-k-i-1} \subset \RP^{\infty}$ (классифицирующее
отображение обозначается также как и соответствующий
когомологический класс). Выберем натуральное $b$ достаточно
большим, и рассмотрим отображение
 $\kappa_M : M^{n-k} \to \RP^b$, для которого полный прообраз проективного подпространства
$\RP^{b-n+k+i} \subset \RP^b$ коразмерности $n-k-i$ совпадает с
подмногообразием $Q^i$.

Рассмотрим отображение $g: P^{i+1} \to \RP^{b-n+k+i}$, ограничение
которого на границу $\partial P^{i+1} = Q^i$ совпадает с
отображением $\kappa_M \vert_{Q^i}$. Рассмотрим  отображение
 $h: U_P \to \RP^b$, которое определено посредством "утолщения" отображения
$g$ с  $P^{i+1}$ на $U_P$, т.е. посредством расширения отображения
$g$ до морфизма дискового расслоения в расслоение регулярной
окрестности подмногообразия $\RP^{b-n+k+i} \subset \RP^b$.
Заметим, что отображение  $M^{n-k} \cup_{Q^i \times D^{n-k-i}} U_P
\to \RP^b$ корректно определено, при этом ограничение этого
отображения на подмногообразие
 $T^{n-k} \subset
M^{n-k} \cup U_P$ не пересекает подмногообразие
 $\RP^{b-n-k+i} \subset \RP^b$.
Поскольку открытое многообразие $\RP^{b} \setminus \RP^{b-n-k+i}$
гомотопически эквивалентно многообразию $\RP^{n-k-i-1}$, требуемая
ретракция $\kappa_T$ порядка $i$ в заданном классе когомологий
построена.

Докажем обратное утверждение, а именно, докажем, что если
скошенно-оснащенное многообразие $T^{n-k}$ допускает ретракцию
порядка $i$ относительно характеристического класса когомологий
скошенного оснащения, то оснащенное многообразие
 $(Q^i,\Psi_Q, A_Q)$, определенное как прообраз точки при отображении ретракции порядка $i-1$,
  ограничивает как скошенно-оснащенное многообразие в коразмерности $(n-i)$,
   т.е. определяет нулевой элемент $\delta(y)$ в группе кобордизма
$Imm^{sf;\kappa_1,\dots,\kappa_j}(i,n-i)$.

Пусть  $(W^{n-k+1},\Xi_W)$ -- скошенно-оснащенное многообразие с
границей
 $\partial W^{n-k+1} = M^{n-k} \cup M_1^{n-k}$, которое определяет нормальный кобордизм
между скошенно-оснащенными многообразиями
 $M^{n-k}$ и $M_1^{n-k}$. Кроме  того, предположим, что многообразие
 $M_1^{n-k}$ допускает ретракцию порядка  $i$. Отображение ретракции мы обозначим через
$\kappa'_{M_1}: M_1^{n-k} \to \RP^{n-k-i-1}$. Если натуральное
число $b$ достаточно велико (например, $b > n-k+1$), то
отображение $F: W^{n-k+1} \to \RP^b$, для которого $F
\vert_{M^{n-k}} = \kappa_M$, $F \vert_{M_1}=\kappa'_{M_1}$, т.е.
$Im(F(M^{n-k})) \subset \RP^{n-k-i} \subset \RP^b$, $Im
(F(M_1^{n-k})) \subset \RP^{n-k-i-1}$ корректно определено.
Рассмотрим подмногообразие $\RP^{b-n-k+i} \subset \RP^b$, которое
пересекает подмногообразие $\RP^{n-k-i} \subset \RP^b$
(подмногообразие $\RP^{n-k-i} \subset \RP^b$ содержит образ
многообразия $M^{n-k}$ при отображении $\kappa_M$) по точке $pt
\in \RP^{n-k-i}\setminus \RP^{n-k-i-1}$ и не пересекает
подмногообразие $\RP^{n-k-i-1} \subset \RP^{n-k-i}$
(подмногообразие $\RP^{n-k-i-1} \subset \RP^{n-k-i}$ содержит
образ многообразия $M_1^{n-k}$ при отображении $\kappa'_{M_1}$).

Обозначим через $P^{i+1}$ подмногообразие $F^{-1}(\RP^{b-n+k+i})$.
По построению $\partial P^{i+1} = Q^i$. Определим скошенное
оснащение $\Xi_P$ размерности  $(n-i)$ как прямую сумму оснащений
подмногообразия $P^{i+1} \subset W^{n-k+1}$ и оснащения $\Xi_W$,
которое получается в результате ограничения оснащения c
$W^{n-k+1}$ на подмногообразие  $P^{i+1} \subset W^{n-k+1}$.

Заметим, что ограничение оснащения $\Xi_W$ на $\partial W^{n-k+1}
=Q^{n-k}$ совпадает с оснащением $\Xi_Q$. Аналогичное утверждение
справедливо для семейства коциклов
 $A_P=\{\kappa_P, \kappa_1, \dots, \kappa_j\}$, полученных в
 результате ограничения семейства коциклов с $W^{n-k+1}$ на $P^{i+1}$.
 Следовательно, элемент $y$, представленный
тройкой $(Q^i, \Xi_Q, A_Q)$, лежит в ядре гомоморфизма
 $\delta$. Лемма  35 доказана.

\subsubsection*{Доказательство Леммы 33}

По Лемме  35, препятствие к ретракции порядка $i$ в классе
кобордизма тройки $(M^{n-k}, \Xi_M, A_M)$, представляющей заданный
элемент $x$ в группе кобордизма $Imm^{sf;\kappa_1, \dots,
\kappa_j}(n-k,k)$, совпадает (в предположении, что ретракция
порядка $i-1$ уже построена) с препятствием к ретракции того же
порядка $i$ для многообразия $M_1^{n-k'}$, которое является
образом элемента $x$ в группе $Imm^{sf;\kappa_1, \dots,
\kappa_j}(n-k',k')$ при гомоморфизме $J_{k'}^{sf}:
Imm^{sf,\kappa_1, \dots, \kappa_j}(n-k,k) \to Imm^{sf, \kappa_1,
\dots, \kappa_j}(n-k',k')$. При построении последовательности
ретракций все более высокого порядка будем сравнивать, в
соответствии с Леммой 31, необходимые и достаточные условия
ретракции следующего порядка для элемента из исходной группы
$Imm^{sf;\kappa_1, \dots, \kappa_j}(n-k,k)$ и для образа этого
элемента в группе $Imm^{sf;\kappa_1, \dots, \kappa_j}(n-k',k')$.
На каждом шаге конструкции оба препятствия совпадают. Лемма 33
доказана.
\[  \]

Для доказательства теоремы о ретракции нам потребуется
конструкция, принадлежащая У.Кошорке (см.[K]), полного препятствия
к построению послойного мономорфизма двух векторных расслоений.

Предположим, что  $\alpha \to Q^q$, $\beta \to Q^q$ -- пара
векторных расслоений над многообразием $Q^q$ (вообще говоря, $Q^q$
является многообразием с границей) $dim(\alpha)=a$,
$dim(\beta)=b$, $dim(Q^q)=q$, при этом предполагается, что
выполнено $2(b-a+1)<q$. Пусть $u: \alpha \to \beta$ -- морфизм
общего положения, обозначим через $\Sigma \subset Q^q$ --
подмногообразие, определенное формулой
$$\Sigma = \{ x \in Q^q \vert ker(u_x: \alpha_x \to \beta_x) \ne 0 \},$$
т.е. подмногообразие точек базы с вырождением слоев при морфизме
$u$. Заметим, что при рассматриваемых размерностных ограничениях,
для морфизма $u_x$ общего положения мы имеем $rk(u_x) \ge a-1$.
Коразмерность подмногообразия $\Sigma \subset Q^q$ при этом равна
$b-a+1$.

Опишем нормальное расслоение к подмногообразию $\Sigma \subset
Q^q$. Обозначим через $\lambda: E(\lambda) \to \Sigma$ -- линейное
расслоение, определенное как поле ядер морфизма $u$ над особым
подмногообразием $\Sigma \subset Q^q$. Определено, тем самым,
включение расслоений $\varepsilon: \lambda \subset \alpha$.
Обозначим через $\Lambda_{\alpha}$ -- расслоение над $\Sigma$,
определенное как ортогональное дополнение подрасслоения
$\varepsilon(\lambda) \subset \alpha$. Естественный морфизм
(послойный мономорфизм) $v_x: \Lambda_{\alpha} \subset \beta$
определен. Расслоение $\Lambda_{\beta}$ определено как расслоение,
слои которого ортогональны над точками $x \in \Sigma$ в слоях
$\beta_x$ расслоения $\beta$ подпространствам
$v_x(\Lambda_{\alpha_x})$. Нормальное расслоение $\nu(\Sigma)$  к
подмногообразию $\Sigma \subset Q$ определено по формуле

\begin{eqnarray}\label{16}
\nu(\Sigma) = \lambda \otimes \Lambda_{\beta}.
\end{eqnarray}

В работе [K2] (в которой приводится ссылка на более ранние работы
того же автора) определяется группа кобордизма вложенных
многообразий коразмерности $b-a+1$, причем нормальное расслоение
подмногообразия, представляющего элемент этой группы кобордизма,
определено по формуле ($\ref{16}$). Указанная структура
нормального расслоения сохраняется при кобордизме. Для
произвольного морфизма $u: \alpha \to \beta$ общего положения
определяется представитель, задающий элемент в описанной группе
кобордизма, который служит полным препятствием для построения
послойного мономорфизма указанных расслоений.

Если $Q^q$ -- многообразие с границей, причем морфизм расслоений,
заданных на всем $Q^q$, определен над границей $\partial Q$, то
многообразие $S_{\partial}$ критических слоев над границей служит
границей подмногообразия $\Sigma \subset Q^q$, $\partial \Sigma=
\Sigma_{\partial}$ c нормальным расслоением, определенным по
аналогичной формуле.

Пусть $\nu \to \RP^{2^k-1}$ -- векторное расслоение размерности
$(n+1-2^k)$, $2^k < n+2$ над вещественным проективным
пространством, изоморфное сумме Уитни $(n+1-2^k)$ экземпляров
линейного нетривиального расслоения $\kappa_{\RP}$ над
$\RP^{2^k-1}$. Обозначим расслоение $\nu \oplus \kappa$ через
$\bar \nu$. Определена естественная проекция $\pi: \bar \nu \to
\nu$ c ядром $\kappa$. Заметим, что расслоение $\bar \nu \to
\RP^{2^k-1}$ является нормальным расслоением к многообразию
$\RP^{2^k-1}$ размерности $(n-2^k+2)$. По теореме Коэна [C]
стандартное проективное пространство $\RP^{2^k-1}$ погружается в
евклидово пространство размерности $2^{k+1}-2-k$.

\subsubsection*{Замечание}
Чтобы ослабить размерностные ограничения в рамках данного
доказательства, можно воспользоваться более сильным результатом о
том, что проективное пространство $\RP^{2b(k)-1}$ погружается в
$\R^{2b(k)-1-4k}$, где $b(k)$-- число Радона-Гурвица, равное
соответствующей степени числа $2$, причем $b(k)<2^k$. Это
утверждение эквивалентно существованию семейства $(n+3-2^{k+1}+k)$
независимых сечений расслоения $\bar \nu$.
\[  \]

Определим понятие стандартного семейства $(n+3-2^{k} + k)$ сечений
(особенности семейства сечений допускаются) расслоения $\nu$.

\subsubsection*{Определение стандартного семейства сечений расслоения $\nu$}
Выберем невырожденное семейство сечений $\bar \Psi = \{\bar
\psi_1, \dots, \bar \psi_s\}$, $s=(n+3-2^{k+1}+k)$ расслоения
$\bar \nu$. Рассмотрим семейство сечений  общего положения $\Psi =
\{\psi_1, \dots, \psi_s\}$, которое получено в результате проекции
семейства $\bar \Psi$ при $\pi: \bar \nu \to \nu$. Семейство
$\Psi$ назовем стандартным.
\[  \]

Обозначим через $\Sigma \subset \RP^{2^k-1}$ подмногообразие точек
базы, для которых стандартное семейство сечений имеет особенности.
Обозначим через $\nu_{\Sigma}$ нормальное расслоение к этому
подмногообразию.

\begin{lemma}
 Многообразие $\Sigma \subset \RP^{2^k-1}$ является $k$-мерным подмногообразием, причем его нормальное
 расслоение $\nu_{\Sigma}$ снабжено скошенным оснащением $\Xi_{\Sigma}$, характеристический класс
 $\kappa_{\Sigma}$ которого совпадает с ограничением характеристического класса $\kappa \in H^1(\RP^{2^k-1};\Z/2)$
 на подмногообразие $\Sigma$, $\kappa_{\Sigma} = \kappa \vert_{\Sigma \subset \RP^{2^k-1}}$.
 \end{lemma}

 \subsubsection*{Доказательство Леммы 36}

Вычислим нормальное расслоение $\nu_{\Sigma}$ по теореме Коршорке.
Более того, определим скошенное оснащение этого расслоения. Пусть
$\lambda \subset \bar \Psi$ -- расслоение ядер семейства сечений
над особым подмногообразием $\Sigma$. Поскольку семейство сечений
является стандартным, $\lambda = \kappa_{\Sigma}$. Расслоение
$\bar \Psi/\lambda$ является обратным к $\lambda$. С другой
стороны, расслоение $\nu \vert_{\Sigma}$ само обращает расслоение
$\kappa_{\Sigma}$. Следовательно, ортогональное дополнение к $\bar
\Psi/\lambda \vert_{\Sigma}$ в расслоении $\nu$ является
тривиальным расслоением размерности $(2^k-k-1)$. По теореме
Кошорке расслоение $\nu_{\Sigma}$ естественно изоморфно расслоению
$(2^{k}-k-1) \varepsilon \otimes \kappa = (2^{k}-k-1) \kappa$.
Этот изоморфизм определяет скошенное оснащение расслоения
$\nu_{\Sigma}$ с характеристическим классом $\kappa_{\Sigma}$.
Лемма 36 доказана.
\[  \]

Пусть теперь многообразие $M^m$ имеет размерность
$dim(M)=m=2^{k}-2$ и является скошенно-оснащенным в коразмерности
$n-dim(M)$, где $n=-2(mod 2^{k})$. Рассмотрим отображение
 $\kappa_M: M^m \to \RP^{2^{k}-1}$, заданное характеристическим
  классом
оснащения многообразия $M^m$. Предположим, что $M^m$ снабжено
набором когомологических классов $A_M= \{\kappa_M, \kappa_1,
\dots, \kappa_j\}$,  при этом определен элемент  $\alpha \in
Imm^{sf,\kappa_1, \dots, \kappa_j}(2^{k}-2,n-2^{k}+2)$.

Рассмотрим расслоение $\bar \nu \to \RP^{2^k-1}$, $\bar \nu =
(n-2^k+2)\kappa_{\RP}$ размерности $\dim(\bar \nu)=n-2^k+2$.
Заметим, что нормальное расслоение $\nu_M$ той же размерности
изоморфно обратному образу расслоения $\bar \nu$ при отображении
$\kappa_M: M^m \to \RP^{2^k-1}$, $\nu_M = \kappa_M^{\ast}(\bar
\nu)$. Рассмотрим отображение $c_M: M^{m} \to \RP^{2^k-1} \times
\prod_{i=1}^j \RP^{\infty}$, определенное по набору классов
когомологий $A_M$. Рассмотрим отображение проекции $\pi:
\RP^{2^k-1} \times \prod_{i=1}^j \RP^{\infty} \to \RP^{2^k-1}$ на
первый сомножитель. Композиция $\pi \circ l: M^m \to \RP^{2^k-1}$,
очевидно, совпадает с отображением $\kappa_M$. Определим
подрасслоение $\phi_M \subset \nu_M$ коразмерности 1 (т.е.
размерности $n+1-2^k$) такое, что $\phi_M = \kappa_M^{\ast}(\nu)$.

Индуцируем набор сечений $\Xi_M=\{\xi_1, \dots, \xi_s\}$,
$s=n+3-2^{k+1}+k$ расслоения $\phi_M$ из стандартного набора
сечений $\Psi = \{\psi_1, \dots \psi_s\}$ расслоения $\nu$
отображением $\kappa_M$. Обозначим через $N^{k-1} \subset M^m$
подмногообразие размерности $(k-1)$ особых семейств сечений.
Многообразие $N^{k-1}$ снабжено набором когомологических классов
$A_N$, полученный ограничением соответствующих классов набора
$A_M$ на это подмногообразие. При этом класс $\kappa_M \vert_N$
обозначается через $\kappa_N$.

Нормальное расслоение $\nu_N$ изоморфно сумме Уитни $\nu_N = \nu_M
\vert_N \oplus \nu_{N \subset M}$, где $\nu_M$--нормальное
расслоение многообразия $M^m$,  $\nu_{N \subset M}$ -- нормальное
расслоение подмногообразия $N^{k-1}$ внутри $M^m$ снабжено
скошенным оснащением с характеристическим классом $\kappa_N$.
Заметим, что расслоение $\nu_M \vert_N$ также снабжено скошенным
оснащением с тем же характеристическим классом, что доставляет
скошенное оснащение $\Xi_N$ многообразия $N^{k-1}$ в коразмерности
$(n-k+1)$.

\begin{lemma}
Тройка $(N^{k-1}, \Xi_{N}, A_{N})$ определяет элемент $ x_{k-1}
\in Imm^{sf,\kappa_1, \dots, \kappa_j}(k-1,n-k+1)$, который
является полным препятствием к построению ретракции порядка
$(k-1)$ для исходного элемента $x \in Imm^{sf,\kappa_1, \dots,
\kappa_j}(2^{k}-2,n-2^{k}+2)$.
\end{lemma}

\subsubsection*{Доказательство Леммы 37}

Рассмотрим многообразие $\Sigma^k \subset \RP^{2^{k}-1}$ особых
сечений семейства $\Psi$, которое ниже будет переобозначено через
$\Sigma_0^{k}$. Это многообразие снабжено естественной
стратификацией (фильтрацией):
\begin{eqnarray}\label{17}
\emptyset \subset \Sigma_{k}^0 \subset \dots \subset
\Sigma_{1}^{k-1} \subset \Sigma^{k}_0 \subset \RP^{2^{k}-1}.
\end{eqnarray}

Не ограничивая общности, можно считать, что каждое многообразие в
фильтрации ($\ref{17}$) является связным. Подмногообразие
$\Sigma_i$, $dim(\Sigma_i)=k-i$ этой фильтрации определен как
особое подмногообразие подсемейства сечений, состоящего из первых
$(s-i)$-сечений. Из прямых вычислений вытекает, что
фундаментальный класс
 $[\Sigma_i]$ соответствующего члена фильтрации представлен образующей в
$H_{k-i}(\RP^{2^{k}-1};\Z/2)$. Действительно, это следует из того,
что указанный цикл двойственен характеристическому классу
$w_{n+1-2^{k}-k+i}((n+1-2^{k})\kappa_{\RP^{2^{k}-1}})$.

Предположим, что отображение $\kappa_M : M^m \to \RP^{2^{k}-1}$
является трансверсальным вдоль построенной стратификации в образе
и прообраз стратификации  $(\ref{17})$ мы обозначим через
\begin{eqnarray}\label{18}
 N_{k-1}^0 \subset N_{k-2}^1 \subset \dots \subset N_{0}^{k-1}
\subset M^m,
\end{eqnarray}
причем старшее
многообразие фильтрации $N_0^{k-1}$ совпадает с многообразием
$N^{k-1}$, которое было построено ранее.

Проведем рассуждения  индукцией по параметру $i$, $i=0, \dots,
k-1$. Рассуждения в случае $i=0$ элементарны. Четномерное
многообразие $M^m$ является ориентированным, поэтому представляет
нулевой класс гомологий в $\RP^{2^k-1}$.

Предположим, что образ отображения $\kappa_M : M^m \to
\RP^{2^{k}-1}$ уже лежит в стандартном проективном подпространстве
$\RP^{2^{k}-2-i} \subset \RP^{2^{k}-1}$. Тогда
$N^{i-1}_{k-i}=\emptyset$. Заметим, что страт $\Sigma_{k-i}^i$
можно выбрать в классе изотопии так, что множество точек
пересечения $\Sigma_{k-i}^i \cap \RP^{2^{k}-1-i}$ состоит из
единственной точки. Действительно, индекс пересечения
фундаментальных классов подмногообразий $\RP^{2^{k}-1-i},
\Sigma_{k-i}^{i}$ в многообразии $\RP^{2^{k}-1}$ нечетен и может
принимать произвольные нечетные значения. Применяя прием Уитни,
можно выбрать $\Sigma_{k-i}^{i}$ в классе изотопии так, что
геометрически указанные подмногообразия пересекаются по минимально
возможному множеству точек, т.е. точка пересечения всего одна.

Оснащенное подмногообразие $N_{k-i-1}^i$ является регулярным
прообразом точки при отображении в остов $\RP^{2^{k}-2-i}$ и
представляет элемент $x_i$ из группы кобордизма $Imm^{sf;\kappa_1,
\dots, \kappa_j}(i,n-i)$. Тогда по Лемме 32, поскольку $x_i=0$,
(что сразу следует из предположения $x_{k-1}=0$), существует
нормальный кобордизм отображения $c_M$ такой, что в результате
получается отображение (которое мы обозначим через $c'_M$), для
которого образ $\pi \circ c'_M$ лежит в подпространстве
$\RP^{2^{k}-2-i}$. Далее переходим от $i$ к $(i+1)$ и повторяем
аналогичные рассуждения. В предположении равенства $x_{k-1}=0$ в
классе кобордизма можно построить ретракцию порядка $(k-1)$. Лемма
37 доказана.
\[  \]

Основным при доказательстве теоремы ретракции является следующее
предложение.

\begin{proposition}
Пусть $n=-2(mod(2^k))$, $n>2^k$ и пусть для заданного элемента $x
\in Imm^{sf;\kappa_1, \dots, \kappa_j}(2^k-2,n-2^k+2)$
рассматривается препятствие $x_{k-1} \in Imm^{sf; \kappa_1, \dots,
\kappa_j}(k-1,n-k+1)$ к построению ретракции порядка $(k-2)$ в
классе кобордизма $x$. Тогда элемент $x_{k-1}$ лежит в образе
гомоморфизма трансфера, т.е. существует элемент $y_{k-1} \in
Imm^{sf,\kappa_1, \dots, \kappa_{j+1}}(k-1,n-k+1)$, для которого
$r_{j+1}(y_{k-1})=x_{k-1}$.
\end{proposition}

\subsubsection*{Доказательство Предложения 38}

Рассмотрим нормальное расслоение $\nu_M=(n-2^{k}+2)\kappa_M$ к
многообразию $M^{m}$, $dim(M^m)=m=2^{k}-2$ и подрасслоение $\phi_M
\subset \nu_M$, $\phi =  (n-2^{k}+1)\kappa_M$ коразмерности 1,
которое ортогонально линейному подрасслоению $\kappa_M \subset
\nu_M$. Существует невырожденное семейство сечений
 $\{\xi_1, \dots, \xi_s \}$,
$s=n+3+k - 2^{k+1}$, расслоения $\phi$.

Действительно, расслоение $\phi_M$ является нормальным расслоением
к многообразию (с краем), которое определено как  тотальное
пространство дискового расслоения $D(\kappa_M)$. Поскольку
$dim(D(\kappa_M))=2^k-1$, $\alpha(2^k-1)=k$, по теореме Коэна [C]
многообразие $D(\kappa_M)$ погружается в евклидово пространство
$\R^{2^{k+1}-2-k}$ (это утверждение сформулировано как Гипотеза 39
и прямое доказательство частного случая этой гипотезы приводится в
Следствии 40). В частности, расслоение $\phi_M$ допускает
семейство $s$ невырожденных сечений. Не ограничивая общности,
можно считать, что скошенно-оснащенное многообразие $M^{2^k-2}$
выбрано в классе нормального кобордизма.

Рассмотрим, кроме того, семейство сечений $\Psi=\{\psi_1, \dots,
\psi_s\}$, которое индуцировано из стандартного семейства сечений
расслоения $(n+1-2^k) \kappa_{\RP^{2^{k}-1}}$ при отображении
$\kappa_M$. Определено поднятие  семейства $\Psi$ в невырожденное
семейство сечений $\bar \Psi=\{\bar \psi_1, \dots, \bar \psi_s\}$
расслоения $\nu_M$. Семейство $\bar \Psi$ проектируется семейство
$\Psi$ при проекции $\nu_M \to \phi_M$ вдоль подрасслоения
$\kappa_M$. Поднятие невырожденного семейства   $\{\xi_1, \dots,
\xi_s \}$ в невырожденное семейство сечений расслоения $\nu_M$
обозначим через $\{\bar \xi_1, \dots, \bar \xi_s \}$.

Рассмотрим, кроме того, семейство сечений общего положения $\bar
X= \{\bar \chi_1, \dots \bar \chi_s\}$ на многообразии
$M^{2^{k}-1} \times I$ расслоения $\nu(M^{2^{k}-1}) \times I$,
которое совпадает с семейством $\{\bar \psi_1, \dots, \bar
\psi_s\}$ на компоненте границы $M^{2^{k}-1} \times \{1\}$ и с
семейством сечений $\{\bar \xi_1, \dots, \bar \xi_s\}$ на
компоненте $M^{2^{k}-1} \times \{0\}$ границы. Семейство сечений
$\{\bar \chi_1, \dots, \bar \chi_s \}$ имеет, вообще говоря,
особенности. Подмногообразие особых сечений обозначим через
$V^{k-1} \subset M^{m} \times I$. Это замкнутое многообразие
размерности $k-1$. Аналогично введем обозначение $\{\chi_1, \dots,
\chi_s \}$ для проекции рассматриваемого семейства сечений в
семейство сечений расслоения $\phi_M \times I$ (обозначим это
расслоение через $\phi_{M \times I}$). Семейство сечений
$\{\chi_1, \dots, \chi_s \}$
 также
рассматривается с заданными граничными условиями, определенными
проекциями семейств $\{\bar \psi_1, \dots, \bar \psi_s\}$, $\{\bar
\xi_1, \dots, \bar \xi_s\}$.

Пусть $K^{k} \subset M^{m} \times I$, $m=2^k-2$ -- подмногообразие
размерности $k$, определенное как подмногообразие особенностей
семейства сечений $\{\chi_1 \dots, \chi_s\}$. Это подмногообразие
имеет компоненту границы $K^{k} \cap (M^{m} \times \{1\})$,
которую мы обозначим через $N^{k-1} \subset M^m \times \{1\}$.
Согласно Лемме 33, тройка $( N^{k-1}, \Xi_N, A_N)$, где $\Xi_N$ --
скошенное оснащение нормального расслоения подмногообразия
$N^{k-1}$, определенное по теореме Кошорке, и семейство
когомологических классов $A_N = A_M \vert_{N^{k-1}}$, представляет
элемент $x_{k-1}$ в группе $Imm^{\kappa_1, \dots, \kappa_j;
sf}(k-1,n-k+1)$, который определяет полное препятствие к
построению ретракции порядка $(k-1)$ в классе кобордизма $x
=[(M^m, \Xi_M, \kappa_M)]$. Рассмотрим поднятие семейства сечений
$\{\chi_1 \dots, \chi_s\}$ расслоения $\phi_{M \times I}$ до
семейства сечений $\{\bar \chi_1 \dots \bar \chi_s\}$ общего
положения расслоения $(n+2-2^{k})\kappa_{M \times I}$. Граничные
условия записываются в виде:

\begin{eqnarray}\label{19}
\bar \chi_l \vert_{M^m \times \{0\}} = \bar \xi_l, \qquad \bar
\chi_l \vert_{M^m \times \{1\}} = \bar \psi_l, \qquad l=1, \dots,
s.
\end{eqnarray}
 Заметим,
что определено включение $V^{k-1} \subset K^{k}$ подмногообразий
коразмерности 1, поскольку над $V^{k-1}$ семейство сечений $\bar
X= \{ \bar \chi_1 \dots \bar \chi_s \}$ вырождается.
Подмногообразие $V^{k-1}$ замкнуто, поскольку на границе
$\partial(M^m \times I)$ семейство сечений $\bar X$ не имеет
особенностей. Применим формулу $(\ref{16})$ для вычисления
стабильного нормального расслоения к подмногообразию $K^{k}$ и к
подмногообразию $V^{k-1} \subset K^{k}$. Обозначим через $\bar
\lambda$ -- линейное расслоения ядер семейства сечений $\bar X$
над $V^{k-1}$.

Докажем, что нормальное расслоение подмногообразия $V^{k-1}
\subset M^{m} \times I$ в объемлющем многообразии изоморфно
расслоению $\varepsilon \oplus  (m-k+1)\lambda$. Действительно,
ограничение нормального расслоения $\nu(M^m)\vert_{V}$ изоморфно
расслоению $(n+2-2^{k})\kappa_{M \times I}$, и, кроме того,
изоморфно тривиальному расслоению $(n+2-2^{k})\varepsilon$, причем
изоморфизм выбирается каноническим (не зависящим ни от $M^m$ ни от
$V^{k-1}$). Расслоение $\bar \Lambda$, дополнительное к линейному
расслоению $\bar \lambda$ в тривиальном расслоении линейных
комбинаций векторов сечений, представляет расслоение $-\bar
\lambda$, обратное к расслоению $\lambda$. Следовательно,
ортогональное дополнение к образу этого расслоения в расслоении
$\nu_M \vert_V$ изоморфно расслоению $\bar \lambda \oplus (2^{k} -
k -1)\varepsilon$. Искомое нормальное расслоение подмногообразия
$V^{k-1} \subset M^{2^{k}-2} \times I$ (тем самым, и стабильное
нормальное расслоение многообразия $V^{k-1}$) изоморфно расслоению
$\bar \lambda \otimes(\bar \lambda \oplus (2^{k}-k-1)\varepsilon)
= \varepsilon \oplus (2^{k}-k-1)\bar \lambda$.

Обозначим через $U_V \subset K^{k}$--регулярную окрестность
подмногообразия $V^{k-1} \subset K^{k}$. Поле ядер семейства
$\{\chi_1, \dots, \chi_s\}$, ограниченное на $V^{k-1}$, (а,
следовательно, и на $U_V$) задано линейным расслоением $\lambda$,
ограниченным с $K^{k}$ на $V^{k-1}$, (на $U_V$) и это расслоение
обозначается через $\bar \lambda$. С учетом того, что ограничение
каждого из расслоений $(2^{k-1})\bar \lambda$, $(2^{k-1})\kappa_M$
на подмногообразие $V^{k-1}$ канонически изоморфно тривиальному
расслоению, ортогональное дополнение к расслоению сечений в
расслоении $\phi \times I \vert_{V^{k-1}}$ изоморфно расслоению
$\bar \lambda \oplus (2^{k-1}-k)\varepsilon \oplus
(2^{k-1}-1)\kappa_M$. При этом по аналогичным вычислениям,
нормальное расслоение к подмногообразию $K^{k} \subset M^{2^{k}-2}
\times I$ (внутри объемлющего многообразия) изоморфно расслоению
$\varepsilon \oplus (2^{k-1}-1-k)\bar \lambda \oplus
(2^{k-1}-1)\kappa_M \otimes \lambda$. Откуда следует, что
стабильный класс изоморфизма рассматриваемого расслоения
вычисляется по формуле $(-1-k)\bar \lambda -(\kappa_M \otimes \bar
\lambda)$.

 Отсюда, в частности, вытекает, что линейное
 нормальное расслоение подмногообразия $V^{k-1} \subset K^{k}$, совпадающее с линейным расслоением
 подмногообразия $V^{k-1}$ в регулярной окрестности $U_V$, изоморфно
расслоению $\bar \lambda \otimes \kappa_M$. Обозначим через
$Q^{k-1}$ -- границу регулярной окрестности $U_V$, т.е. $Q^{k-1} =
\partial U_V$.
 Многообразие $K^{k} \setminus U_V$ имеет компоненты
границы $\partial(K^{k} \setminus U_V)= Q^{k-1} \cup K^{k-1}$ и
является скошенно-оснащенным в коразмерности $n-k$ с
характеристическим классом оснащения $\kappa_{M\times I} \vert_{K
\setminus U_V}$. Действительно, нормальное расслоение этого
подмногообразия было описано выше, а полное нормальное расслоение
получается из нормального расслоения подмногообразия в результате
суммы Уитни с тривиальным $(n+2-2^k)$--мерным расслоением, которое
канонически изоморфно расслоению $(n+2-2^k)\kappa_{M \times I}$.

Тройка $(N^{k-1}, \Xi_N, A_N)$ ($Q^{k-1}, \Xi_Q, A_Q)$), состоящая
из многообразия $N^{k-1}$ ($Q^{k-1}$) вместе со скошенным
оснащением $\Xi_K$ ($\Xi_Q$) и семейством  $A_K$ ($A_Q$)
когомологических классов, которое  определено как ограничение
семейства $A_{M\times I}=\{\kappa_1, \dots, \kappa_j\}$ на
рассматриваемое подмногообразие c объемлющего многообразия
$M^{2^{k}-2} \times I$, определяют соответствующие элементы в
группе кобордизма $Imm^{sf;\kappa_1, \dots, \kappa_j}(k-1,n-k+1)$.
Для краткости эти элементы будем обозначать через $[K]$, $[Q]$, не
указывая дополнительную структуру скошенного оснащения и набора
классов когомологий.

По аналогичным соображениям многообразие $K^{k} \setminus U_V$,
снабженно набором когомологических классов
$\{\kappa_1\vert_{N\setminus U_V}, \dots, \kappa_j \vert_{N
\setminus U_V}\}$ и скошенным оснащением в коразмерности $(n-k+1)$
c характеристическим классом $\kappa_{M^m\times I}
\vert_{K^{k}\setminus U_V}$, определяет скошенно-оснащенный
кобордизм между $[K]$ и $[Q]$. По построению $[K]=\alpha_{k-1}$
представляет полное препятствие к ретракции порядка $(k-1)$
исходного элемента $(M^{2^{k}-2}, \Xi_M, A_M)$ относительно
характеристического класса $\kappa_M$. Элемент $[Q]$ получен в
результате гомоморфизма трансфера из некоторого элемента
$\beta_{k-1} \in Imm^{sf;\kappa_1, \dots, \kappa_j,
\kappa_{j+1}}(k-1,n-k+1)$. Элемент $y_{k-1}$ задан
скошенно-оснащенным многообразием $V^{k-1}$ и набором
когомологических классов $A_V$, который индуцирован ограничением
набора соответствующих классов с $M^{2^{k}-2} \times I$. Заметим,
что многообразие $Q^{k-1}$ служит двулистным накрывающим над
$V^{k-1}$. При этом когомологический класс $\kappa_{j+1} \in
H^1(V^{k-1};\Z/2)$, относительно которого берется трансфер,
совпадает с $\kappa_{M \times I} \vert_V$. Характеристический
класс $\kappa_Q$ скошенного оснащения $\Xi_Q$ многообразия
$Q^{k-1}$  получен при трансфере характеристического класса $\bar
\lambda = \kappa_V$ скошенного оснащения $\Xi_V$. (Заметим, что
этот же характеристический класс $\kappa_Q$ индуцируется из класса
когомологий $\kappa_{j+1} = \kappa_M \vert_{V}$ при накрытии
$Q^{k-1} \to V^{k-1}$.) Предложение 38 доказано.

\subsubsection*{Доказательство Теоремы 29 о ретракции}

Рассмотрим натуральное число $\psi$, определенное в Предложении 32
для значения параметра ретракции равного $d-1$. Тогда гомоморфизм
$(\ref{104})$ полного трансфера на группе $Imm^{sf;\kappa_1,
\dots, \kappa_{\psi}}(d-1,n-d+1)$ оказывается тривиальным.
 Определим натуральное число
$l(d) = \exp_2(\exp_2 \dots \exp_2(d)\dots +1)$, где число
итераций функции $\exp_2(x+1) = 2^{x+1}$ равно $\psi$, а начальное
значение $x=d-1$.

Пусть $l'$-- произвольная степень двойки, не меньшая, чем $l(d)$.
Определим $n=l'-2$. Докажем, что произвольный элемент группы
$Imm^{sf}(n-1,1)$ допускает ретракцию порядка $d$.

Положим $n_0=l(d)-2$, по предположению $n_0 \le n$.  Определим
последовательность из $\psi$ чисел: $2n_1 = log(n_0+2)-2$, $2n_2=
log_2(n_1+2)-2$, $\dots,$
 $2n_{\psi}=log_2(n_{\psi-1} +2)-2$. Все
 указанные числа являются натуральными, при этом $n_{\psi}=d-1$.

Пусть $x \in Imm^{sf}(n-k,k)$-- произвольный элемент, $x_{d-1} \in
Imm^{sf}(d-1,n-d+1)$--полное препятствие к построению ретракции
порядка $d-1$ для исходного элемента $x$. Рассмотрим образ этого
элемента $x_{n_0} \in Imm^{sf}(n_0,n-n_0)$. Заметим, что
препятствие $x_{d-1}$ является также препятствием к ретракции
порядка $(d-1)$ элемента $x_{n_0}$.

Рассмотрим препятствие  $x_{n_1} \in Imm^{sf}(n_1,n-n_1)$ к
построению ретракции порядка $n_1$ для $\alpha_{n_0}$. По
Предложению 34, существует элемент $y_{n_1} \in
Imm^{sf;\kappa_1}(n_1,n-n_1)$, представленный скошенно-оснащенным
многообразием $X^{n_1}_1$ с характеристическим классом
$\kappa_{X^{n_1}_1} \in H^1(X^{n_1}_1;\Z/2)$ скошенного оснащения,
который переходит при трансфере относительно когомологического
класса $\kappa_1=\kappa_N \vert_{X_1}$ в элемент $x_{n_1}$.

Рассмотрим препятствие $y_{n_2} \in Imm^{sf;\kappa_1}(n_2,n-n_2)$
к ретракции порядка $n_2$ для элемента $y_{n_1}$. Снова по
Предложению 34, существует элемент $z_{n_2} \in Imm^{sf; \kappa_1,
\kappa_2}(n_2,n-n_2)$ такой, что элемент $r_2(z_{n_2})$
препятствует к ретракции порядка $n_2$ для элемента $y_{n_1}$.
Элемент $r_{tot}(z_{n_2})=r_1 \circ
r_2(z_{n_2})=J^{sf}_{n_2}(x_{n_0})=x_{n_2}$ служит препятствием к
ретракции порядка $(n_2)$ для элемента $x_{n_1}$. Этот же элемент
$x_{n_2}$ служит полным препятствием к ретракции порядка $(n_2)$
для исходного элемента $x$.

Рассуждая по индукции, получаем, что препятствие к построению
ретракции порядка $(d-1)$ для элемента $x$ получено в результате
применения полного $\psi$-кратного трансфера к некоторому элементу
$\varepsilon \in Imm^{sf;\kappa_1, \dots,
\kappa_{\psi}}(d-1,n-d+1)$. По Предложению 32
$r_{tot}(\varepsilon)=0$. Следовательно, $x_{d-1}=0$. Теорема 29 о
ретракции доказана.
\[  \]

\subsection*{Добавление}

Обозначим через $b=b(k)$-- произвольную степень двойки, кратную
степени $2^{k}$. При этом предпологается, что $n+2 \ge 2^{b(k)}$.
Обозначим натуральное число $(n-2b(k)+k+3)$ через $s$.

\begin{conjecture}
Произвольный класс кобордизма $\alpha \in Imm^{sf,\kappa_1, \dots,
\kappa_j}(b(k)-1,n-b(k)+1)$, представлен многообразием
$M^{b(k)-1}$ со скошенным оснащением $\Xi_M$ и набором
характеристических классов $\kappa_{M} \in H^1(M^{b(k)-1};\Z/2)$,
$\kappa_i \in H^1(M^{b(k)-1};\Z/2)$, $i=1, \dots, j$ таким, что
многообразие $M^{b(k)-1}$ погружается в пространство
$\R^{2b(k)-2-k}$. Эквивалентно, найдется семейство  $\{\xi_1,
\dots, \xi_{s}\}$ сечений нормального расслоения
$\nu(M)=(n-b(k)+1)\kappa_M$, которое не имеет вырождений.
\end{conjecture}

При доказательстве Предложения 38 мы воспользовались следующим
следствием Гипотезы 39.

\begin{corollary}
Произвольный класс кобордизма $\alpha \in Imm^{sf; \kappa_1,
\dots, \kappa_j}(b(k)-1,n-b(k)+1)$,  представлен многообразием
$M^{b(k)-1}$ со скошенным оснащением $\Xi_M$ и набором
характеристических классов $\kappa_{M} \in H^1(M^{b(k)-1};\Z/2)$,
$\kappa_i \in H^1(M^{b(k)-1};\Z/2)$, $i=1, \dots, j$. При этом
найдется семейство $\{\xi_1, \dots, \xi_{s}\}$ сечений нормального
расслоения $\nu(M)=(n-b(k)+1)\kappa_M$, которое не имеет
вырождений над подмногообразием  $N^{b(k)-2} \subset M^{b(k)-1}$,
представляющим эйлеров класс расслоения $\kappa_M$.
\end{corollary}

Мы докажем Следствие 40 в предположении $j=0$. При этом мы
воспользуемся основным результатом работы [A-E], согласно которому
естественный гомоморфизм $Imm^{fr}(k,1) \to Imm^{sf}(k,b(k)-1)$
является тривиальным (это утверждение доказано в [A-E] даже при
более слабых размерностных ограничениях).

\subsubsection*{Замечание}
Гипотеза 39 (для произвольного $j$ вытекает из результата Коэна
[C], согласно которому многообразие $M^{b(k)-1}$ допускает
погружение в пространство  $\R^{2b(k)-k-2}$. Из существования
такого погружения вытекает, что нормальное расслоение $\nu(M)$
размерности $dim(\nu)=n+1-b(k)$ имеет, по меньшей мере,
$(n-2b(k)+k+3)$ невырожднных сечения. Над подмногообразием
$N^{b(k)-2} \subset M^{b(k)-1}$ это семейство сечений подавно не
вырождается.
\[  \]

\subsubsection*{Доказательство Следствия 40 при $j=0$}

Доказательство основано на теореме о
погружаемости стандартного проективного пространства
$\RP^{b(k)-1}$ в евклидово пространство $\R^{2b(k)-k-2}$.

Обозначим для краткости  $dim(M)=b(k)-1$ через $m$. Тогда
$dim(N)=m-1$. Обозначим расслоение $\nu_M \vert_{N^{m-1}}$ через
$\nu_N$, нормальное расслоение к многообразию $N^{m-1}$ через
$\nu_N$, расслоение $\kappa_M \vert_N$ через  $\kappa_N$. Заметим,
что $(n-b(k)+2)\kappa_N = \nu_N$ и что $\nu_N \oplus \kappa_N =
\bar \nu_N$. Определена проекция $\pi: \bar \nu_N \to \nu_N$
расслоений с ядром $\kappa_N$.

Заметим, что расслоение $(n-b(k)+2)\kappa_{\RP^m}$ является
обратным расслоением к касательному расслоению $T(\RP^{m})=b(k)
\kappa_{\RP^m}$ и, следовательно, является нормальным расслоением
к $\RP^m$ размерности $(n-b(k)+2)$. Рассмотрим  невырожденное
(стандартное) $s$--семейство сечений расслоения
$(n-b(k)+2)\kappa_{\RP^m}$.  Это семейство определено, поскольку
многообразие $\RP^m$ погружается в евклидово пространство
$\R^{2b(k)-k-2}$. Рассмотрим отображение $\kappa_N: N^{m-1} \to
\RP^{m}$ и индуцируем этим отображением невырожденное
$s$--семейство общего положения сечений
 расслоения $\bar
\nu_N$ и семейство проекций этого семейства, которое является семейством
сечений (вообще говоря, вырожденных) расслоения $\nu_N$:
\begin{eqnarray}\label{xi}
\{\bar \psi_1, \dots, \bar \psi_s\}, \qquad \{ \psi_1, \dots,
\psi_s\}_{N^{m-1} \times \{0\}}.
\end{eqnarray}

  Не ограничивая
общности, в классе нормального кобордизма многообразия $N^{m-1}$
выберем представитель, имеющий гомотопический $k$-тип пространства
$\RP^{\infty}$, индуцированный отображении $\kappa_N$.

Будем рассуждать по индукции c параметром $i=0, \dots, k$.  Пусть построена гомотопия
\begin{eqnarray}\label{chi}
\{\bar
\chi_1, \dots, \bar \chi_s\}_{N^{m-1} \times [0,i]}, \qquad \{
\chi_1, \dots, \chi_s\}_{N^{m-1} \times [0,i]},
\end{eqnarray}
с граничными условиями:
\begin{eqnarray}\label{chi0}
\{\bar \chi_1, \dots, \bar \chi_s\}_{N^{m-1} \times \{0\}}= \{\bar
\psi_1, \dots, \bar \psi_s\}_{N^{m-1} \times \{0\}},
\end{eqnarray}

\begin{eqnarray}\label{chii}
\{\bar \chi_1, \dots, \bar \chi_s\}_{N^{m-1} \times \{i\}}= \{\bar
\psi_1, \dots, \bar \psi_s\}_{N^{m-1} \times \{i\}}.
\end{eqnarray}

При этом выполнены нижеследующие условия 1-3.

--1. Подсемейство $\{\bar \psi_1, \dots, \bar
\psi_{s-k+i+1}\}_{N^{m-1} \times \{i\}}$ семейства $(\ref{chii})$
является невырожденным.

--2. В семействе проекций $\{ \psi_1, \dots, \psi_s\}_{N^{m-1}
\times \{i\}}$ семейства сечений $(\ref{chii})$ подсемейство $\{
\psi_1, \dots, \psi_{s-k+i}\}_{N^{m-1} \times \{i\}}$ является
невырожденным.

--3. Гомотопия $(\ref{chi})$ имеет, быть может, вырождение с тривиальным линейным расслоением ядер.
\[  \]

Обозначим через $L^i \subset N^{m-1} \times \{i\}$ -- $i$-мерное
подмногообразие базы, над которыми подсемейство $\{\bar \psi_1,
\dots, \bar \psi_{s-k+1+i}\}_{N^{m-1} \times \{i\}}$ пересекается
с подрасслоением $\kappa_N$ (т.е. семейство проекций $\{\xi_1,
\dots, \xi_{s-k+1+i}\}$  вырождается).

\subsubsection*{Базовый шаг индукции}

При $i=0$ фундаментальный класс многообразия $L^0$ представляет
гомологический эйлеров класс в группе $H_0(N^{m-1};\Z/2)$. Этот
гомологический класс двойственен в смысле Пуанкаре
когомологическому классу $w_{b(k)-2}((n-b(k)+1)\kappa_N =
w_1(\kappa_N)^{b(k)-2} \in H^{b(k)-2)}(N^{m-1};\Z/2)$ и полностью
определяется характеристическим числом
$$ \langle w_1^{m-1)}(\kappa_N);[N^{m-1}] \rangle = \langle w_1^{m}(\kappa_M);[M^m] \rangle. $$
Указанное характеристическое число равно нулю, поскольку
многообразие $N^{b(k)-2}$ ориентированно и имеет
$1$-гомотопический тип $\RP^{m}$. Тем самым, при $i=0$ в классе
изотопии невырожденного семейства сечений $\{\bar \psi_1, \dots,
\bar \psi_{s}\}_{N^{m-1} \times \{0\}}$ найдется семейство для
которого подсемейство проекций $\{ \psi_1, \dots, \psi_{s-k+1}
\}_{N^{m-1} \times \{0\}}$ не вырождается, что обеспечивает
основание индукции.

Для доказательства индуктивного предположим, что гомотопия $(\ref{chi})$ уже построена
при $i=i_0-1$ и построим такую же гомотопию при $i=i_0$.
Вначале построим
регулярную гомотопию
\begin{eqnarray}\label{chiii}
\{ \bar \chi_1, \dots, \bar \chi_{s-k+i_0+1} \}_{N^{m-1} \times [i_0,i_0+1]}
\end{eqnarray}
заданного семейства сечений $\{\bar \psi_1, \dots, \bar
 \psi_{s-k+i_0+1}\}_{N^{m-1} \times \{i_0\}}$ к семейству сечений
$\{\bar \psi_1, \dots, \bar
 \psi_{s-k+i_0+1}\}_{N^{m-1} \times \{i_0+1\}}$
 (условие регулярности
на компоненте границы $N^{m-1} \times \{i_0\}$ выполнено по
индуктивному предположению), для которой семейство проекций
\begin{eqnarray}\label{projchi}
\{ \psi_1, \dots, \psi_{s-k+i_0+1} \}_{N^{m-1} \times \{ i_0+1 \}}
\end{eqnarray}
не имеет вырождений.

Рассмотрим $i_0$-мерное подмногообразие $L^{i_0} \subset N^{m-1}
\times \{i_0\}$, состоящее из точек, в которых подсемейство
$\{\psi_1, \dots, \psi_{s-k+i_0+1}\}_{N^{m-1} \times \{i_0\}}$
вырождается. Из индуктивного предположения следует, что
ограничение класса $\kappa_N$ на подмногообразие $L^{i_0}$
тривиально, следовательно, подмногообразие $L^{i_0} \subset
N^{m-1} \times \{i_0\}$ является оснащенным. Обозначим оснащение
многообразия $L^{i_0}$ через $\Xi_L$.

Согласно основному результату [A-E]  существует $(i_0+1)$-мерное
скошенно-оснащенное в коразмерности $(b(k)-1)$ многообразие $(\Gamma^{i_0+1},\Psi_{\Gamma})$,
c краем $\partial
\Gamma^{i_0+1}$, диффеоморфным многообразию $L^{i_0}$, при этом ограничение скошенного оснащения
$\Psi_{\Gamma}$
на границу $L^{i_0}$ является оснащением, которое совпадает с заданным
оснащением $\Xi_L$ на крае $L^{i_0}$.

Поскольку $N^{m-1}$ имеет гомотопический
$k$-тип $\RP^{\infty}$ и нормальное расслоение $\nu_N$ снабжено канонической тривиализацией на любом
подмногообразии в $N^{m-1}$ размерности не более $k$, скошенно-оснащенное многообразие
$\Gamma^{i_0+1}$ можно рассматривать как
подмногообразие $\Gamma^{i_0+1} \subset N^{m-1} \times [i_0,i_0+1]$
cо стабильным скошенным оснащением $\Psi_{\Gamma}$,
 ограничение которого на $\partial
\Gamma^{i_0+1}$ совпадает с подмногообразием $L^{i_0} \subset N^{m-1}
\times \{i_0\}$ и при этом ограничение $\kappa_N$ на подмногообразие
$\Gamma^{i_0+1}$ совпадает с характеристическим классом скошенного
оснащения $\Xi_{\Gamma}$.

Регулярная гомотопия $(\ref{chiii})$ определена по поднятию проекции
\begin{eqnarray}\label{projchi1}
\{ \chi_1, \dots, \chi_{s-k+i_0+1} \}_ {N^{m-1} \times [i_0,
i_0+1] },
\end{eqnarray}
 для которой
 подмногообразие
$(\Gamma^{i_0+1},L^{i_0}) \subset (N^{m-1} \times [i_0,i_0+1],
N^{m-1} \times \{i_0\})$ служит подмногообразием особых сечений с
тривиальным линейным расслоением ядер.

Теперь построим гомотопию $(\ref{chi})$ при $i=i_0+1$. На отрезке
$[0,i_0]$ гомотопия уже построена по предположению индукции. На
отрезке $[i_0,i_0+1]$ на подсемействе первых $(s-k+i_0+1)$
векторов гомотопия совпадает с гомотопией $(\ref{chiii})$.
 Заметим, что особые слои гомотопии
$$\{ \bar \chi_1, \dots, \bar \chi_{s-k+i_0+2} \}_{N^{m-1} \times
[0,i_0]}$$ образуют многообразие с краем, которое обозначим через
\begin{eqnarray}\label{delta}
(\Delta^{i_0+1},K^{i_0}) \subset (N^{m-1} \times [0,i_0],N^{m-1}
\times \{i_0\}).
\end{eqnarray}
В силу условия 3, подмногообразие $(\ref{delta})$ является
оснащенным. Обозначим оснащение этого подмногообразия через
$\Psi_{\Delta}$. Граница $K^{i_0}$ многообразия $\Delta^{i_0+1}$
является подмногообразием особенностей семейства сечений
\begin{eqnarray}\label{psii0}
\{\bar \psi_1, \dots, \bar \psi_{s-k+i_0+2}\}_{N^{m-1} \times
\{i_0\}}.
\end{eqnarray}
 Воспользуемся тем, что многообразие особых сечений семейства
$\{\bar \psi_1, \dots \bar \psi_{s-k+i_0+1}\}_{N^{m-1} \times
\{i_0\}}$ заключено в шар $U \subset N^{m-1}$. Будем строить
гомотопию $(\ref{chi})$ на подсемействе первых $(s-k+i_0+2)$
векторов на многообразии $N^{m-1} \times [i_0,i_0+1]$, неподвижную
на границе $\partial U$ шара $U$, и совпадающую с гомотопией
$(\ref{chiii})$ на подсемействе первых $(s-k+i_0+1)$ векторов.
Указанную гомотопию построим так, что многообразие особенностей
семейства сечений $(\ref{psii0})$ перестраивается вдоль
оснащенного многообразия с краем
$$ (\Delta_1^{i_0+1},K^{i_0}) \subset (N^{m-1} \times
  [i_0,i_0+1],N^{m-1} \times \{i_0\}),
$$
 которое зеркально оснащенному подмногообразию
$(\ref{delta})$ и совпадает с этим многообразием на общей границе
$K^{i_0} \subset N^{m-1} \times \{i_0\}$.

Гомотопия $(\ref{chi})$ при $i=i_0+1$ на семействе первых
$(s-k+i_0+2)$ векторов построена. Продолжим построенную гомотопию
на все семейство $s$ векторов над многообразием $N^{m-1} \times
[i_0,i_0+1]$, обеспечив выполнение граничных условий на $N^{m-1}
\times \{i_0\}$ и выполнение условия 3.

Индуктивное предположение доказано. Следствие 40 в предположении
$j=0$ доказано.

\[  \]
\[  \]

Московская Область, Троицк, 142190, ИЗМИРАН.

 pmakhmet@izmiran.rssi.ru

\end{document}